\documentclass[aos,preprint,tikz,border=2pt,png]{imsart}

 \usepackage{amssymb, latexsym}
 \usepackage{amsmath} 
 \usepackage{amsthm}
 \usepackage{bm}
 \usepackage[mathscr]{eucal}
 \usepackage{enumerate}

\usepackage{amsmath,mathtools,amssymb,amsfonts, amsthm}

 \theoremstyle{plain}
\RequirePackage{amsthm,amsmath, amsfonts, amssymb, graphics, graphicx, subfigure, color, multirow}
\RequirePackage{natbib}
\usepackage{appendix}
\usepackage{mathtools}
\usepackage{graphicx}

\usepackage{algpseudocode}
\usepackage{algorithm}
\usepackage{comment}

\newtheorem{lem}{Lemma}
\newtheorem{cor}{Corollary}

\newtheorem{thm}{Theorem}
\newtheorem{prop}{Proposition}

\newtheorem{rmk}{Remark}

\newcommand\independent{\protect\mathpalette{\protect\independenT}{\perp}}
\def\independenT#1#2{\mathrel{\rlap{$#1#2$}\mkern2mu{#1#2}}}
\newcommand{\epf}{\hfill $\Box$}

\newcommand{\vM}{\boldsymbol{M}}

\newcommand{\vm}{\boldsymbol{m}}

\newcommand{\bbE}{\mathbb{E}}
\newcommand{\bbP}{\mathbb{P}}
\newcommand{\bbR}{\mathbb{R}}

\newcommand{\vA}{\boldsymbol{A}}

\newcommand{\vB}{\boldsymbol{B}}

\newcommand{\vC}{\boldsymbol{C}}
\newcommand{\vD}{\boldsymbol{D}}
\newcommand{\vE}{\boldsymbol{E}}
\newcommand{\vF}{\boldsymbol{F}}
\newcommand{\vG}{\boldsymbol{G}}

\newcommand{\vI}{\boldsymbol{I}}

\newcommand{\vV}{\boldsymbol{V}}

\newcommand{\vx}{\boldsymbol{x}}
\newcommand{\vX}{\boldsymbol{X}}
\newcommand{\vW}{\boldsymbol{W}}
\newcommand{\vw}{\boldsymbol{w}}
\newcommand{\vy}{\boldsymbol{y}}
\newcommand{\vY}{\boldsymbol{Y}}
\newcommand{\vz}{\boldsymbol{z}}
\newcommand{\vZ}{\boldsymbol{Z}}
\newcommand{\vT}{\boldsymbol{T}}

\newcommand{\balpha}{\boldsymbol{\alpha}}
\newcommand{\bbeta}{\boldsymbol{\beta}}
\newcommand{\bdelta}{\boldsymbol{\delta}}
\newcommand{\bgamma}{\boldsymbol{\gamma}}

\newcommand{\bLambda}{\boldsymbol{\Lambda}}

\newcommand{\bSigma}{\boldsymbol{\Sigma}}

\newcommand{\bepsilon}{\boldsymbol{\epsilon}}

\newcommand{\boldeta}{\boldsymbol{\eta}}
\newcommand{\bxi}{\boldsymbol{\xi}}

\begin{document}

 \begin{titlepage}

 \begin{center}

\textbf{Sparse Sliced Inverse Regression  Via Lasso}\\

 \vspace{10ex}

Qian Lin \footnote{Email: qianlin88@gmail.com},
Zhigen Zhao \footnote{Co-first author. Email: zhaozhg@temple.edu},
 and Jun S. Liu\footnote{Corresponding author. Email:jliu@stat.harvard.edu} \\

Center of Statistical Science, Tsinghua University, \\
Department of Statistical Science, Temple University\\ 
Department of Statistics, Harvard University \\

 \end{center}

{\bf Abstract:} For multiple index models, it has recently been shown that the sliced inverse regression (SIR) is consistent for estimating the sufficient dimension reduction (SDR) space if and only if $\rho=\lim\frac{p}{n}=0$, where $p$ is the dimension and $n$ is the sample size. Thus, when $p$ is of the same or a higher order of $n$, additional assumptions such as sparsity must be imposed in order to ensure consistency for SIR.  
By constructing artificial response variables made up from top eigenvectors of the estimated conditional covariance matrix, we introduce a simple Lasso regression method to obtain an estimate of the SDR space. The resulting algorithm, Lasso-SIR, is shown to be consistent and  achieve the optimal convergence rate under certain sparsity conditions when $p$ is of order $o(n^2\lambda^2)$, where $\lambda$ is the generalized signal-to-noise ratio. 
We also demonstrate the superior performance of Lasso-SIR compared with existing approaches  via extensive numerical studies and several real data examples.  

\end{titlepage}

\section{Introduction}\label{sec:intro}
 
 \setlength{\parindent}{5ex}

Dimension reduction and variable selection have become indispensable steps for modern-day data analysts  in dealing with the ``big data," where thousands or even millions of features  are often available for only hundreds or thousands of samples. 
With these ultra high-dimensional data, an effective modeling strategy is to assume that  only a few features and/or a few linear combinations of these features carry the information that researchers are interested in. One can consider the following {\it multiple index model} \citep{li1991sliced}:
\begin{align}\label{model:multiple}
 y=f(\bbeta_{1}^{\tau}\vx,\bbeta_{2}^{\tau}\vx,...,\bbeta_{d}^{\tau}\vx,\epsilon),
\end{align}
where $\vx$ follows a $p$-dimensional elliptical distribution with mean zero and covariance matrix $\bSigma$, the $\bbeta_i$'s are unknown projection vectors, $d$ is unknown but is assumed to be much smaller than $p$, and the error $\epsilon$ is independent of $\vx$ and has mean 0.
When $p$ is very large, it is reasonable to further restrict each $\bbeta_i$ to be  a sparse vector.

Since the introduction of the sliced inverse regression (SIR) method (\cite{li1991sliced}), many methods have been proposed to estimate the space spanned by $(\bbeta_1,\cdots,\bbeta_d)$ with few assumptions on the link function 
 $f(\cdot)$. 
 Assume the multiple index model (\ref{model:multiple}),
 the objective of all the SDR ( Sufficient Dimension Reduction, \cite{cook1998regressiongraph}) methods is to find the minimal subspace $\mathcal{S}\subseteq \mathbb{R}^p$ such that $y\independent \vx \mid P_{\mathcal{S}}\vx$, where $P_{\mathcal{S}}$ stands for the projection operator to the subspace $\mathcal{S}$.   When the dimension of $\vx$ is moderately large, all the SDR methods, including SIR, are proven to be successful \citep{xia2002adaptive, ni2005note, li2006sparse, li2007sparse, zhu2006sliced}.
However, these methods were previously known to work well when the sample size $n$ grows much faster than the dimension $p$, an assumption that becomes inappropriate for many modern-day datasets, such as those from biomedical researches. 
It is important to have a thorough investigation of ``the behavior of these SDR estimators when $n$ is not large relative to $p$", as raised by \cite{cook2012estimating}.

\cite{lin2015consistency} made an attempt to address the aforementioned challenge for SIR. They showed that, under mild conditions, the SIR estimate of the central space is consistent if and only if $\rho_n\stackrel{def}{=}p/n$ goes to zero as $n$ grows. Additionally, they showed that the convergence rate of the SIR estimate  of the central space (without any sparsity assumption) is $\rho_n$. When $p$ is greater than $n$,   certain constraints must be imposed in order for SIR to be consistent.  The sparsity assumption, i.e.,  
the number of active variables $s$ must be an order of magnitude smaller than $n$ and $p$, appears to be a reasonable one.
In a follow-up work, \cite{neykov2015support} studied the sign support recovery problem of the single index model ($d=1$), suggesting that the correct optimal convergence rate for estimating the central space might be $\frac{s\log(p)}{n}$, a speculation that is partially confirmed in \cite{lin2016minimax}.
It is shown that, for multiple index models with bounded dimension $d$ and the identity covariance matrix,
the optimal rate for estimating the central space is $\frac{ds+s\log(p/s)}{n\lambda}$, where $s$ is the number of active covariates and  $\lambda$ is the 
smallest
non-zero eigenvalue of $var(\bbE[\vx|y])$. 
 They further showed that the Diagonal-Thresholding algorithm  proposed in \cite{lin2015consistency} achieves the optimal rate for the single index model with the identity covariance matrix. 

 \paragraph{The main idea.} In this article, we introduce an efficient Lasso variant of SIR for the multiple index model (\ref{model:multiple}) with a general covariance matrix $\bSigma$.
Consider first the single index model: $y=f(\bbeta^{\tau}\vx,\epsilon)$. Let $\boldeta$ be the eigenvector associated with the largest eigenvalue of ${var}(\bbE[\vx|y])$. Since ${\bbeta\propto\bSigma^{-1}\boldeta}$, there are two immediate ways to estimate the space spanned by $\bbeta$. 
The first approach, as discussed in \cite{lin2015consistency}, estimates $\bSigma^{-1}$ and $\boldeta$ separately (see Algorithm \ref{alg:benchmarkI}).
The second one avoids a direct estimation of $\bSigma^{-1}$  by solving the following penalized least square problem: $\|\frac{1}{n}\vX\vX^{\tau}\bbeta-\boldeta\|_{2}^{2}+\mu\|\bbeta\|_{1}$, where $\vX$ is the $p\times n$ covariate matrix formed by the $n$ samples (see Algorithm \ref{alg:benchmarkII}).
However, similar to most $L_{1}$-penalization methods  for nonlinear models, theoretical underpinning of this approach has not been well understood.
Since these two approaches provide good estimates compared with earlier approaches (e.g.,\cite{li1991sliced, li2006sparse, li2007sparse}) as shown in \cite{lin2015consistency} and Supplementary Materials, we set the two approaches as benchmarks for comparisons.

We note that an eigenvector $\widehat{\boldeta}$ of $\widehat{var}(\bbE[\vx|y])$,  where $\widehat{var}(\bbE[\vx|y])$ is an estimate of the conditional covariance matrix $var(\bbE[\vx|y])$ using  SIR \citep{li1991sliced}, must be a linear combination of the column vectors of $\vX$. 
Thus, we can construct an artificial response vector $\widetilde{\vy}\in \bbR^{n}$ such that $\widehat{\boldeta}=\frac{1}{n}\vX\widetilde{\vy}$, 
and  estimate $\bbeta$ by 
solving another penalized least square problem:
$\frac{1}{2n}\|\widetilde{\vy}-\vX^{\tau}\bbeta\|^{2}_{2}+\mu\|\bbeta\|_{1}$ (see Algorithm \ref{alg:SIM}). We call this algorithm ``Lasso-SIR'', which is computationally very efficient. 
In Section \ref{sec:Theorems}, we further show that the convergence rate of the estimator resulting from Lasso-SIR  is ${\frac{s\log(p)}{n\lambda}},$ which is optimal  if $s=O(p^{1-\delta})$ for some positive constant $\delta$.
Note that Lasso-SIR  can be easily extended to other regularization and SDR methods, such as SCAD (\cite{fan2001variable}), Group Lasso (\cite{yuan2006model}), sparse Group Lasso (\cite{simon2013sparse}),  SAVE (\cite{dennis2000save}), etc. 

\paragraph{Connection to Other work}
Estimating the central space is widely considered as a generalized eigenvector problem in the literature \citep{li1991sliced,li2006sparse,li2007sparse,chen1998can}.
 \cite{lin2016minimax}  explicitly described the similarities and differences between SIR and  PCA (as first studied by \cite{jung2009pca}) under the ``high dimension, low sample size (HDLSS)" scenario. However, after comparing their results with those for Lasso regression, \cite{lin2016minimax}  advocated that a more appropriate prototype of SIR  (at least for the single index model) should be the linear regression. In the past three decades, tremendous efforts have been put into the study of  linear regression models $y=\vx^\tau\bbeta+\epsilon$ for HDLSS data.
By imposing the $L_{1}$ penalty on the regression  coefficients, the Lasso approach  \citep{tibshirani1996regression} produces a sparse estimator of $\bbeta$, which turns out to be rate optimal \citep{raskutti2011minimax}. 
Because of apparent limitations of linear models, there are many attempts to build flexible and computationally friendly semi-parametric models, such as the projection pursuit regression \citep{friedman1981projection, chen1991estimation}, sliced inverse regression \citep{li1991sliced}, MAVE \citep{xia2002adaptive}. 
However, none of these methods work under the HDLSS setting. Existing theoretical results for HDLSS data mainly focus on linear regressions \citep{raskutti2011minimax} and submatrix detections \citep{butucea2013detection}, and are not applicable to index models. In this paper, we provide a new framework for the theoretical investigation of  regularized SDR methods for HDLSS data.



The rest of the paper is organized as follows. After briefly reviewing SIR, we present the  Lasso-SIR algorithm in Section \ref{sec:sparseSIR}. The consistency of the Lasso-SIR estimate and its connection to the Lasso regression are presented in Section \ref{sec:Theorems}. 
Numerical simulations and real data applications are reported in Sections \ref{sec:numerical} and \ref{sec:realdata}. Some potential extensions are briefly discussed in Section \ref{sec:discussion}. To improve the readability, we defer all the proofs and brief reviews of some existing results to the appendix.

\section{Sparse SIR for High Dimensional Data}\label{sec:sparseSIR}

\noindent
{\bf\emph{Notations.}}
We adopt the following notations throughout this paper.
For a matrix $\vV$, we call the space generated by its column vectors the column space and denote it by $col(\vV)$. The $i$-th row and $j$-th column of the matrix are denoted by $\vV_{i,*}$ and $\vV_{*,j}$, respectively.
For (column) vectors $\vx$ and $\bbeta$ $\in$ $\mathbb{R}^{p}$, we denote their inner product $\langle \vx, \bbeta \rangle$ by $\vx(\bbeta)$, and the $k$-th entry of $\vx$ by $\vx(k)$. For two positive numbers $a,b$, we use $a\vee b$ and $a\wedge b$ to denote $\max\{a,b\}$ and  $\min\{a,b\}$ respectively;
We use $C$, $C'$, $C_1$ and $C_2$ to denote generic absolute constants, though the actual value may vary from case to case.
For two sequences $\{a_{n}\}$ and $\{b_{n}\}$, we denote $a_{n}\succ b_{n}$ and $a_{n}\prec b_{n}$ if there exist positive constants $C$ and $C'$ such that $a_{n} \geq Cb_{n}$ and $a_{n} \leq C'b_{n}$, respectively. We denote $a_{n}\asymp b_{n}$ if both $a_{n}\succ b_{n}$ and $a_{n}\prec b_{n}$ hold.
The $(1,\infty)$ norm and $(\infty, \infty)$ norm  of matrix $A$ are defined as
$
\|A\|_{1,\infty}=\max_{1\leq j \leq p} \sum_{i=1}^{p}|A_{i,j}|$  and $\max_{1\leq i, j \leq n}\|A_{i,j}\|$ respectively.
To simplify discussions, we assume that $\frac{s\log(p)}{n\lambda}$ is sufficiently small. We emphasize again that our covariate data $X$ is a $p\times n$ instead of the traditional $n\times p$ matrix. 

\noindent
{\bf\emph{A brief review of Sliced Inverse Regression (SIR).}} 
In the multiple index model (\ref{model:multiple}), the matrix $\vB$ formed by the vectors $\bbeta_{1},...,\bbeta_{d}$ is not identifiable. However, $col(\vB)$, the space spanned by the columns of $\vB$  is uniquely defined. 
Given $n$ $i.i.d.$ samples $(y_{i},\vx_{i})$, $i=1,\cdots,n$,   
 SIR \citep{li1991sliced} first divides the data into $H$ equal-sized slices according to the order statistics $y_{(i)}$, $i=1,\ldots, n$. To ease notations and arguments, we assume that $n=cH$ and $\bbE[\vx]=0$, and  re-express the data as $y_{h,j}$ and $\vx_{h,j}$, where  $h$ refers to the slice number and $j$ refers to the order number of a sample in the $h$-th slice, i.e.,
$y_{h,j}=y_{(c(h-1)+j)}, \ \vx_{h,j}=\vx_{(c(h-1)+j)}. $
Here $\vx_{(k)}$ is the concomitant of $y_{(k)}$.
Let the sample mean in the $h$-th slice 
be denoted by $ \overline{\vx}_{h,\cdot}$,
then  $\bLambda\triangleq var(\bbE[\vx|y])$ can be estimated by:
\begin{equation}\label{eqn:lambda}
\widehat{\bLambda}_{H}=\frac{1}{H}\sum_{h=1}^{H}\bar{\vx}_{h,\cdot}\bar{\vx}_{h,\cdot}^{\tau}=\frac{1}{H}\vX_{H}\vX^{\tau}_{H} 
 \end{equation}
where $\vX_{H}$ is a $p\times H$ matrix formed by the $H$ sample means, i.e., $\vX_H=(\bar{\vx}_{1,\cdot},\ldots, \bar{\vx}_{H,\cdot})$.
Thus, $col(\bLambda)$ is estimated by ${col(\widehat{\vV}_{H})}$, where ${\widehat{\vV}_{H}}$ is the matrix formed by the top $d$ eigenvectors of ${\widehat{\bLambda}_{H}}$.   The ${col(\widehat{\vV}_{H})}$ was shown to be a consistent estimator of $col(\bLambda)$ under a few technical conditions when $p$ is fixed \citep{duan1991slicing, hsing1992asymptotic, zhu2006sliced,li1991sliced, lin2015consistency}, which are summarized in the online supplementary file.  
Recently, \cite{lin2015consistency,lin2016minimax} showed that ${col(\widehat{\vV}_{H})}$  is consistent for $col(\bLambda)$ if and only if $\rho_n=\frac{p}{n}\rightarrow 0$ as $n\rightarrow \infty$, when the number of slices $H$ can be chosen as a fixed integer independent of $n$ and $p$ when the dimension $d$ of the central space is bounded.
When $\vx$'s distribution is elliptically symmetric, \cite{li1991sliced} showed that 
\begin{align}\label{eqn:beta-eta}
\bSigma col(\vB)=col(\bLambda),
\end{align}
and thus our goal is to recover $col(\vB)$ by solving the above equation. It is shown in \citep{lin2015consistency} that when $\rho_n\to 0$, $\widehat{col(\vB)} =\widehat{\bSigma}^{-1} col(\widehat{\vV}_{H})$ consistently estimate $col(\vB)$ where $\widehat{\bSigma}=\frac{1}{n}\vX\vX^{\tau}$ is the sample covariance matrix of $\vX$.
However, this simple approach breaks down when $\rho_n\not\rightarrow 0$, especially when $p\gg n$.
Although  stepwise methods \citep{zhong2012correlation,jiang2013sliced} can work under HDLSS settings,  the sparse SDR algorithms proposed in \cite{li2007sparse} and \cite{li2006sparse} appeared to be ineffective.  Below we  describe two intuitive non-stepwise methods for HDLSS scenarios, which will be used as  benchmarks in our simulation studies to measure the performance of newly proposed SDR algorithms.


\medskip

 {\bf Diagonal Thresholding-SIR.} 
When $p\gg n$, the Diagonal Thresholding (DT) screening method \citep{lin2015consistency} proceeds by  marginally screening all the variables via the diagonal elements of $\widehat{\bLambda}_{H}$ and then applying SIR  to those retained variables to obtain an estimate of $col(\vB)$.
The procedure is shown to be
 consistent if  the number of nonzero entries in each row of $\bSigma$ is bounded. 

\floatname{algorithm}{Algorithm}
\begin{algorithm}[H]
\caption{(DT-SIR)\label{alg:benchmarkI}}
\begin{algorithmic}[1]
\vspace*{1mm}
\\ Use the magnitudes of the diagonal elements of $\widehat{\bLambda}_{H}$ to select the set of important predictors $\mathcal{I}$, with $|\mathcal{I}|=o(n)$
\\ Apply SIR  to the data $(\vy, \vx_{\mathcal{I}})$ to estimate a subspace $\widehat{\mathcal{S}}_{\mathcal{I}}$.
\\ Extend $\widehat{\mathcal{S}}_{\mathcal{I}}$ to a subspace in $\bbR^{p}$ by filling in $0's$ for unimportant predictors. 
\end{algorithmic}
\end{algorithm}

\smallskip

{\bf Matrix Lasso.}
We can bypass the estimation and inversion of $\bSigma$ by solving an $L_1$ penalization problem. Since \eqref{eqn:beta-eta} holds at the population level, a  reasonable estimate of $col(\vB)$ can be obtained by solving a sample-version of the equation with an appropriate regularization term to cope with the high dimensionality. Let $\widehat{\boldeta}_{1},\cdots,\widehat{\boldeta}_{d}$ be the eigenvectors associated with the largest $d$ eigenvalues of $\widehat{\bLambda}_{H}$.  Replacing $\bSigma$ by its sample version $\frac{1}{n}\vX\vX^{\tau}$ and imposing an $L_1$ penalty, we obtain  a penalized sample version of  \eqref{eqn:beta-eta}:  
 \begin{align}\label{eqn:sample}
\|\frac{1}{n}\vX\vX^{\tau}\bbeta-\widehat{\boldeta}_{i}\|_{2}^{2}+\mu_{i}\|\bbeta\|_{1}
\end{align}
for some appropriate $\mu_{i}$'s.

\begin{algorithm}[H]
\caption{(Matrix Lasso)\label{alg:benchmarkII}
\vspace*{1mm}}
\begin{algorithmic}[1]
\\ Let $\widehat{\boldeta}_{1},...,\widehat{\boldeta}_{d}$ be the eienvectors associated with the largest $d$ eigenvalues of $\widehat{\bLambda}_{H}$;
\\ For $1\leq i \leq d$, let $\widehat{\bbeta}_{i}$ be the minimizer of equation \eqref{eqn:sample};
\\ Estimate the central space $col(\vB)$ by $col(\widehat{\bbeta}_{1},\ldots,\widehat{\bbeta}_{d})$.
\end{algorithmic}
\end{algorithm}
This simple procedure can be easily implemented  to produce sparse estimates of $\bbeta_{i}$'s. Empirically it works reasonably well, so we set it as another benchmark to compare with. Since we later observed that its numerical performance was consistently worse than that of  our main algorithm, Lasso-SIR, we did not further investigate its theoretical properties.

\medskip

\noindent 
{\bf\emph{The Lasso-SIR algorithm.}}
First consider the single index model
\begin{equation}\label{model:single}
y=f(\vx^{\tau}\bbeta_{0},\bepsilon).
\end{equation}
Without loss of generality,  we assume that $(\vx_i,y_i), \ i=1,\ldots, n,$  are
arranged in a way such that $y_1\le y_2\le \cdots\le y_n$. Construct an $n\times H$ matrix $\vM=\vI_{H}\otimes \mathbf{1}_{c}$, where $\mathbf{1}_{c}$ is the $c\times 1$ vector with all entries being 1. Then, according to the definition of $\vX_{H}$, we can write $\vX_{H}=\vX\vM/c$.
Let $\widehat{\lambda}$ be the largest eigenvalue of $\widehat{\bLambda}_{H}=\frac{1}{H}\vX_{H}\vX_{H}^{\tau}$ and let $\widehat{\boldeta}$ be the corresponding eigenvector of length 1. That is,
$$\hat{\lambda} \widehat{\boldeta} = \frac{1}{H}\vX_{H} \vX^{\tau}_{H}\widehat{\boldeta} = \frac{1}{n c}\vX\vM\vM^{\tau} \vX^{\tau} \widehat{\boldeta} .$$ 
Thus, by defining
\begin{align}\label{eqn:response:single}
\widetilde{\vy} = \frac{1}{c\widehat{\lambda}} \vM\vM^{\tau}\vX^{\tau}\widehat{\boldeta}
\end{align}
we have $\hat{\boldeta}= \frac{1}{n}\vX\widetilde{\vy}$.
Note that a key in estimating the central space $col(\bbeta)$ of SIR is the equation $\boldeta \propto \bSigma\bbeta$. 
If approximating $\boldeta$ and $\bSigma$ by $\widehat{\boldeta}$ and $\frac{1}{n}\vX\vX^{\tau}$ respectively, this equation can be written as $\frac{1}{n}\vX \widetilde{\vy} \propto \frac{1}{n}\vX\vX^{\tau}\bbeta$. To recover a sparse vector $\widehat{\bbeta}\propto \bbeta$, one can consider the following optimization problem
\[
\min||\bbeta||_{1}, \quad \textrm{subject to}\quad || \vX(\widetilde{\vy} - \vX^{\tau}\bbeta)||_{\infty}\leq \mu,
\]
which is known as the Dantzig selector \citep{candes2007dantzig}. A related formulation is the Lasso regression, where $\bbeta$ is estimated by the minimizer of 
\begin{align}\label{eqn:optimization:1}
\mathcal{L}_{\bbeta}=\frac{1}{2n}\|\widetilde{\vy}-\vX^{\tau}\bbeta\|_{2}^{2}+\mu\|\bbeta\|_{1}.
\end{align}
As shown by \cite{bickel2009simultaneous}, the Dantzig selector is asymptotically equivalent to the Lasso for linear regressions. We thus  propose and study the Lasso-SIR algorithm in this paper.

\floatname{algorithm}{Algorithm}

\begin{algorithm}[H]
\caption{(Lasso-SIR-1: for single index models)\label{alg:SIM}}
\begin{algorithmic}[1]
\vspace*{1mm}
\\ Let $\widehat{\lambda}$ and $\widehat{\boldeta}$ be the first eigenvalue and eigenvector of $\widehat{\bLambda}_{H}$, respectively;
\vspace*{1mm}
\\ Let $\widetilde{\vy}=\frac{1}{c\widehat{\lambda}} \vM\vM^{\tau}\vX^{\tau}\widehat{\boldeta}$ and solve the Lasso optimization problem
\begin{align*}
\widehat{\bbeta}(\mu)=\arg\min\mathcal{L}_{\bbeta}, \mbox{ where }\mathcal{L}_{\bbeta}=\frac{1}{2n}\|\widetilde{\vy}-\vX^{\tau}\bbeta\|_{2}^{2}+\mu\|\bbeta\|_{1},
\end{align*}
 where $\mu=C\sqrt{\frac{\log(p)}{n\widehat{\lambda}}}$ for sufficiently large constant $C$;
\\ Estimate $P_{\bbeta}$ by $P_{\widehat{\bbeta}(\mu)}$.
\end{algorithmic}
\end{algorithm}
There is no need to estimate the inverse of $\bSigma$ in Lasso-SIR. Moreover, since the optimization problem \eqref{eqn:optimization:1} is well studied for linear regression models \citep{tibshirani1996regression, efron2004least, friedman2010regularization}, we may formally ``transplant'' their results to the index models. Practically, we use the R package {\it glmnet} to solve the optimization problem where the tuning parameter $\mu$ is chosen using cross-validation. 

Last but not least, Lasso-SIR can be easily generalized to the multiple index model (\ref{model:multiple}). 
Let $\widehat{\lambda}_{i}, 1\leq i\leq d,$ be the $d$-top eigenvalues of $\widehat{\bLambda}_{H}$ and $\widehat{\boldeta}=(\widehat{\boldeta}_1, \cdots, \widehat{\boldeta}_d)$ be the corresponding eigenvectors. Similar to the definition of the {\it ``pseudo response variable''} for the single index model, we define a multivariate pseudo response $\widetilde{\vY}$ as
\begin{align}\label{eqn:response:multiple}
\widetilde{\vY}  = \frac{1}{c}\vM\vM^\tau \vX^{\tau}\widehat{\boldeta} \ diag(\frac{1}{\widehat{\lambda}_1},\cdots,\frac{1}{\widehat{\lambda}_d}).
\end{align}
We then apply the Lasso on each column of the pseudo response matrix to produce the corresponding estimate.

\begin{algorithm}[H]
\caption{(Lasso-SIR: for multiple index model)\label{alg:MIM1}}
\begin{algorithmic}[1]
\vspace*{1mm}
\\ Let $\widehat{\lambda}_{i}$ and $\widehat{\boldeta}_{i}$, $i=1,\cdots,d$ be the top $d$ eigenvalues and eigenvectors of $\widehat{\bLambda}_{H}$ respectively.  
\vspace*{1mm}
\\ Let $\widetilde{\vY}  = \frac{1}{c}\vM\vM^\tau \vX^{\tau}\widehat{\boldeta} \ diag(\frac{1}{\widehat{\lambda}_1},\cdots,\frac{1}{\widehat{\lambda}_d})$. For each $1\leq i \leq d$, solve the Lasso optimization problem
\begin{align*}
\widehat{\bbeta}_{i}=\arg\min\mathcal{L}_{\bbeta,i}\mbox{ where }\mathcal{L}_{\bbeta,i}=\frac{1}{2n}\|\widetilde{\vY}_{*,i}-\vX^{\tau}\bbeta\|_{2}^{2}+\mu_{i}\|\bbeta\|_{1},
\end{align*}
 where $\mu_{i}=C\sqrt{\frac{\log(p)}{n\widehat{\lambda}_{i}}}$ for sufficiently large constant $C$;
\\ Let $\widehat{\vB}$ be the matrix formed by $\widehat{\bbeta}_{1},\cdots,\widehat{\bbeta}_{d}$. The estimate of $P_{\vB}$ is given by $P_{\widehat{\vB}}$.
\end{algorithmic}
\end{algorithm}

  The number of directions $d$ plays an important role when implementing  Algorithm \ref{alg:MIM1}. A common practice is to locate the maximum gap among the ordered eigenvalues of the matrix $\widehat{\bLambda}_{H}$, which does not work well under HDLSS settings. In Section \ref{sec:Theorems}, we show that there exists a gap among the adjusted eigenvalues $\hat{\lambda}_{i}^a=\hat{\lambda}_i||\hat{\beta}_i||_2$ where $\hat{\beta}_i$ is the $i$-th output of Algorithm \ref{alg:MIM1}. Motivated by this, we estimate $d$ according to the following algorithm:

\begin{algorithm}[H]
\caption{Estimation of the number of directions $d$ \label{alg:choose:d}}
\begin{algorithmic}[1]
\vspace*{1mm}
\\ Apply Algorithm \ref{alg:MIM1} by setting $d=H$;
\vspace*{1mm}
\\ For each $i$, calculate $\hat{\lambda}_i^a = \hat{\lambda}_i ||\hat{\beta}_i||_2$;
\\ Apply the k-means method on $\hat{\lambda}_i^a$ with $k$ being 2 and the total number of points in the cluster with larger $\hat{\lambda}^a$ is the estimated value of $d$.
\end{algorithmic}
\end{algorithm}

\begin{rmk}
\normalfont In another paper that the authors are preparing, it is shown that the Lasso-SIR algorithm works on the joint distribution of $(\vX,Y)$ and is thus not tied to the single or multiple index models. We choose the single/multiple index models to have a clear representation of the central subspace $\mathcal{S}$, i.e., $\mathcal{S}=span\{\bbeta_{1},...,\bbeta_{d}\}$ .
\end{rmk}
\begin{rmk} \normalfont
When dealing with real data, we suggest that the users employ quantile normalization to transform each covariate when X is not normally distributed.
When $p$ is too large and beyond our bound of $n=O(\sqrt{p})$, as required by our provided R-package (see Section 7 for its downloading information), the user can first conduct variable screening based on DT-SIR, which is also included in this package.

\end{rmk}

\section{Consistency of Lasso-SIR}\label{sec:Theorems}
For simplicity, we assume that $\vx\sim N(0,\bSigma)$. The normality assumption can be relaxed to elliptically symmetric distributions with sub-Gaussian tail; however, this will make technical arguments unnecessarily tedious and is not the main focus of this paper. 
 From now on, we assume that $d$, the dimension of the central space, is bounded; thus we can  assume that $H$, the number of slices, is  a large enough but finite integer \citep{lin2016minimax, lin2015consistency}.
In order to prove the consistency, we need the following technical conditions:
\begin{itemize}
\item[${\bf A1)}$] There exist constants $C_{\min}$ and $C_{\max}$ such that $0<C_{\min}<\lambda_{\min}(\bSigma) \leq \lambda_{\max}(\bSigma)< C_{\max}$;
\end{itemize}
\begin{itemize}
\item
[${\bf A2)}$] There exists a constant $\kappa\geq 1$, such that 
$$0<\lambda=\lambda_{d}(var(\bbE[\vx|y])\leq ... \leq \lambda_{1}(var(\bbE[\vx|y])\leq \kappa\lambda\leq \lambda_{max}(\bSigma);$$
\end{itemize}

\begin{itemize}
\item
[${\bf A3)}$] The central curve $\vm(y)=\bbE[\vx|y]$ satisfies the sliced stability condition.
\end{itemize}

Condition A1 is commonly imposed in the analyses of high-dimensional linear regression models.
Condition A2 is merely a refinement of the coverage condition that is commonly imposed in the SIR literature, i.e., $rank(var(\bbE[\vx|y]))$=$d$.
For single index models, there is a more intuitive explanation of condition A2.
Since $rank(var(\bbE[\vx|y]))=1$, 
condition  A2 is simplified to
$0<\lambda=\lambda_1\leq \lambda_{max}(\bSigma)$ which is a direct corollary of the total variance decomposition identity ( i.e., $var(\vx)=var(\bbE[\vx|y])+\bbE[var(\vx|y)]$). 
We may treat $\lambda$ as a generalized $SNR$ and A2 simply requires that the generalized SNR is non-zero. 
Condition A3 is a property of the central curve, or equivalently, a regularity condition on the link function $f(\cdot)$ and the noise $\bepsilon$ introduced in \cite{lin2015consistency}. 


\begin{rmk}[Generalized $SNR$ and eigenvalue bound] 
\normalfont
Recall that the  signal-to-noise ratio ($SNR$) for the linear model $y=\bbeta^{\tau}\vx+\epsilon$, where $\vx \sim N(0,\bSigma)$  and $\epsilon\sim N(0,1)$,
is defined as
\begin{align*}
SNR=\frac{E[(\bbeta^{\tau}\vx)^{2}]}{\bbE[y^{2}]}=\frac{\|\bbeta\|_{2}^{2}\bbeta_{0}^{\tau}\bSigma\bbeta_{0}}{1+\|\bbeta\|_{2}^{2}\bbeta_{0}^{\tau}\bSigma\bbeta_{0}}.
\end{align*}
where $\bbeta_{0}=\bbeta/\|\bbeta\|_{2}$.
A simple calculation shows that
\begin{equation*}
var(\bbE[\vx|y])=\frac{\bSigma\bbeta\bbeta^{\tau}\bSigma}{\bbeta_{0}^{\tau}\bSigma\bbeta_{0}\|\bbeta\|_{2}^{2}+1}, \mbox{ \ \  and  \ \ } \lambda(var(\bbE[\vx|y]))=\frac{\bbeta_{0}^{\tau}\bSigma\bSigma\bbeta_{0}\|\bbeta\|_{2}^{2}}{\bbeta_{0}^{\tau}\bSigma\bbeta_{0}\|\bbeta\|_{2}^{2}+1},
\end{equation*}
where $\lambda(var(\bbE[\vx|y]))$ is the unique non-zero eigenvalue of $var(\bbE[\vx|y])$. This leads to the following identity for the linear model:
$$
\lambda (var(\bbE[\vx|y]))=\frac{\bbeta_{0}^{\tau}\bSigma\bSigma\bbeta_{0}}{\bbeta_{0}^{\tau}\bSigma\bbeta_{0}} SNR.
$$
Thus, in a  multiple index model we call $\lambda$, the smallest non-zero eigenvalue of $var(\bbE[\vx|y])$, the model's generalized $SNR$ . 
\end{rmk}

\medskip

\begin{thm}[Consistency of Lasso-SIR for Single Index Models]\label{thm:consistency} Assume that $n\lambda=p^{\alpha}$ for some $\alpha>1/2$ and  that conditions A1-A3 hold for  the single index model, $y=f(\bbeta^{\tau}_{0}\vx, \epsilon)$, where $\bbeta_{0}$ is a unit vector.
Let $\widehat{\bbeta}(\mu)$ be the output of Algorithm \ref{alg:SIM},
  then 
\begin{align*}
\|P_{\widehat{\bbeta}}-P_{\bbeta_{0}}\|_{F}\leq C_{1}\sqrt{ \frac{s\log(p)}{n\lambda} }
\end{align*}
holds with probability at least $1-C_{2}\exp(-C_{3}\log(p))$ for some constants $C_{2}$ and $C_{3}$. 
\end{thm}
 When no sparsity on $\boldeta$ is assumed, the condition $\alpha>1/2$ is necessary. This condition can be relaxed if a certain sparsity structure is imposed on $\boldeta$ or $\bSigma$ such that $\bSigma\bbeta$ becomes sparse. 
Next, we state the theoretical result regarding the multiple index model (\ref{model:multiple}).

\begin{thm}[Consistency of Lasso-SIR]\label{thm:consistency:multiple:version1}
 Assume that $n\lambda=p^{\alpha}$ for some $\alpha>1/2$, where $\lambda$ is the smallest nonzero eigenvalue of $var(\bbE[\vx|y])$, and that conditions A1-A3 hold for the multiple index model \eqref{model:multiple}.
 Assume further that the dimension $d$ of the central subspace  is known. Let $\widehat{B}$ be the output of Algorithm \ref{alg:MIM1},
then
\begin{align*}
\|P_{\widehat{\vB}}-P_{\vB}\|_{F}\leq C_{1}\sqrt{ \frac{s\log(p)}{n\lambda} }
\end{align*}
holds with probability at least $1-C_{2}\exp(-C_{3}\log(p))$ for some constants $C_{2}$ and $C_{3}$. 
\end{thm}

\cite{lin2016minimax} have shown that the lower bound of the risk $\bbE\|P_{\widehat{\vB}}-P_{\vB}\|^{2}_{F}$ is $\frac{s\log(p/s)}{n\lambda}$ when (i) $d=1$, or (ii) $d(>1)$ is finite and $\lambda>c_0>0$. 
This implies that if $s=O(p^{1-\delta})$ for some positive constant $\delta$, the Lasso-SIR algorithm achieves the optimal rate, i.e., we have the following corollary. 

\begin{cor}
Assume that conditions A1-A3 hold. If $n\lambda=p^{\alpha}$ for some $\alpha>1/2$ and $s=O(p^{1-\delta})$,  then Lasso-SIR estimate $P_{\widehat{\vB}}$ achieves the minimax rate when (i) $d=1$, or (ii) $d(>1)$ is finite and $\lambda>c_0>0$. 
\end{cor}

\begin{rmk} [ ] 
\normalfont
Consider the linear regression $y=\bbeta^\tau\vx+\epsilon$, where $\vx\sim N(0,\bSigma), \epsilon\sim N(0,1)$. 
It is shown in \cite{raskutti2011minimax} that the lower bound of the minimax rate of the $l_{2}$ distance between any estimator and the true $\bbeta$ is $\frac{s\log(p/s)}{n}$ and the convergence rate of Lasso estimator $\widehat{\bbeta}_{Lasso}$ is $\frac{s\log(p)}{n}$. Namely, the Lasso estimator is rate optimal for linear regression when $s=O(p^{1-\delta}) $ for some positive constant $\delta$.
 A simple calculation shows that $\lambda(var(\bbE[\vx|y]))\sim \|\bbeta\|^{2}_{2}$, if $\|\bbeta\|_{2}$ is bounded away from $\infty$. Consequently,
\begin{align}
\|P_{\widehat{\bbeta}_{Lasso}}-P_{\bbeta}\|_{F}
\leq 4\frac{\|\widehat{\bbeta}_{Lasso}-\bbeta\|_{2}}{\|\bbeta\|_{2}}\leq C\sqrt{\frac{s\log(p)}{n\lambda(var(\bbE[\vx|y]))}}
\end{align}
holds with high probability.
In other words, if we treat Lasso as a dimension reduction method (where $d=1$ and the link function is linear), the projection matrix $P_{\widehat{\bbeta}_{Lasso}}$ based on Lasso is rate optimal. The Lasso-SIR has extended the Lasso to the non-linear multiple index models. 
This justifies a statement in \cite{chen1998can}, stating that "SIR should be viewed as an alternative or generalization of the multiple linear regression".
The connection also justifies a speculation in \cite{lin2016minimax} that "a more appropriate prototype of the high dimensional SIR problem should be the sparse linear regression rather than the sparse PCA and the generalized eigenvector problem".
\end{rmk}

Determining the dimension $d$ of the central space is a challenging problem for SDR, especially for HDLSS cases. 
If we want to discern signals (i.e., the true directions) from noises (i.e., the other directions) simply  via the eigenvalues $\widehat{\lambda}_{i}$ of $\widehat{\bLambda}_{H}$, $i=1,...,H$, we face the problem that all these $\widehat{\lambda}_{i}$'s are of order $p/n$, but the gap between the signals and noises is of order $\lambda$ ($\leq C_{\max}$).
With the Lasso-SIR, we can bypass this difficulty by using the adjusted eigenvalues $\hat{\lambda}^{a}_{i}=\widehat{\lambda}_{i}\|\widehat{\bbeta}_{i}\|_{2}$, $i=1,\ldots,H$. To this end, we have the following theorem.
\begin{thm}\label{thm:3}
Let $\widehat{\bbeta}_{i}$ be the output of Algorithm \ref{alg:MIM1} for $i=1,\ldots,H$. 
Assume that $n\lambda=p^{\alpha}$ for some $\alpha>1/2$, $s\log(p)=o(n\lambda)$, and $H>d$, then, for some constants $C_{1},C_{2}$ $C_{3}$ and $C_{4}$, 
\begin{align*}
\hat{\lambda}_{i}^{a}
\geq C_{1}\sqrt{\lambda}-C_{2}\sqrt{\frac{s\log(p)}{n}} , \ \mbox{for } 1\leq i \leq d, \mbox{ and }  \\
\hat{\lambda}_{i}^{a}
\leq C_{3}\frac{\sqrt{p\log(p)}}{n\lambda}\sqrt{\lambda}+C_{4}\sqrt{\frac{s\log(p)}{n}}, \ \mbox{for } d+1\leq i\leq H,
\end{align*}
hold with probability at least $1-C_{5}\exp(-C_{6}\log(p))$ for some constants $C_{5}$ and $C_{6}$. 
\end{thm}
Theorem \ref{thm:3} states that, if  $s\log(p)\vee (p\log(p))^{1/2}=o(n\lambda)$, there is a clear gap between signals and  noise. The Lasso-SIR algorithm then provides us the rate optimal estimation of the central space. It can be easily verified that $p^{1/2}$ dominants $s\log(p)$ if $s<p^{1/2}$ and  $s\log(p)$ dominants $p^{1/2}$ if $s>p^{1/2}$.
The region $s^{2}=o(p)$  and the region $p=o(s^{2})$ are often referred to as the ``highly sparse" and ``moderately sparse" regions \citep{ingster2010detection}, respectively. These two scenarios should be treated  differently in high dimensional SIR and SDR frameworks, just like what has been done in high dimensional linear regression (\cite{ingster2010detection}).

\section{Simulation Studies}\label{sec:numerical}
\subsection{Single index models}\label{sec:nume:single}
Let $\bbeta$ be the vector of coefficients and let $\mathcal{S}$ be the active set; namely, $\beta_i=0, \forall i \in\mathcal{S}^c$. 
Furthermore, for each $i\in\mathcal{S}$, we simulated independently $\beta_i\sim N(0,1)$. 
Let $\vx$ be the design matrix with each row following $N(0, \bSigma)$. 
We consider two types of covariance matrices: (i) $\bSigma = (\sigma_{ij})$ where $\sigma_{ii}=1$ and $\sigma_{ij}=\rho^{|i-j|}$; 
and (ii) $\sigma_{ii}=1$, $\sigma_{i,j}=\rho$ when $i,j\in\mathcal{S}$ or $i,j\in \mathcal{S}^c$, and $\sigma_{i,j}=0.1$ when $i\in\mathcal{S}, j\in\mathcal{S}^c$ or vice versa. 
The first one represents a covariance matrix which is essentially sparse and we choose $\rho$ among 0, 0.3, 0.5, and 0.8. 
The second one represents a dense covariance matrix with $\rho$ chosen as 0.2. 
In all the simulations, we set $n=1,000$ and let $p$ vary among 100, 1,000, 2,000, and 4,000. 
For all the settings, the random error $\bepsilon$ follows $N(0,\vI_{n})$.
For single index models, we consider the following model settings:
\begin{itemize}
\item[I.] $\vy=  \vx\bbeta  + \bepsilon $ where $\mathcal{S}=\{1,2,\cdots,10\}$;
\item[II.] $\vy=  (\vx\bbeta)^3/2 + \bepsilon $ where $\mathcal{S}=\{1,2,\cdots,20\}$;
\item[III.]  $\vy=  \sin(\vx\bbeta)*exp(\vx\bbeta) + \bepsilon $ where $\mathcal{S}=\{1,2,\cdots,10\}$;
\item[IV.] $\vy= \exp( \vx\bbeta/10 ) + \bepsilon $ where $\mathcal{S}=\{1,2,\cdots,50\}$;
\item[V.] $\vy= \exp(\vx\bbeta  + \bepsilon ) $ where $\mathcal{S}=\{1,2,\cdots,7\}$.
\end{itemize}

\begin{table}
  \caption{ \label{tab:sim3} Estimation error for the first type covariance matrix with $\rho=0.5$.   }
  \begin{tabular}{|c|c|c|c|c|c|c|c|}
    \hline
    &  p & Lasso-SIR  & DT-SIR & Lasso & M-Lasso &Lasso-SIR(Known $d$) & $\hat{d}$ \\
    \hline
    \hline
    \multirow{4}{*}{I}          & 100  &  0.12 ( 0.02 )&  0.47 ( 0.11 )& 0.11 ( 0.02 )& 0.19 ( 0.08 )& 0.12 ( 0.02 )& 1 \\ & 1000  &  0.18 ( 0.02 )&  0.65 ( 0.14 )& 0.15 ( 0.02 )& 0.26 ( 0.02 )& 0.18 ( 0.02 )& 1 \\ & 2000  &  0.2 ( 0.02 )&  0.74 ( 0.15 )& 0.16 ( 0.02 )& 0.3 ( 0.03 )& 0.2 ( 0.02 )& 1 \\ & 4000  &  0.23 ( 0.09 )&  0.9 ( 0.17 )& 0.18 ( 0.01 )& 0.39 ( 0.09 )& 0.23 ( 0.03 )& 1 \\  
    \hline
    \multirow{4}{*}{II}    & 100  &  0.07 ( 0.01 )&  0.6 ( 0.1 )& 0.23 ( 0.03 )& 0.27 ( 0.31 )& 0.07 ( 0.01 )& 1 \\     & 1000  &  0.12 ( 0.02 )&  0.78 ( 0.11 )& 0.31 ( 0.04 )& 0.17 ( 0.02 )& 0.12 ( 0.02 )& 1 \\     & 2000  &  0.15 ( 0.02 )&  0.86 ( 0.13 )& 0.34 ( 0.05 )& 0.2 ( 0.03 )& 0.15 ( 0.02 )& 1 \\     & 4000  &  0.2 ( 0.04 )&  0.99 ( 0.15 )& 0.37 ( 0.05 )& 0.28 ( 0.06 )& 0.19 ( 0.03 )& 1 \\ 
    \hline 
    \multirow{4}{*}{III}  & 100  &  0.21 ( 0.03 )&  0.55 ( 0.12 )& 1.25 ( 0.19 )& 0.26 ( 0.11 )& 0.21 ( 0.03 )& 1 \\     & 1000  &  0.28 ( 0.04 )&  0.74 ( 0.14 )& 1.32 ( 0.18 )& 0.51 ( 0.04 )& 0.27 ( 0.04 )& 1 \\     & 2000  &  0.35 ( 0.17 )&  0.87 ( 0.17 )& 1.34 ( 0.14 )& 0.66 ( 0.14 )& 0.31 ( 0.05 )& 1.1 \\     & 4000  &  0.46 ( 0.28 )&  1 ( 0.25 )& 1.33 ( 0.16 )& 0.83 ( 0.22 )& 0.39 ( 0.1 )& 1.1 \\   
    \hline
    \multirow{4}{*}{IV}    & 100  &  0.46 ( 0.05 )&  0.92 ( 0.09 )& 0.78 ( 0.12 )& 0.58 ( 0.06 )& 0.45 ( 0.04 )& 1 \\     & 1000  &  0.62 ( 0.22 )&  1.07 ( 0.18 )& 0.87 ( 0.11 )& 0.78 ( 0.22 )& 0.59 ( 0.04 )& 1.1 \\     & 2000  &  0.71 ( 0.34 )&  1.22 ( 0.26 )& 0.89 ( 0.12 )& 0.94 ( 0.31 )& 0.59 ( 0.04 )& 1.3 \\     & 4000  &  0.71 ( 0.26 )&  1.3 ( 0.18 )& 0.91 ( 0.13 )& 1 ( 0.22 )& 0.63 ( 0.04 )& 1.2 \\  
    \hline
    \multirow{4}{*}{V}      & 100  &  0.12 ( 0.02 )&  0.37 ( 0.1 )& 0.42 ( 0.18 )& 0.15 ( 0.02 )& 0.12 ( 0.02 )& 1 \\     & 1000  &  0.2 ( 0.03 )&  0.55 ( 0.15 )& 0.55 ( 0.22 )& 0.41 ( 0.05 )& 0.2 ( 0.05 )& 1 \\     & 2000  &  0.38 ( 0.34 )&  0.8 ( 0.29 )& 0.6 ( 0.24 )& 0.67 ( 0.27 )& 0.29 ( 0.18 )& 1.2 \\     & 4000  &  0.78 ( 0.51 )&  1.22 ( 0.31 )& 0.77 ( 0.25 )& 1.06 ( 0.41 )& 0.48 ( 0.31 )& 1.5 \\  
    \hline
  \end{tabular}
\end{table}

The goal is to estimate $col(\bbeta)$, the space spanned by $\bbeta$. As in \cite{lin2015consistency}, the estimation error is defined as  $\mathcal{D}(\widehat{col(\bbeta)}, col(\bbeta))$, where $\mathcal{D}(M,N)$, the distance between two subspaces $M, N\subset \mathcal{R}^p$, is defined as  the Frobenius norm of $P_{M}-P_{N}$ where $P_{M}$ and $P_{N}$ are the projection matrices associated with these two spaces. The methods we compared with are DT-SIR, matrix Lasso (M-Lasso), and Lasso. 
The number of slices $H$ is chosen as 20 in all simulation studies. 
The number of directions $d$ is chosen according to Algorithm \ref{alg:choose:d}. 
Note that both benchmarks (i.e., DT-SIR and M-Lasso) require  the knowledge of $d$ as well. 
To be fair, we use the $\widehat{d}$ estimated based on Algorithm \ref{alg:choose:d} for both benchmarks. For comparison, we have also included the estimation error of Lasso-SIR assuming $d$ is known. For each $p$, $n$, and $\rho$, we replicate the above steps 100 times to calculate the average estimation error for each setting. We tabulated the results for the first type of covariance matrix with $\rho=0.5$ in Table \ref{tab:sim3}  and put the results for other settings in Tables \ref{tab:sim1}-\ref{tab:sim5} in the online supplementary file. 
The average of estimated directions $\hat{d}$ is reported in the last column of these tables.

The simulation results in Table~\ref{tab:sim3} show that Lasso-SIR  outperformed both DT-SIR and M-Lasso under all settings. The performance of DT-SIR has become worse when the dependence is stronger and denser. The reason is that this method is based on the diagonal threshold and is only supposed to work well for the diagonal covariance matrix. Overall, Algorithm \ref{alg:choose:d} provided a reasonable estimate of $d$ especially for moderate covariance matrix. When assuming $d$ is known, the performances of both DT-SIR and M-Lasso are inferior to Lasso-SIR, and are thus not reported.

 Under Setting I when the true model is linear, Lasso performed the best among all the methods, as expected. However, the difference between Lasso and Lasso-SIR is not significant, implying that Lasso-SIR does not sacrifice much efficiency without the knowledge of the underlying linearity. On the other hand, when the models are not linear (Case II-VI), Lasso-SIR worked much better than Lasso. We observed that  Lasso performed better than Lasso-SIR for Setting V when $\rho$=0.8 (Supplemental Materials) or when the covariance matrix is dense. One explanation is that Lasso-SIR tends to overestimate $d$ under these conditions while Lasso used the actual $d$. If assuming known $d=1$, Lasso-SIR's  estimation error is smaller than that of Lasso.

The results, reported in the supplementary material, for the other values of $\rho$ are similar to what we observed when $\rho=0.5$. The Lasso-SIR performed the best when compared to its competitors.


\subsection{Multiple index models}\label{sec:nume:multiple}
Let $\bbeta$ be the $p\times 2$ matrix of coefficients and $\mathcal{S}$ be the active set. Let $\vx$ be simulated similarly as in Section \ref{sec:nume:single}, and denote $\vz=\vx\bbeta$. Consider the following settings:

\begin{enumerate}
\item [VI.] $y_i= | z_{i2}/4+2|^3 * sgn( z_{i1})  +\epsilon_i$ where $\mathcal{S}=\{1,2,\cdots, 7\}$ and $\beta_{1:4,1}=1, \beta_{5:7,2}=1$, and $\beta_{i,j}=0$ otherwise;
\item [VII.] $y_i = z_{i1} *exp( z_{i2}) + \epsilon_i$ where $\mathcal{S}=\{1,2,\cdots, 12\}$ and $\beta_{1:7,1}, \beta_{8:12,2}\sim N(0,1)$, and $\beta_{i,j}=0$ otherwise;
\item [VIII.] $y_i=z_{i1}*exp(z_{i2}+\epsilon_i)$ where $\mathcal{S}=\{1,2,\cdots, 12\}$ and $\beta_{1:7,1}, \beta_{8:12,2}\sim N(0,1)$, and $\beta_{i,j}=0$ otherwise;
\item [IX.] $y_i=z_{i1}*( 2+z_{i2}/3)^2 + \epsilon_i$ where $\mathcal{S}=\{1,2,\cdots,12\}$ and $\beta_{1:8,1}=1, \beta_{9:12,2}=1$ and $\beta_{i,j}=0$ otherwise.
\end{enumerate}

\begin{table}
  \caption{\label{tab:sim8} Estimation error for the first type covariance matrix with $\rho=0.5$.   }
  \begin{tabular}{|c|c|c|c|c|c|c|}
    \hline
    &  p & Lasso-SIR  & DT-SIR & M-Lasso &Lasso-SIR(Known $d$) &$\hat{d}$ \\
    \hline
    \hline
    \multirow{4}{*}{VI}    & 100  &  0.26 ( 0.06 )&  0.57 ( 0.15 )& 0.31 ( 0.05 )& 0.26 ( 0.05 )& 2 \\      & 1000  &  0.33 ( 0.07 )&  0.74 ( 0.17 )& 0.62 ( 0.04 )& 0.33 ( 0.07 )& 2 \\      & 2000  &  0.36 ( 0.11 )&  0.92 ( 0.18 )& 0.73 ( 0.07 )& 0.38 ( 0.08 )& 2 \\      & 4000  &  0.44 ( 0.14 )&  1.12 ( 0.25 )& 0.87 ( 0.1 )& 0.42 ( 0.09 )& 2 \\     
    \hline
    \multirow{4}{*}{VII} &   100  &  0.32 ( 0.04 )&  0.67 ( 0.11 )& 0.42 ( 0.04 )& 0.32 ( 0.04 )& 2 \\      & 1000  &  0.6 ( 0.28 )&  0.93 ( 0.22 )& 1.02 ( 0.2 )& 0.66 ( 0.3 )& 2.1 \\      & 2000  &  0.95 ( 0.44 )&  1.18 ( 0.27 )& 1.35 ( 0.32 )& 0.83 ( 0.35 )& 2.3 \\      & 4000  &  1.17 ( 0.38 )&  1.43 ( 0.31 )& 1.47 ( 0.33 )& 1.08 ( 0.33 )& 2.1 \\  
    \hline
    \multirow{4}{*}{VIII}    & 100  &  0.29 ( 0.09 )&  0.61 ( 0.11 )& 0.34 ( 0.08 )& 0.25 ( 0.03 )& 2 \\      & 1000  &  0.37 ( 0.08 )&  0.82 ( 0.14 )& 0.69 ( 0.13 )& 0.35 ( 0.07 )& 2 \\      & 2000  &  0.54 ( 0.35 )&  1 ( 0.25 )& 0.92 ( 0.28 )& 0.47 ( 0.22 )& 2.2 \\      & 4000  &  0.88 ( 0.45 )&  1.37 ( 0.26 )& 1.27 ( 0.31 )& 0.71 ( 0.37 )& 2.5 \\  
        \hline
    \multirow{4}{*}{IX}    & 100  &  0.43 ( 0.06 )&  0.74 ( 0.12 )& 0.48 ( 0.05 )& 0.43 ( 0.07 )& 2 \\      & 1000  &  0.47 ( 0.09 )&  0.91 ( 0.15 )& 0.91 ( 0.05 )& 0.48 ( 0.09 )& 2 \\      & 2000  &  0.58 ( 0.23 )&  1.11 ( 0.23 )& 1.12 ( 0.16 )& 0.5 ( 0.1 )& 2.1 \\      & 4000  &  0.57 ( 0.18 )&  1.25 ( 0.22 )& 1.23 ( 0.1 )& 0.56 ( 0.11 )& 2 \\  
    \hline
    \hline
  \end{tabular}
\end{table}

For the multiple index models, we compared both benchmarks (DT-SIR and M-Lasso) with Lasso-SIR. Lasso is not applicable for these cases and is thus not included.
Similar to Section \ref{sec:nume:single}, we tabulated the results for the first type covariance matrix with $\rho=0.5$ in Table \ref{tab:sim8} and put the results for others in Tables \ref{tab:sim6}-\ref{tab:sim10} in the online supplementary file. 

For the identity covariance matrix ($\rho=0$), there was little difference between performances of Lasso-SIR and DT-SIR. However, Lasso-SIR was substantially better than DT-SIR in other cases. Under all settings, Lasso-SIR worked much better than the matrix Lasso. For the dense covariance matrix $\bSigma_2$, Algorithm \ref{alg:choose:d} tended to underestimate $d$, which is worthy of further investigation.

The results, reported in the supplementary material, for the other values of $\rho$ are similar to what we observed when $\rho=0.5$. The Lasso-SIR performs the best when compared to its competitors.

There are other sparse inverse regression method, such as the Sparse SIR, given in \cite{li2006sparse}. In \cite{lin2015consistency}, we have shown that the DT-SIR outperforms this method. We thus did not include the numerical comparison. For the reason of completeness, we have included the numerical results of comparing Lasso-SIR and  Sparse SIR in Section D of the online supplementary file, showing that Lasso-SIR is better than Sparse-SIR.


\subsection{Discrete responses}
We consider the following simulation settings where for the response variable $Y$ is discrete. 

\begin{enumerate}
\item[X.] $\vy= 1( \vx\bbeta  + \bepsilon >0 ) $ where $\mathcal{S}=\{1,2,\cdots,10\}$;
\item[XI.] $\vy= 1(exp(\vx\bbeta) + \epsilon >0 ) $ where $\mathcal{S}=\{1,2,\cdots,7\}$;
\item[XII.] $\vy= 1( (\vx\bbeta)^3/2 + \epsilon) $ where $\mathcal{S}=\{1,2,\cdots,20\}$;
\item[XIII.] Let $\vz=\vx\bbeta$ where $\mathcal{S}=\{1,2,\cdots, 12\}$, $\bbeta$ is a $p$ by $2$ matrix with  $\beta_{1:7,1}, \beta_{8:12,2}\sim N(0,1)$ and $\beta_{i,j}=0$ otherwise. The response $y_i$ is
\begin{eqnarray*}
y_i = \left\{ \begin{array}{cc} 1, & \textrm{ if $z_{i1}+\epsilon_{i1}<0 $},\\
2, & \textrm{ if $z_{i1}+\epsilon_{i1}>0 $ and $z_{i2}+\epsilon_{i2} <0 $},\\
3, & \textrm{ if $z_{i1}+\epsilon_{i1}>0 $ and $z_{i2}+\epsilon_{i2}>0$},
\end{array}\right.
\end{eqnarray*}
where $\epsilon_{ij}\sim N(0,1)$.
\end{enumerate}

In settings X, XI, and XII, the response variable is dichotomous, and $\beta_i\sim N(0,1)$ when $i\in\mathcal{S}$ and $\beta_i=0$ otherwise. Thus the number of slices $H$ can only be  2. For Setting XIII where the response variable is trichotomous, the number of slices $H$ is chosen as 3. The number of direction $d$ is chosen as $H-1$ in all these simulations.

\begin{table}
  \caption{\label{tab:sim13} Estimation error for the first type covariance matrix with $\rho=0.5$.   }
  \begin{tabular}{|c|c|c|c|c|c|}
    \hline
    &  p & Lasso-SIR  & DT-SIR & M-Lasso &Lasso \\
    \hline
    \hline
    \multirow{4}{*}{X}     & 100  &  0.22 ( 0.03 ) &  0.66 ( 0.05 )& 0.26 ( 0.03 )& 0.2 ( 0.03 ) \\          & 1000  &  0.26 ( 0.04 ) &  1.21 ( 0.03 )& 0.52 ( 0.03 )& 0.28 ( 0.03 ) \\          & 2000  &  0.27 ( 0.03 ) &  1.33 ( 0.02 )& 0.59 ( 0.02 )& 0.29 ( 0.04 ) \\          & 4000  &  0.28 ( 0.04 ) &  1.39 ( 0.02 )& 0.65 ( 0.03 )& 0.3 ( 0.04 ) \\   
    \hline
    \multirow{4}{*}{XI}    & 100  &  0.32 ( 0.07 ) &  0.83 ( 0.07 )& 0.6 ( 0.17 )& 0.33 ( 0.07 ) \\          & 1000  &  0.43 ( 0.1 ) &  1.32 ( 0.02 )& 1.07 ( 0.05 )& 0.45 ( 0.09 ) \\          & 2000  &  0.45 ( 0.09 ) &  1.38 ( 0.01 )& 1.15 ( 0.04 )& 0.46 ( 0.09 ) \\          & 4000  &  0.49 ( 0.12 ) &  1.41 ( 0.01 )& 1.2 ( 0.05 )& 0.51 ( 0.12 ) \\  
    \hline
    \multirow{4}{*}{XII}     & 100  &  0.24 ( 0.03 ) &  0.63 ( 0.05 )& 0.52 ( 0.35 )& 0.22 ( 0.03 ) \\          & 1000  &  0.33 ( 0.03 ) &  1.18 ( 0.04 )& 0.53 ( 0.03 )& 0.32 ( 0.03 ) \\          & 2000  &  0.37 ( 0.05 ) &  1.3 ( 0.04 )& 0.62 ( 0.03 )& 0.35 ( 0.03 ) \\          & 4000  &  0.4 ( 0.04 ) &  1.38 ( 0.03 )& 0.68 ( 0.03 )& 0.39 ( 0.04 ) \\    
    
    \hline
    \multirow{4}{*}{XIII}   & 100  &  0.38 ( 0.06 ) &  1.09 ( 0.06 )& 0.61 ( 0.05 )& 1.07 ( 0.02 ) \\          & 1000  &  0.39 ( 0.07 ) &  1.79 ( 0.02 )& 1.12 ( 0.05 )& 1.08 ( 0.02 ) \\          & 2000  &  0.38 ( 0.07 ) &  1.91 ( 0.02 )& 1.24 ( 0.04 )& 1.09 ( 0.03 ) \\          & 4000  &  0.42 ( 0.07 ) &  1.98 ( 0.01 )& 1.32 ( 0.03 )& 1.1 ( 0.03 ) \\     
    \hline
    \hline
  \end{tabular}
\end{table}

Similar to the previous two sections, we calculated the average estimation errors for Lasso-SIR (Algorithm \ref{alg:MIM1}), DT-SIR, M-Lasso, and generalized-Lasso based on 100 replications and reported the result in Table \ref{tab:sim13} for the first type covariance matrix with $\rho=0.5$ and the results for other cases in Tables \ref{tab:sim11}-\ref{tab:sim15} in online supplementary file.
It is clearly seen that Lasso-SIR performed much better than DT-SIR and M-Lasso under all settings and the improvements were very significant. The generalized Lasso performed as good as Lasso-SIR for the dichotomous response; however, it performed substantially worse for Setting XIII.

\section{Applications to Real Data}\label{sec:realdata}

\noindent
{\bf\emph{Arcene Data Set.}} 
We first apply the methods to a two-class classification problem, which aims to distinguish between cancer patients and normal subjects from using their mass-spectrometric measurements. The data were obtained by the National Cancer Institute (NCI) and the Eastern Virginia Medical School (EVMS) using the SELDI technique, including samples from 44 patients with ovarian and prostate cancers and 56 normal controls. The dataset was downloaded from the UCI machine learning repository (\cite{Lichman:2013}), where a detailed description can be found. It has also been used in the NIPS 2003 feature selection challenge (\cite{guyon2004result}).  
For each subject, there are 10,000 features where 7,000 of them are real variables and 3,000 of them are random probes. There are 100 subjects in the validation set.

After standardizing $\vX$, we estimated the number of directions $d$ as 1 using Algorithm \ref{alg:choose:d}. We then applied Algorithm \ref{alg:SIM} and the sparse PCA to calculate the direction of $\bbeta$ and the corresponding components, followed by a logistic regression model. We  applied the fitted model to the validation set and calculated the probability of each subject being a cancer patient. We also fitted a Lasso logistic regression model to the training set and applied it to the validation set to calculate the corresponding probabilities. 

\begin{figure}
  \centering
\includegraphics[height=50mm,width=50mm]{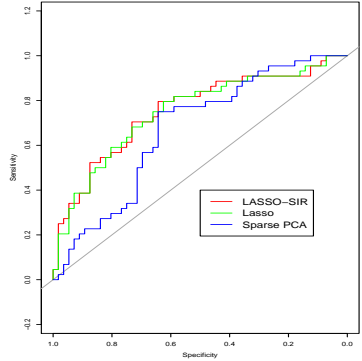}
\caption{ROC curve of various methods for Arcene Data set.
}\label{fig:arcene}
\end{figure}

In Figure \ref{fig:arcene}, we plot the Receiver Operating Characteristic (ROC) curves for various methods. Lasso-SIR, represented by the red curve, was slightly better than Lasso (insignificant) and the sparse PCA, represented by the green and blue curves respectively. The areas under these three curves are 0.754, 0.742, and 0.671, respectively. 

\noindent
{\bf\emph{HapMap.}}
In this section, we analyzed a data set with a continuous response. We consider the gene expression data from 45 Japanese and 45 Chinese from the international ``HapMap'' project (\cite{thorisson2005international, thorgeirsson2010sequence}). The total number of probes is 47,293. According to \cite{thorgeirsson2010sequence}, the gene {\it CHRNA6} is the subject of many nicotine addiction studies. Similar to \cite{fan2015discoveries}, we treat the mRNA  expression of {\it CHRNA6} as the response $Y$ and expressions of other genes as the covariates. Consequently, the number of dimension $p$ is 47,292, much greater than  the number of subjects $n$=90.

We first applied Lasso-SIR to the data set with $d$ being chosen as 1 according to Algorithm \ref{alg:choose:d}. The number of selected variables was 13. Based on the estimated coefficients $\bbeta$ and $\vX$, we calculated the first component and the scatter plot between the response $Y$ and this component, showing a moderate linear relationship between them. We then fitted a linear regression between them. The R-sq of this model is 0.5596 and the mean squared error of the fitted model 0.045.

We also applied Lasso to estimate the direction $\bbeta$. The tuning parameter $\lambda$ is chosen as 0.1215 such that the number of selected variables is also 13. When fitting a regression model between $Y$ and the component based on the estimated $\bbeta$, the R-sq is 0.5782 and the mean squared error is 0.044. There is no significant difference between these two approaches. This confirms the message that Lasso-SIR performs as good as Lasso when the linearity assumption is appropriate.

We have also calculated a direction and the corresponding components based on the sparse PCA \citep{zou2006sparse}. We then fitted a regression model. The R-sq is only 0.1013 and the mean squared error is 0.093, significantly worse than the above two approaches.

\medskip

\noindent
{\bf\emph{Classify Wine Cultivars.}}
We investigate the popular wine data set which has been used to compare various classification methods. This is a three-class classification problem. The data, available from the UCI machine learning repository (\cite{Lichman:2013}),  consists of 178 wines grown in the same region in Italy under three different cultivars. For each wine, the chemical analysis was conducted and the quantities of 13 constituents are obtained, which are Alcohol, Malic acid, Ash, Alkalinity of ash, Magnesium, Total Phenols, Flavanoids, Nonflavanoid Phenols, Proanthocyanins, Color intensity, Hue, OD280/OD315 of diluted wines, and Proline. One of the goals is to use these 13 features to classify the cultivar. 

\begin{figure}
  \centering
  \includegraphics[width=60mm, height=60mm]{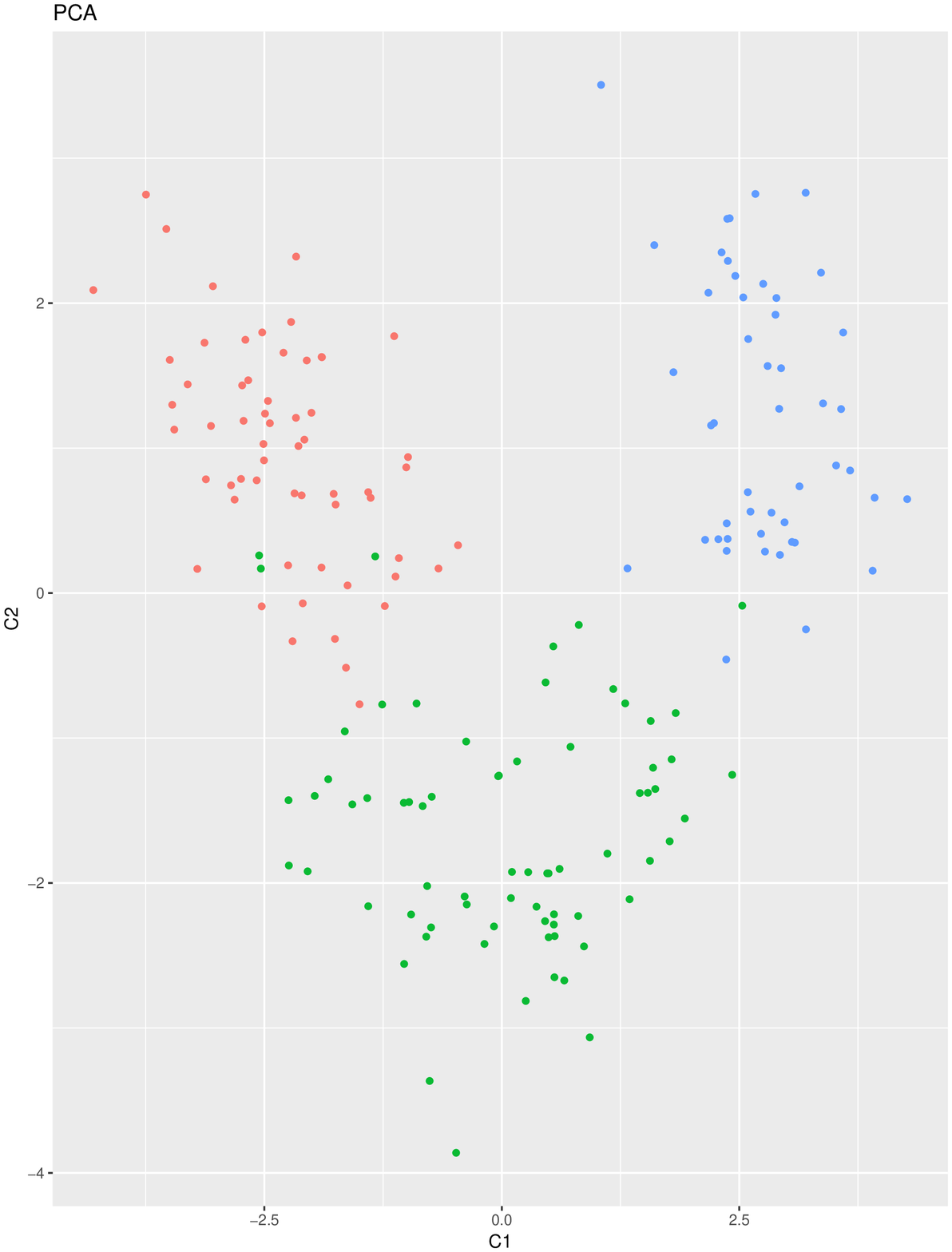}
  \includegraphics[width=60mm, height=60mm]{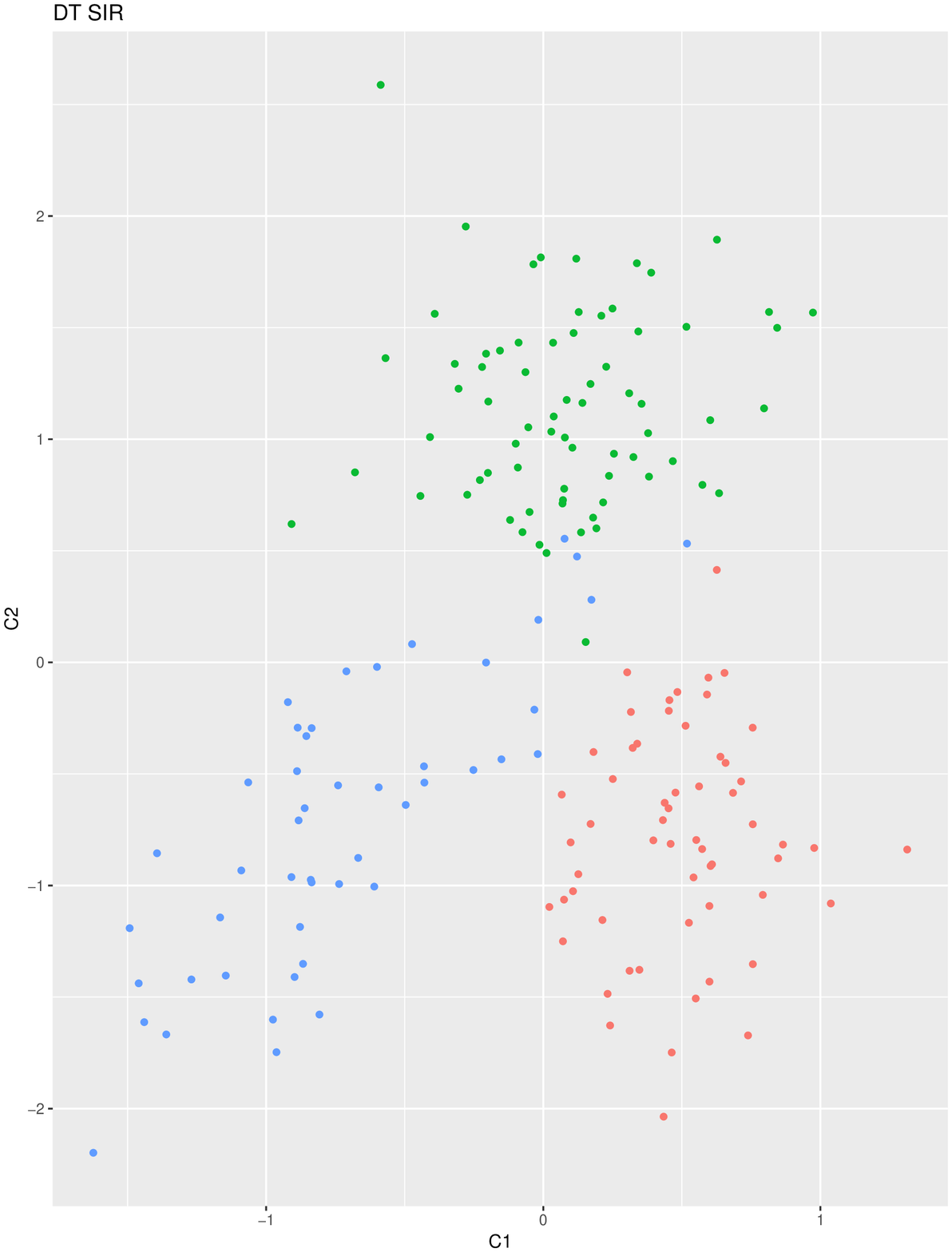}
  \includegraphics[width=60mm, height=60mm]{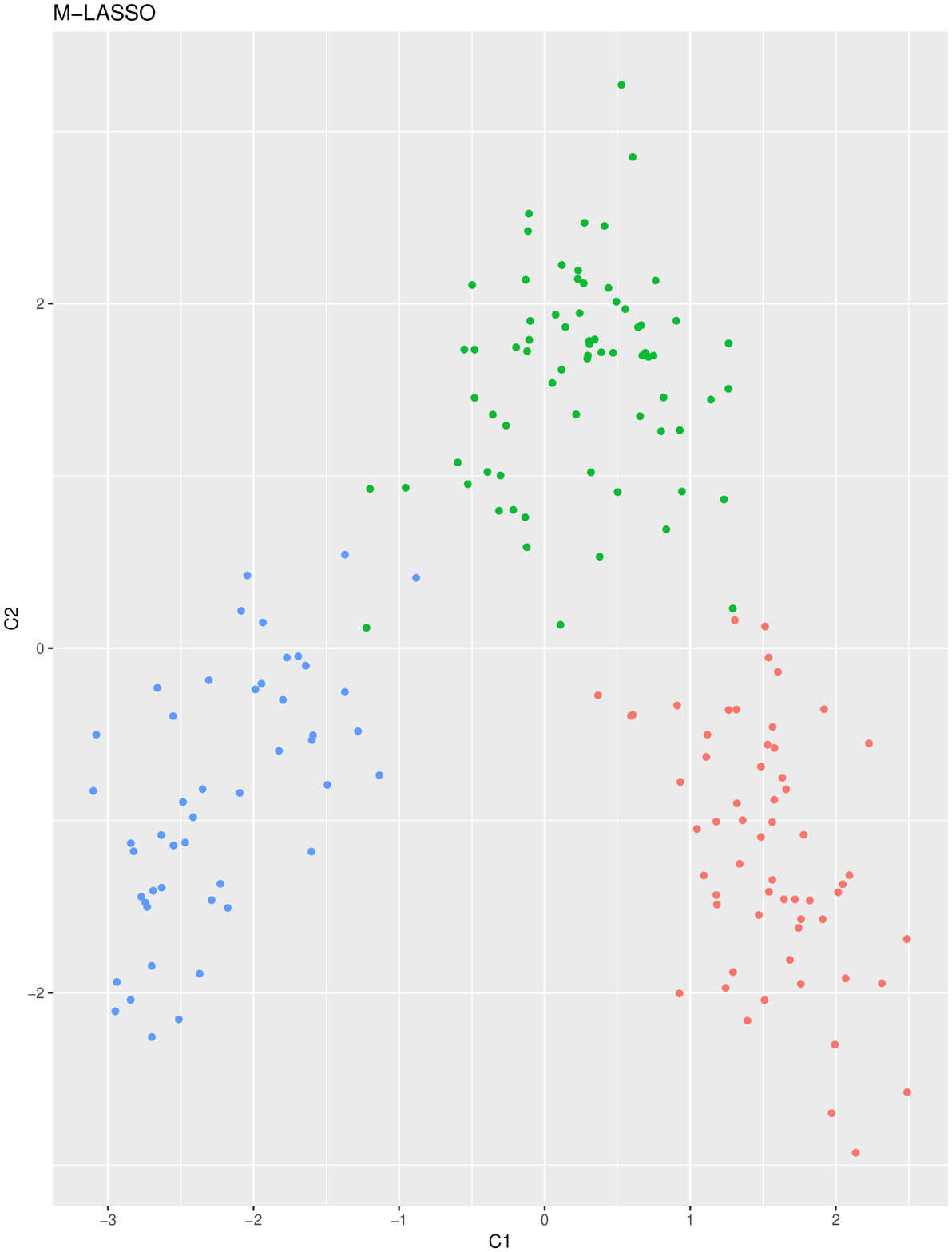}
  \includegraphics[width=60mm, height=60mm]{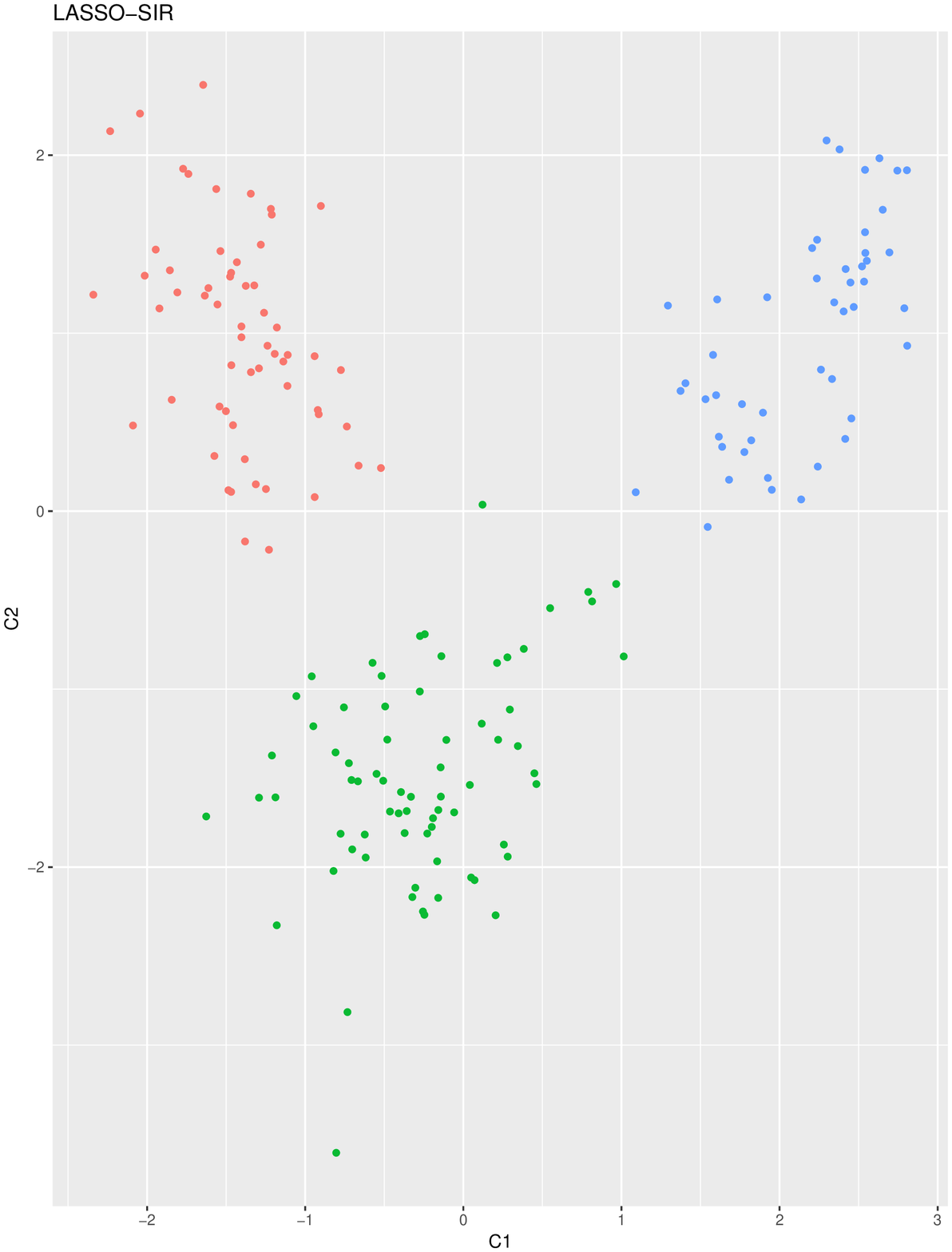}
\caption{We plotted the second component versus the first component for all the wines, which are labeled with different colors, representing different cultivars (1--red, 2--green, 3--blue). The four methods for calculating the directions are PCA, DT-SIR, M-Lasso, and Lasso-SIR from top-left to bottom-right. It is clearly seen that Lasso-SIR offered the best separation among these three groups.
}\label{fig:wine:1}
\end{figure}

The number of directions $d$ is chosen as 2 according to Algorithm \ref{alg:choose:d}. We tested PCA, DT-SIR, M-Lasso, and Lasso-SIR, for  obtaining these two directions. In Figure \ref{fig:wine:1}, we plotted the projection of the data onto the space spanned by two components. The colors of the points correspond to three different cultivars. 
It is clearly seen that  Lasso-SIR provided the best separation of the three cultivars. When using one vertical and one horizontal line to classify three groups, only one subject would be wrongly classified.

\section{Discussion}\label{sec:discussion}
Researchers have made some attempts to extend Lasso to non-linear regression models in recent years (e.g.,\cite{Plan:Vershynin:2016}, \cite{neykov2016support}). However, these approaches are not efficient enough for SDR problems.
In comparison, Lasso-SIR introduced in this article is an efficient high-dimensional variant of SIR \citep{li1991sliced} for obtaining a sparse solution to the estimation of the SDR subspace for multiple index models. We showed that Lasso-SIR is rate optimal if $n\lambda=p^{\alpha}$ for some $\alpha>1/2$, where $\lambda$ is the smallest nonzero eigenvalue of  $var(\bbE[\vx | y])$.
This technical assumption on $n$, $\lambda$, and $p$ is slightly disappointing from the ultra-high dimensional perspective. 
We believe that this technical assumption arises from an intrinsic limitation in estimating the central subspace,
i.e., some further sparsity assumptions on either
$\bSigma$ or $var(\bbE[\vx|y])$ or both are needed to show the consistency of any estimation method.
We will address such extensions in our future researches. 

Cautious reader may find that the  concept of ``pseudo-response variable'' is not essential for developing the theory of the Lasso-SIR algorithm. 
However, by re-formulating the SIR method as a linear regression problem using the pseudo-response variable,  
we can formally consider the model selection consistency, regularization path and many others for multiple index models. In other words, the Lasso-SIR does not only provide an efficient high dimensional variant of SIR, but also extends the rich theory developed for Lasso linear regression in the past decades to the semi-parametric index models. 

The R-package, {\it LassoSIR}, is available on CRAN (\url{https://cran.r-project.org/package=LassoSIR}).

{\section{Acknowledgement}\label{sec:acknow}

  Jun S. Liu is partially supported by the NSF Grants  DMS-1613035 and DMS-1713152, and NIH  Grant R01 GM113242-01. Zhigen Zhao is partially supported by the NSF Grant IIS-1633283.
\newpage  

\bibliographystyle{plainnat}
\bibliography{sir}

\newpage

\begin{appendices}
\section{Appendix: Sketch of Proofs of Theorems 1 , 2 and 3}\label{append:proof} 
We assume the condition {\bf A1)}， {\bf A2)} and {\bf A3)} hold throughout of the rest of the paper.  In particular, the sliced stability condition {\bf A3)} requires  that $H>d$ is a large enough but finite integer. We denote the SIR estimate of $\bLambda=var(\bbE[\vx\mid y])$  by $\widehat{\bLambda}_{H}=\frac{1}{H}\vX_{H}\vX_{H}^{\tau}$  (see e.g., \eqref{eqn:lambda} ) and its eigenvector of unit length associated to the $j$-th eigenvalue $\widehat{\lambda}_{j}$ by $\widehat{\boldeta}_{j}$. To avoid unnecessary confusion, we assume that $\frac{s\log(p)}{n\lambda}$ and $\frac{\sqrt{p}}{n\lambda}$ are sufficiently small. We call an event $\Omega$ happens with high probability if $\bbP(\Omega^{c})\leq C_{1}\exp\left(-C_{2}\log(p) \right)$ for some absolute constants $C_{1}$ and $C_{2}$.

\subsection{Assistant Lemmas}

\subsubsection{Concentration Inequalities}

\begin{lem}\label{lem:deviation:bound:diagonal}
Let $d_{1}, \ldots, d_{p}$ be positive constants. We have the following statements:
\begin{itemize}
\item[i)]For $p$  $i.i.d.$ standard normal random variables $x_{1},\ldots,x_{p}$, there exist constants $C_{1}$ and $C_{2}$ such that for any sufficiently small $a$ , we have
\begin{align}\label{inline:tempppp}
\bbP\left(\left|\frac{1}{p}\sum_{i}d_{i}(x_{i}^{2}-1)\right|>a\right)
\leq C_{1}\exp(-\frac{p^{2}a^{2}}{C_{2}\sum_{j}d_{j}^{4}}).
\end{align}
\item[ii)]For $2p$  $i.i.d.$ standard normal random variables $x_{1},\ldots,x_{p},$ $y_{1},\cdots,$ and $y_{p}$, there exist constants $C_{1}$ and $C_{2}$ such that for any sufficiently small $a$ , we have
\begin{align}
\bbP\left(|\frac{1}{p}\sum_{i}d_{i}x_{i}y_{i}|>a\right)
\leq C_{1}\exp(-\frac{p^{2}a^{2}}{C_{2}\sum_{j}d_{j}^{4}}).
\end{align}
\end{itemize}

\proof  ii) is a direct corollary of i). We put the proof of i) in the supplementary materials.

\end{lem}

\subsubsection{Sine-Theta Theorem}

\begin{lem}[Sine-Theta Theorem] \label{lem:sin_theta} 
Let $\vA$ and $\vA+\vE$ be symmetric matrices satisfying
\begin{equation*}
\vA=[\vF_{0},\vF_{1}]\left[ \begin{array}{cc}
\vA_{0} & 0\\
0 & \vA_{1}
\end{array} \right]
\left[ \begin{array}{c}
\vF^{\tau}_{0}\\
\vF^{\tau}_{1}
\end{array} \right],
\quad
\vA+\vE=[\vG_{0},\vG_{1}]\left[ \begin{array}{cc}
\bLambda_{0} & 0\\
0 & \bLambda_{1}
\end{array} \right]
\left[ \begin{array}{c}
\vG^{\tau}_{0}\\
\vG^{\tau}_{1}
\end{array} \right]
\end{equation*}
where $[\vF_{0},\vF_{1}]$ and $[\vG_{0},\vG_{1}]$ are orthogonal matrices. If the eigenvalues of $\vA_{0}$ are contained in an interval (a,b) , and the eigenvalues of $\bLambda_{1}$ are excluded from the interval $(a-\delta,b+\delta)$ for some $\delta>0$, then
\[
\|\vF_{0}\vF_{0}^{\tau}-\vG_{0}\vG_{0}^{\tau}\|\leq \frac{\min(\|\vF_{1}^{\tau}\vE\vG_{0}\|, \|\vF_{0}^{\tau}\vE\vG_{1}\|)}{\delta},
\]
and 
\[
\frac{1}{\sqrt{2}}\|\vF_{0}\vF_{0}^{\tau}-\vG_{0}\vG_{0}^{\tau}\|_{F} \leq \frac{\min(\|\vF_{1}^{\tau}\vE\vG_{0}\|_{F}, \|\vF_{0}^{\tau}\vE\vG_{1}\|_{F})}{\delta}.
\]
\end{lem}

\subsubsection{Restricted Eigenvalue Properties}
We briefly review the restricted eigenvalue (RE) property,    
 which was first introduced in \cite{raskutti2010restricted}.
Given a set $S\subset	[p]=\{1,...,p\}$, for any positive number $\alpha$, define the set $\mathcal{C}(S,\alpha)$ as 
\begin{align*}
\mathcal{C}(S,\alpha)=\{\theta\in \bbR^{p}~|~\|\theta_{S^{c}}\|_{1}\leq \alpha \|\theta_{S}\|_{1}\}.
\end{align*}
We say that a sample matrix $\vX\vX^{\tau}/n$ satisfies the restricted eigenvalue condition over $S$ with parameter $(\alpha,\kappa)\in [1,\infty)\times (0,\infty)$ if 
\begin{align}\label{cond:re}
\frac{1}{n}\theta^{\tau}\vX\vX^{\tau}\theta\geq \kappa^{2}\|\theta\|_{2}^{2}, \quad \forall \theta\in \mathcal{C}(S,\alpha).
\end{align}
If (\ref{cond:re}) holds uniformly for all the subsets $S$ with cardinality $s$, we say that $\vX\vX^{\tau}/n$ satisfies the restricted eigenvalue condition of order $s$ with parameter $(\alpha, \kappa)$. 
Similarly, we say that the covariance matrix $\bSigma$ satisfies the RE condition over $S$ with parameter $(\alpha,\kappa)$ if $\|\bSigma^{1/2}\theta\|_{2}\geq \kappa \|\theta\|$ for all $\theta \in \mathcal{C}(S,\alpha)$. 
Additionally, if this condition holds uniformly for all the subsets $S$ with cardinality $s$, we say that $\bSigma$ satisfies the restricted eigenvalue condition of order $s$ with parameter $(\alpha, \kappa)$. 
The following Corollary is borrowed from \cite{raskutti2010restricted}.

\begin{cor}\label{cor:RE}
Suppose that $\bSigma$ satisfies the RE condition of order $s$ with parameter $(\alpha,\kappa)$. 
Let $\vX$ be the $p\times n$ matrix formed by $n$ $i.i.d$ samples from $N(0,\bSigma)$.
For some universal positive constants $\mathfrak{a}_{1},\mathfrak{a}_{2} $ and $\mathfrak{a}_{3}$, if the sample size satisfies 
\begin{align*}
n>\mathfrak{a}_{3}\frac{(1+\alpha)^{2}\max_{i \in [p]}\bSigma_{ii}}{\kappa^{2}}s\log(p),
\end{align*}
then the matrix $\frac{1}{n}\vX\vX^{\tau}$ satisfies the RE condition of order $s$ with parameter $(\alpha, \frac{\kappa}{8})$ with probability at least $1-\mathfrak{a}_{1}\exp\left(-\mathfrak{a}_{2}n\right)$.
\end{cor}

It is clear that  $\lambda_{\min}(\bSigma)\geq C_{\min}$ implies that $\bSigma$ satisfies the RE condition of any order $s$ with parameter  $(3, \sqrt{C_{\min}})$.   
Thus,   we have the following proposition.
\begin{prop}\label{prop:second}
For some universal constants $\mathfrak{a}_{1}$, $\mathfrak{a}_{2}$ and $\mathfrak{a}_{3}$, if the sample size satisfies that 
$
n>\mathfrak{a}_{1}s\log(p),
$
then the matrix $\frac{1}{n}\vX\vX^{\tau}$ satisfies the RE condition for any order $s$ with parameter $(3,\sqrt{C_{\min}}/8)$ with probability at least $1-\mathfrak{a}_{2}\exp\left(-\mathfrak{a}_{3}n \right)$.
\end{prop}

\subsubsection{The sliced approximation Lemma}\label{append:discussion}
Let $\vx \in \mathbb{R}^{p}$ be a sub-Gaussian random variable.  For any unit vector $\bbeta \in \mathbb{R}^{p}$, let $\vx(\bbeta)=\langle\vx,\bbeta\rangle$ and $\vm(\bbeta)=\langle\vm,\bbeta\rangle=\bbE[\vx(\bbeta)\mid y]$. 
In order to get the deviation properties of the statistics $var_{H}(\vx(\bbeta))$, \cite{lin2015consistency} has introduced the sliced stable condition, i.e., the condition $A3$ in this paper.  For the exact definition and more discussion, we refer to \cite{lin2015consistency}. 

\begin{lem}\label{lem:new_statistic:deviation} Let $\vx \in \mathbb{R}^{p}$ be a sub-Gaussian random variable.  Assume that $\bbE[\vx|y]$ is sliced table with respect to $y$. For any unit vector $\bbeta \in \mathbb{R}^{p}$, let $\vx(\bbeta)=\langle\vx,\bbeta\rangle$ and $\vm(\bbeta)=\langle\vm,\bbeta\rangle=\bbE[\vx(\bbeta)\mid y]$, we have the following:
\begin{itemize}
\item[i)] If $var(\vm(\bbeta))= 0$ , there exist positive constants $C_{1}, C_{2}$ and $C_{3}$  such that  for any $b$ and sufficiently large $H$,  we have
\begin{equation*} \label{eqn:deviation:outSets}
\begin{aligned}
\bbP(var_{H}(\vx(\bbeta))>b) 
\leq  C_{1}\exp\left(-C_{2}\frac{nb}{H^{2}} +C_{3}\log(H) \right).
\end{aligned}
\end{equation*}
\item[ii)] If $var(\vm(\bbeta))\neq 0$ , there exist positive constants $C_{1}, C_{2}$ and $C_{3}$ such that,  for any $\nu>1$,  we have
\[
|var_{H}(\vx(\bbeta))-var(\vm(\bbeta))|\geq \frac{1}{2\nu}var(\vm(\bbeta))
\]
with probability at most 
\[
C_{1}\exp\left( -C_{2}\frac{n~var(\vm(\bbeta))}{H^{2}\nu^{2}}+C_{3}\log(H)\right).
\]
where we choose $H$ such that $H^{\vartheta}>C_{4}\nu$ for some sufficiently large constant $C_{4}$.
\end{itemize}
\end{lem}
The following proposition is a direct corollary.
\begin{prop}\label{prop:knowledge}
There exist positive constants $C_{1}$,$C_{2}$ and $C_{3}$, such that  
\begin{align}
\|\bbeta^{\tau}\bLambda_{H}\bbeta-\bbeta^{\tau}var(\bbE[\vx|y])\bbeta\|_{2}
\geq \frac{1}{2\nu}\bbeta^{\tau}var(\bbE[\vx|y])\bbeta
\end{align}
with probability at most $C_{1}\exp\left(-C_{2}\frac{n\lambda}{H^{2}\nu^{2}} +C_{3}\log(H)\right)$.
\proof It follows from Lemma \ref{lem:new_statistic:deviation} and the fact that for any $\bbeta \in col(var(\bbE[\vx|y]))$, $var(\vm(\bbeta)))\geq \lambda$.\epf
\end{prop}

\subsubsection{Properties of $\widehat{\eta}_{j}$ 's. }
\begin{prop} \label{prop:key} 
Recall that $\widehat{\boldeta}_{j}$ is the eigenvector associated to the $j$-th eigenvalue of $\widehat{\bLambda}_{H}, j=1,...,H$. If $n\lambda=p^{\alpha}$ for some $\alpha>1/2$, there exist positive constants $C_{1}$ and $C_{2}$ such that
\begin{itemize}
\item[i)] for $j=1,....,d$, one has
\begin{align}\label{true:projection}
\|P_{\bLambda}\widehat{\boldeta}_{j}\|_{2}\geq C_{1} \sqrt{\frac{\lambda}{\widehat{\lambda}}_{j}}
\end{align}
\item[ii)] for $j=d+1,...,H$,  one has
\begin{align}\label{true:projection}
 \|P_{\bLambda}\widehat{\boldeta}_{j}\|_{2}\leq C_{2}\frac{\sqrt{p\log(p)}}{n\lambda} \sqrt{\frac{\lambda}{\widehat{\lambda}}_{j}}
\end{align}
\end{itemize}
hold with high probability.

\noindent{\it Remark:} \normalfont This result might be of independent interest.
In order to justify that the sparsity assumption for the high dimensional setting is necessary,
 \cite{lin2015consistency} have shown that for single index models, $\bbE[\angle(\boldeta_{1},\widehat{\boldeta}_{1})]=0$ if and only if $\lim\frac{p}{n}=0$.  
Proposition \ref{prop:key} states that the projection of $\sqrt{\widehat{\lambda}_{j}}\widehat{\boldeta}_{j}$, $j=1,\ldots,d$, onto the true direction is non-zero if $n\lambda>p^\alpha$ where $\alpha>1/2$. 

\medskip

\proof Let $\vx=\vz+\vw$ be the orthogonal decomposition with respect to $col(var(\bbE[\vx|y]))$ and its orthogonal complement. 
We define two $p\times H$ matrices $\vZ_{H}=(\vz_{1,\cdot},\ldots,\vz_{H,\cdot})$ and $\vW_{H}=(\vw_{1,\cdot},\ldots,\vw_{H,\cdot})$ whose definition are similar to the definition of $\vX_{H}$. 
We then have the following decomposition
\begin{align}\label{eqn:decomposition}
\vX_{H}=\vZ_{H}+\vW_{H}.
\end{align}
By definition, we know that 
$\vZ_{H}^{\tau}\vW_{H}=0$ and
$y\independent \vw$.
Let  $\bSigma_{1}$ be the covariance matrix of $\vw$, then $\vW_{H}=\frac{1}{\sqrt{c}}\bSigma_{1}^{1/2}\vE_{H}$ where  $\vE_{H}$ is a $p\times H$ matrix with $i.i.d.$ standard normal entries. 

For sufficiently large $\nu_{1} $ and $\mathfrak{a}$, Lemma \ref{lem:new_statistic:deviation} implies that 
\begin{align}
\Omega_{1}&=\Big\{ \omega \mid (1-\frac{\kappa}{2\nu_{1}})\lambda \leq \lambda_{\min}( \frac{1}{H}\vZ_{H}^{\tau}\vZ_{H} )\leq \lambda_{\max}( \frac{1}{H}\vZ_{H}^{\tau}\vZ_{H} ) \leq (1+\frac{1}{2\nu_{1}})\kappa\lambda  \Big\} 
\end{align}
happens with high probability and Lemma \ref{lem:deviation:bound:diagonal} implies
\begin{align}
\Omega_{2}&=\Big\{\omega\mid~\left\|\frac{1}{H}\vW_{H}^{\tau}\vW_{H}-\frac{tr(\bSigma_{1})}{n}\vI_{H}\right\|_{F}\leq \mathfrak{a}\frac{\sqrt{p\log(p)}}{n} ~\Big\}
\end{align}
happens with high probability.

For any  $\omega \in \Omega=\Omega_{1}\cap \Omega_{2}$,  we can choose a $p\times p$ orthogonal matrix $T$ and an $H\times H$ orthogonal matrix $S$ such that 
\begin{align}
\frac{1}{H}T\vZ_{H}(\omega)S=\left(\begin{array}{cc}
B_{1} & 0 \\
0 & 0 \\
0 & 0
\end{array} \right) \mbox{ and } 
\frac{1}{H}T\vW_{H}(\omega)S =\left(\begin{array}{cc}
0 & 0\\
B_{2} & B_{3}\\
0 & B_{4} 
\end{array} \right)
\end{align}
where $B_{1}$ is a $d\times d $ matrix, $B_{2}$ is a $d\times d$ matrix, $B_{3}$ is a $d\times (H-d)$ matrix and $B_{4}$ is a $(p-2d)\times (H-d)$ matrix. By definition of the event $\omega$, we have
\begin{equation}\label{Condition:C}
\begin{aligned}
& (1-\frac{\kappa}{2\nu_{1}})\lambda \leq \lambda_{\min}( B_{1}^{\tau}B_{1} )  \leq \lambda_{\max}( B_{1}^{\tau}B_{1} )\leq (1+\frac{1}{2\nu_{1}})\kappa\lambda  \\
& \left\|\left(\begin{array}{cc}
B_{2}^{\tau}B_{2} & B_{2}^{\tau}B_{3} \\
B_{3}^{\tau}B_{2} & B_{3}^{\tau}B_{3}+B_{4}^{\tau}B_{4}
\end{array} \right)-\frac{tr(\bSigma_{1})}{n}\vI_{H}\right\|_{F}\leq \mathfrak{a}\frac{\sqrt{p\log(p)}}{n}
\end{aligned} \quad .
\end{equation}
Proposition \ref{prop:key} follows from the linear algebraic lemma below:
\begin{lem}\label{lem:prop1:2}  Assume that $n\lambda=p^{\alpha}$ for some $\alpha>1/2$.   To avoid unnecessary confusion, we also assume that $\frac{\sqrt{p\log(p)}}{n\lambda}$ is sufficiently small.
Let $\vM=\left(\begin{array}{cc}
B_{1} & 0\\
B_{2} & B_{3}\\
0 & B_{4} 
\end{array} \right)$ be a $p\times H$ matrix, where $B_{1}$ is a $d\times d$ matrix, $B_{2}$ is a $d\times d$ matrix, $B_{3}$ is a $d\times (H-d)$ matrix and $B_{4}$ is a $(p-2d)\times(H-d)$ matrix satisfying \eqref{Condition:C}. 
Let $\widehat{\boldeta}_{j}$ be the eigenvector associated with the $j$-th eigenvalue $\widehat{\lambda}_{j}$ of $\vM\vM^{\tau}$, $j=1,...,H$. 
Then the length of the projection of $\widehat{\boldeta}_{j}$ onto its first $d$-coordinates is at least $C\sqrt{\frac{\lambda}{\widehat{\lambda}_{j}}}$ for $j=1,...,d$ and is at most $C\frac{\sqrt{p\log(p)}}{n\lambda}\sqrt{\frac{\lambda}{\widehat{\lambda}_{j}}}$ for $j=d+1,...,H$.
\proof 
Let us consider the eigen-decompositions of 
\begin{align*}
Q_{1}\triangleq\vM^{\tau}\vM=
\left(\begin{array}{cc}
\vE_{1} & \vE_{2}\\
\vE_{3} & \vE_{4}
\end{array} \right)
\left(\begin{array}{cc}
\vD_{1} & 0\\
0 & \vD_{2}
\end{array} \right)
\left(\begin{array}{cc}
\vE_{1}^{\tau} & \vE_{3}^{\tau}\\
\vE_{2}^{\tau} & \vE_{4}^{\tau}
\end{array} \right)
\end{align*}
where $\vD_{1}$ ( {\it resp.} $\vD_{2}$)  is a $d\times d$ ({\it resp. $(H-d)\times (H-d)$} )  diagonal matrices satisfying that $\lambda_{\min}(\vD_{1})\geq \lambda_{\max}(\vD_{2})$.
\eqref{Condition:C} implies that
\begin{align*}
(1-\frac{\kappa}{2\nu_{1}})\lambda+\frac{tr(\bSigma_{1})}{n}-&\mathfrak{a}\frac{\sqrt{p\log(p)}}{n} \leq \lambda_{\min}(D_{1})\\
&\leq \lambda_{\max}(D_{1}) \leq (1+\frac{1}{2\nu_{1}})\kappa\lambda+\frac{tr(\bSigma_{1})}{n}+\mathfrak{a}\frac{\sqrt{p\log(p)}}{n} .
\end{align*}
On the other hand, we could consider the eigen-decomposition of 
\begin{align*}
Q_{2}\triangleq&\left(\begin{array}{cc}
B_{1}^{\tau}B_{1}+\frac{tr(\bSigma_{1})}{n}\vI_{d} & 0\\
0 & \frac{tr(\bSigma_{1})}{n}\vI_{H-d}
\end{array} \right)=\left(\begin{array}{cc}
\vF_{1} & 0\\
0 & \vF_{2}
\end{array} \right)
\left(\begin{array}{cc}
\vD'_{1} & 0\\
0 & \vD'_{2}
\end{array} \right)
\left(\begin{array}{cc}
\vF_{1}^{\tau} & 0\\
0 & \vF_{2}^{\tau}
\end{array} \right)
\end{align*}
where $\vD_{1}'$ ( {\it resp.} $\vD_{2}'$)  is a $d\times d$ ({\it resp. $(H-d)\times (H-d)$} )  diagonal matrices satisfying that $\lambda_{\min}(\vD_{1}')\geq \lambda_{\max}(\vD_{2}')$.
\eqref{Condition:C} implies that
\begin{align*}
 \frac{tr(\bSigma_{1})}{n}-\mathfrak{a}\frac{\sqrt{p\log(p)}}{n}&\leq
\lambda_{\min}(\vD_{2}') \leq \lambda_{\max}(\vD_{2}') \leq \frac{tr(\bSigma_{1})}{n}+\mathfrak{a}\frac{\sqrt{p\log(p)}}{n}.
\end{align*}
Thus the eigen-gap is of order $\lambda-\mathfrak{a}\frac{\sqrt{p\log(p)}}{n}$ ( which is of order $\lambda$, since $n\lambda=p^{\alpha}$ for some $\alpha>1/2$). From \eqref{Condition:C} , we know that $\|Q_{1}-Q_{2}\|_{F}\leq C\frac{\sqrt{p\log(p)}}{n}$.
The Sine-Theta theorem (see e.g., Lemma \ref{lem:sin_theta}) implies that
\begin{align}
\left\|
\left(
\begin{array}{c}
\vE_{1}\\
\vE_{3}
\end{array}
\right)
\left(
\begin{array}{cc}
\vE_{1}^{\tau} &
\vE_{3}^{\tau}
\end{array}\right)-
\left(
\begin{array}{cc}
\vI_{d} & 0\\
0 & 0
\end{array}
\right)
\right\|_{F}\leq C\frac{\sqrt{p\log(p)}}{n\lambda},
\end{align}
{\it i.e.}, $\|\vE_{3}\vE_{3}^{\tau}\|_{F}\leq C\frac{\sqrt{p\log(p)}}{n\lambda}$.
Similar argument gives us that
$\|\vE_{2}\vE_{2}^{\tau}\|_{F}\leq C\frac{\sqrt{p\log(p)}}{n\lambda}$.

\vspace*{3mm}

Let $\boldeta$ be the (unit) eigenvector associated to the non-zero eigenvalue $\widehat{\lambda}$ of $\vM\vM^{\tau}$. 
Let us write 
$\boldeta^{\tau}=(\boldeta_{1}^{\tau},\boldeta_{2}^{\tau},\boldeta_{3}^{\tau})$
where $\boldeta_{1}, \boldeta_{2} \in \bbR^{d}$ and $\boldeta_{3} \in \bbR^{p-2d}$. 
Let  
$\alpha=(\alpha_{1}^{\tau},\alpha_{2}^{\tau})$
where
$\alpha_{1}=B_{1}^{\tau}\boldeta_{1}+B_{2}^{\tau}\boldeta_{2}\in \bbR^{d}$ and $\alpha_{2}=B_{3}^{\tau}\boldeta_{2}+B_{4}^{\tau}\boldeta_{3}\in \bbR^{H-d}$.
It is easy to verify that $\alpha/\sqrt{\widehat{\lambda}}$ is the (unit) eigenvector associated to the eigenvalue $\widehat{\lambda}$ of $\vM^{\tau}\vM$ and
\begin{align*}
\boldeta_{1}=\frac{B_{1}}{\sqrt{\widehat{\lambda}}}\frac{\alpha_{1}}{\sqrt{\widehat{\lambda}}},
\quad \boldeta_{2}=\frac{B_{2}}{\sqrt{\widehat{\lambda}}}\frac{\alpha_{1}}{\sqrt{\widehat{\lambda}}}
+\frac{B_{3}}{\sqrt{\widehat{\lambda}}}\frac{\alpha_{2}}{\sqrt{\widehat{\lambda}}},
\mbox{ and } \boldeta_{3}=\frac{B_{4}}{\sqrt{\widehat{\lambda}}}\frac{\alpha_{2}}{\sqrt{\widehat{\lambda}}}.
\end{align*}
If $\widehat{\lambda}$ is  among the first $d$ eigenvalues of $\vM^{\tau}\vM$, then $\|\alpha_{1}/\sqrt{\widehat{\lambda}}\|_{2}$ is bounded below by some positive constant. 
Thus $\|\boldeta_{1}\|_{2}\geq C\sqrt{\frac{\lambda}{\widehat{\lambda}}}$ . 
If $\widehat{\lambda}$ is  among the last $H-d$ eigenvalues of $\vM^{\tau}\vM$, then $\|\alpha_{1}/\sqrt{\widehat{\lambda}}\|_{2}=O\left( \frac{\sqrt{p\log(p)}}{n\lambda}\right)$.
Thus $\|\boldeta_{1}\|_{2}\leq O\left(\kappa\sqrt{\frac{\lambda}{\widehat{\lambda}}} \frac{\sqrt{p\log(p)}}{n\lambda} \right)$ .
\epf

\end{lem}

\end{prop}

\subsection{Sketch of Proof of Theorem \ref{thm:consistency}}
We only  sketch some key points of the proof here  and leave the details in the online supplementary files.
Recall that for single index model
$
y=f(\bbeta_{0}^{\tau}\vx, \epsilon)
$
where $\bbeta_{0}$ is a unit vector, we have denoted by $\widehat{\boldeta}$ the eigenvector of $\widehat{\bLambda}_{H}$ associated to the largest eigenvalue $\widehat{\lambda}$.
  Let $\widehat{\bbeta}$ be the minimizer of 
\begin{align*}
\mathcal{L}_{\bbeta}=\frac{1}{2n}\|\widetilde{\vy}-\vX^{\tau}\bbeta\|+\mu\|\bbeta\|_{1},
\end{align*}
where $\widetilde{\vy}\in \bbR^{n}$ such that $\widehat{\boldeta}=\frac{1}{n}\vX\widetilde{\vy}$.
Let $\boldeta_{0}=\bSigma\bbeta_{0}$, $\widetilde{\boldeta}=P_{\eta_{0}}\widehat{\boldeta}$
and $\widetilde{\bbeta}=\bSigma^{-1}\widetilde{\boldeta}\propto \bbeta_{0}$.
Since we are interested in the distance between the directions of $\widehat{\bbeta}$ and $\bbeta_{0}$,  we consider the difference  $\bdelta=\widehat{\bbeta}-\widetilde{\bbeta}$. 
A slight modification of the argument in \cite{bickel2009simultaneous} implies that, if we choose $\mu= C\sqrt{\frac{\log(p)}{n\widehat{\lambda}}}$ for sufficiently large constant $C$,  we have
$
\|\bdelta\|_{2}\leq C_{1}\sqrt{\frac{s\log(p)}{n\widehat{\lambda}}}
$
with high probability.
The detailed arguments are put in the online supplementary file. 
The Proposition \ref{prop:key}, Condition
${\bf A1)}$ and $\widetilde{\bbeta}=\bSigma^{-1}\widetilde{\boldeta}, $
imply that
$
 C_{1}\sqrt{\frac{\lambda}{\widehat{\lambda}}}\leq\|\widetilde{\bbeta}\|_{2}\leq C_{2} 
$
holds with high probability for some constants $C_{1}$ and $C_{2}$.
Thus, we have
\begin{align}
\|P_{\widehat{\bbeta}}-P_{\bbeta_{0}}\|_{F}=\|P_{\widehat{\bbeta}}-P_{\widetilde{\bbeta}}\|_{F}\leq 4\frac{\|\widehat{\bbeta}-\widetilde{\bbeta}\|_{2}}{\|\widetilde{\bbeta}\|_{2}}=4\|\bdelta\|_{2}/\|\widetilde{\bbeta}\|_{2}\leq C\sqrt{\frac{s\log(p)}{n\lambda}}
\end{align}
holds with high probability.
\epf

\subsection{Proof of Theorem \ref{thm:consistency:multiple:version1} }
Recall that $\widehat{\boldeta}_{j}$'s are the (unit) eigenvectors associated to the $j$-th eigenvalues of $\widehat{\bLambda}_{H}$, $j=1,....,d$.
We introduce the following notations,
\begin{align}\label{notation:inline}
\widetilde{\boldeta}_{j}=P_{\bLambda}
\widehat{\boldeta}_{j}, \quad  \widetilde{\bbeta}_{j}=\bSigma^{-1}\widetilde{\boldeta}_{j} \quad \mbox{ and }
\bgamma_{j}=\widetilde{\bbeta}_{j}/\|\widetilde{\bbeta}_{j}\|_{2}.
\end{align}
Applying the argument in Theorem \ref{thm:consistency} on  these eigenvectors, we have
\begin{align}\label{eqn:tempp:distance}
\|\widehat{\bbeta}_{j}-\widetilde{\bbeta}_{j}\|_{2}\leq C\sqrt{\frac{s\log(p)}{n\widehat{\lambda}_{j}}} \mbox{ and }
\|P_{\widehat{\bbeta}_{j}}-P_{\widetilde{\bbeta}_{j}}\|_{F}
\leq C\sqrt{\frac{s\log(p)}{n\lambda}}
\end{align} 
for some constant $C$ hold with high probability. 
Since we assume that $d$ is fixed, if we can prove that
\begin{itemize}
\item[I)]  the lengths of  $\widetilde{\bbeta}_{j}, j=1,...,d,$ are bounded below by $C\sqrt{\frac{\lambda}{\widehat{\lambda}_{j}}}$,
\item[II)] the angles between any two vectors of $\widetilde{\bbeta}_{j}, j=1,...,d,$ are bounded below by some constant,
\end{itemize}
 hold with probability, then the Gram-Schmit process implies that
$
\|P_{\widehat{\vB}}-P_{\vB}\|_{F} \leq C\sqrt{\frac{s\log(p)}{n\lambda}}
$
holds with high probability from \eqref{eqn:tempp:distance}.
It is easy to verify that $I)$ follows from the Proposition \ref{prop:key} , the Condition ${\bf A1)}$ and the definition of $\widetilde{\bbeta}_{j}( ~=
\bSigma^{-1}\widetilde{\boldeta}_{j}), j=1,...,d$.
$II)$ is a direct corollary of the following two statements.

\paragraph{{\bf Statement A. The angles between any two vectors in $\widetilde{\boldeta}_{j}$ $'s$ , $j=1,...,d$ are nearly $\pi/2$.}}  Since $n\lambda=p^{\alpha}$ for some $\alpha>1/2$,
we only need to prove that 
\begin{align}
\left|\cos\left(\angle(\widetilde{\boldeta}_{j}, \widetilde{\boldeta}_{j}) \right)\right|
\leq C\frac{\sqrt{p\log(p)}}{n\lambda}
\end{align}
holds with high probability for any $i\neq j$.
Recall that we have the following decomposition
$
\vX_{H}=\vZ_{H}+\vW_{H}.
$
It is easy to see that $col(\vZ_{H})=col(var(\bbE[\vx|y]))$ and  $\sqrt{n}~cov(\vw)^{-1/2}\vW_{H}$ is identically distributed to a matrix, $\mathcal{E}_{1}$, with all the entries are $i.i.d.$ standard normal random variables. 
Let us choose an orthogonal matrix $T$ such that $\frac{1}{\sqrt{H}}T\vZ_{H}=(\vA^{\tau},0)^{\tau}$ and  
$\frac{1}{\sqrt{H}}T\vW_{H}=(0,\vB^{\tau})^{\tau}$
where $\vA$ is a $d\times H$ matrix and $\vB$ is a $(p-d)\times$ H matrix. 
Thus, $T\widehat{\boldeta}_{j}$ is the eigenvector of $\frac{1}{H}T\vX_{H}\vX_{H}^{\tau}T^{\tau}$ associated with the $j$-th eigenvalue $\widehat{\lambda}_{j}$, $j=1,...,d$. 
If we have $a)$ $\lambda_{\min}(\vA\vA^{\tau})\geq \lambda$, $b)$ $\|P_{col(T\vZ_{H})}(T\widehat{\boldeta}_{j})\|_{2}\geq C\sqrt{\frac{\lambda}{\widehat{\lambda}_{j}}}$ and $c)$ $\|\vB^{\tau}\vB-\mu\vI_{H}\|_{F}\leq C\frac{\sqrt{p\log(p)}}{n}$ for some scalar $\mu>0$,
then the statement {\bf I} is reduced to the following linear algebra lemma.
\begin{lem}
Let $\vA$ be a $d\times H $  matrix $(d<H)$ with $\lambda_{\min}(\vA\vA^{\tau})=\lambda$. Let $\vB$ be a $(p-d)\times H$ matrix such that 
$\|\vB^{\tau}\vB-\mu \vI_{H}\|_{F}^{2}\leq C\frac{\sqrt{p\log(p)}}{n}$.
Let $\widehat{\bxi}_{j}$ be the $j$-th (unit) eigenvector of 
$\vC\vC^{\tau}$ associated with the $j$-th eigenvalue $\widehat{\lambda}_{j}$ where $\vC^{\tau}=(\vA^{\tau},\vB^{\tau})$ and $\widetilde{\bxi}_{j}$ be the projection of $\widehat{\bxi}_{j}$ onto its first $d$-coordinates. 
If $\|\widetilde{\bxi}_{j}\|_{2}\geq C\sqrt{\frac{\lambda}{\widehat{\lambda}_{j}}}$ ,  then for any $i\neq j$,
\begin{align}
\left|\cos\left(\angle(\widetilde{\bxi}_{i}, \widetilde{\bxi}_{j}) \right)\right|
\leq C\frac{\sqrt{p\log(p)}}{n\lambda}.
\end{align}
Thus, $\widetilde{\bxi}_{i}$ $'s$ are nearly orthogonal if $n\lambda=p^{\alpha}$ for some $\alpha>1/2$.
\proof Let $\widehat{\balpha}_{j}=\vC^{\tau}\widehat{\bxi}_{j}$, then $\widehat{\bxi}_{j}=\frac{1}{\widehat{\lambda}_{j}}\vC\balpha_{j}$ and 
$\widetilde{\bxi}_{j}=\frac{1}{\widehat{\lambda}_{j}}\vA\balpha_{j}$. 
It is easy to see that $\|\widehat{\balpha}_{j}\|_{2}=\sqrt{\widehat{\lambda}_{j}}$  and $\|\vC\widehat{\balpha}_{j}\|_{2}\geq \widehat{\lambda}_{j}$. 
Since $\widehat{\balpha}_{j}/\sqrt{\widehat{\lambda}_{j}}$ is also the (unit) eigenvector of
\begin{align*}
\vC^{\tau}\vC=\vA^{\tau}\vA+\mu\vI+(\vB^{\tau}\vB-\mu\vI),
\end{align*}
for any $i\neq j$, we have
\begin{align*}
0&=\widehat{\balpha}_{j}^{\tau}\vC^{\tau}\vC\widehat{\balpha}_{i}=
\widehat{\balpha}_{j}^{\tau}\vA^{\tau}\vA\widehat{\balpha}_{j}
+\mu\widehat{\balpha}_{j}^{\tau}\widehat{\balpha}_{i}
+\widehat{\balpha}_{j}^{\tau}(\vB^{\tau}\vB-\mu\vI)\widehat{\balpha}_{i}\\
&=\widehat{\lambda}_{j}\widehat{\lambda}_{i}
\widetilde{\bxi}_{j}^{\tau}\widetilde{\bxi}_{i}+
\widehat{\balpha}_{j}^{\tau}(\vB^{\tau}\vB-\mu\vI)\widehat{\balpha}_{i}.
\end{align*}
Since 
$\|\vB^{\tau}\vB-tr(\bSigma)\vI_{H}\|_{F}\leq C\frac{\sqrt{p\log(p)}}{n}$ and  $\|\widehat{\bxi}_{j}\|_{2}\geq C\sqrt{\frac{\lambda}{\widehat{\lambda}_{j}}}$, 
 $\forall i\neq j$, we have
\begin{align*}
\left|\frac{\bxi_{j}^{\tau}\bxi_{i}}{\|\bxi_{i}\|_{2}\|\bxi_{j}\|_{2}}\right| \leq C\left|\widehat{\bxi}_{j}^{\tau}\widehat{\bxi}_{i}
\frac{\widehat{\lambda}_{j}^{1/2}\widehat{\lambda}_{i}^{1/2}}{\lambda}\right|
=C\left|\frac{1}{\lambda}\frac{\widehat{\balpha}_{j}^{\tau}}{\widehat{\lambda}_{j}^{1/2}}(\vB^{\tau}\vB-\mu\vI)\frac{\widehat{\balpha}_{i}}{\widehat{\lambda}_{i}^{1/2}}\right|\leq C\frac{\sqrt{p\log(p)}}{n\lambda} . 
\end{align*}
\epf
\end{lem}
Note that $a)$ follows from the Lemma \ref{prop:knowledge}, $b)$ follows from Proposition \ref{prop:key} and $c)$ follows from the Lemma \ref{lem:deviation:bound:diagonal}. Thus statement {\bf A} holds. 

\paragraph{{\bf Statement B. The angles between any two vectors in $\widetilde{\bbeta}_{j}$ $'s$ are bounded away from 0.} }
Since $\widetilde{\bbeta}_{j}
=\bSigma^{-1}\widetilde{\boldeta}_{j}$, we only need to prove that there exists a positive constant $\zeta<1$ such that
\begin{align}
\left| \frac{\widetilde{\boldeta}_{i}^{\tau}\bSigma^{-1}
\bSigma^{-1}\widetilde{\boldeta}_{i}}
{\|\bSigma^{-1}\widetilde{\boldeta}_{i}\|_{2}
\|\bSigma^{-1}\widetilde{\boldeta}_{j}\|_{2}}\right| \leq \zeta.
\end{align}
Let $(\widetilde{\boldeta}_{1}/\|\widetilde{\boldeta}_{1}\|_{2},...,\widetilde{\boldeta}_{d}/\|\widetilde{\boldeta}_{d}\|_{2})=\vT\vM$, where $\vT$ is a $p\times d$ orthogonal matrix. Since $\widetilde{\boldeta}_{j}/\|\widetilde{\boldeta}_{j}\|_{2}'s$ are nearly mutually orthogonal,
we know that $\vM^{\tau}\vM$ is nearly an identity matrix.
Thus, by some continuity argument,  the statement is reduced to the following linear algebra lemma.
\begin{lem}
Let $\vA$ be a $p\times p$ positive definite matrix such that $C_{\min}\leq \lambda_{\min}(\vA)\leq \lambda_{\max}(\vA) \leq C_{\max}$ for some positive constants $C_{\min}$ and $C_{\max}$. There exists constant $0<\zeta<1$ such that for any $p\times d$ orthogonal matrix $\vB$,  we have 
\begin{align}
\left|\frac{\vB_{*,i}^{\tau}\vA^{\tau}\vA\vB_{*,j}}{\|\vA\vB_{*,i}\|_{2}\|\vA\vB_{*,j}\|_{2}}\right|\leq \zeta \quad \forall i \neq j.
\end{align}
\proof When $d$ is finite, without loss of generality, we can assume that $B$ is a $p\times 2$ matrix. Note that the expression on the left side is invariant under orthogonal transformation of $B$. We can simply assume that $B$ is a matrix with the last $p-2$-rows consisting of all zeros. The result follows immediately based on basic calculation.
\epf
\end{lem}

\subsection{Proof of Theorem \ref{thm:3}}
Recall that $\widehat{\boldeta}_{i}$ is the eigenvector associated with the $i$-th eigenvalue  $\widehat{\lambda}_{i}$ of $\widehat{\bLambda}_{H}$,
 $\widetilde{\boldeta}_{i}=P_{\bLambda}\widehat{\boldeta}_{i}$ and $\widetilde{\bbeta}_{i}=\bSigma^{-1}\widetilde{\boldeta}_{i}$, $i=1,\ldots,H$ (see e.g., \eqref{notation:inline}). 
The argument in Theorem 1 implies that, for any $1 \leq i \leq H$,
\begin{align}
\|\widehat{\bbeta}_{i}-\widetilde{\bbeta}_{i}\|_{2}\leq C\sqrt{ \frac{s\log(p)}{n\widehat{\lambda}}}.
\end{align}
The Proposition 1 implies that
\begin{align}
\|\widetilde{\bbeta}_{i}\|_{2}\geq C_{1}\sqrt{\frac{\lambda}{\widehat{\lambda}_{i}}}, 1\leq i \leq d \mbox{ and } 
\|\widetilde{\bbeta}_{i}\|_{2}\leq C_{2}\sqrt{\frac{\lambda}{\widehat{\lambda}_{i}}} \frac{\sqrt{p\log(p)}}{n\lambda}, d+1\leq i \leq H. 
\end{align}
The above two statements give us the desried result in Theorem \ref{thm:3}. \epf

\newpage

\begin{center}
{\LARGE \bf Supplementary Material} 
\end{center}

\hspace{.5in}

\section{\bf Proof of Theorem \ref{thm:consistency}}\label{sec:appendix:proof:consistency}
  Let $T=supp(\bbeta_{0})$. For a vector $\bgamma\in \bbR^{p}$, let $\bgamma_{T}$ and $\bgamma_{T^{c}}$ be the sub-vector consists of the components of $\bgamma$ in $T$ and $T^{c}$ respectively.
 We consider the following events sets
\begin{align*}
\mathsf{E}_{1}&=\left\lbrace~\omega\mid \|P_{\bLambda}\widehat{\boldeta}\|_{2}\geq \mathfrak{b}_{1}\sqrt{\frac{\lambda}{\widehat{\lambda}}} ~ \right\rbrace,\\
\mathsf{E}_{2}&= \left\lbrace~\omega\mid\frac{1}{n}\vX\vX^{\tau} \mbox{ satisfies the RE condition of order $s$} \mbox{ with parameter}  \left(3, \frac{\sqrt{C_{\min}}}{8}\right) ~\right\rbrace ,\\ 
\mathsf{E}_{3}&=\left\{~\omega \mid \|\frac{1}{n}\vX\vX^{\tau}\bdelta\|_{\infty}\leq \mu \mbox{ where } \mu=\mathfrak{b}_{2}\sqrt{\frac{\log(p)}{n\widehat{\lambda}}} ~\right\},\\
\mathsf{E}_{4}&=\left\lbrace~\omega\mid \|\bdelta_{T^{c}}\|_{1}\leq 3 \|\bdelta_{T}\|_{1} ~\right\rbrace
\end{align*}
where $\mathfrak{b}_{1}$ and $\mathfrak{b}_{2}$ are sufficiently large constants to be specified later.

Proposition \ref{prop:key} implies that $\mathsf{E}_{1}$ happens with high probability.
Proposition \ref{prop:second} implies that $\mathsf{E}_{2}$ happens with high probability.
Below, we will show that  $\mathsf{E}_{3}$ happens with high probability ( see Lemma \ref{lem:key1} below) and $\mathsf{E}_{4}$ happens with high probability (see  Lemma \ref{lem:key2} below).  We conclude that  the event $\mathsf{E}=\mathsf{E}_{1}\cap\mathsf{E}_{2}
\cap\mathsf{E}_{3}\cap\mathsf{E}_{4}$ happens with high probability.
Conditioning on $\mathsf{E}$, we have
\begin{align*}
\frac{1}{64}C_{\min}\|\bdelta\|_{2}^{2}\leq \frac{1}{n}\|\vX^{\tau}\bdelta\|^{2}\leq\frac{1}{n} \|\vX\vX^{\tau}\bdelta\|_{\infty}\|\bdelta\|_{1}\leq 4\mu\|\bdelta_{T}\|_{1}\leq 4\sqrt{s}\mu\|\bdelta_{T}\|_{2},
\end{align*}
i.e.,
$
\|\bdelta\|_{2}\leq \frac{256}{C_{\min}}\sqrt{s}\mu=C\sqrt{\frac{s\log(p)}{n\widehat{\lambda}}}.
$
Since  conditioning on $\mathsf{E}$, there exist constants $C_{1}$ and $C_{2}$ such that
$
C_{1}\sqrt{\frac{\lambda}{\widehat{\lambda}}}\leq \|\widetilde{\bbeta}\|_{2}\leq C_{2}
$,
we know that $\|P_{\widehat{\bbeta}}-P_{\bbeta_{0}}\|_{F}=\|P_{\widehat{\bbeta}}-P_{\widetilde{\bbeta}}\|_{F}\leq 4\frac{\|\widehat{\bbeta}-\widetilde{\bbeta}\|_{2}}{\|\widetilde{\bbeta}\|_{2}} \leq  C \sqrt{\frac{s\log(p)}{n\lambda}}$ holds with high probability. \epf

\begin{lem}\label{lem:key1}
Assume that conditions ${\bf A1)},{\bf A2)}$ and ${\bf A3)}$ hold. Let $\mu=A\sqrt{\frac{\log(p)}{n\widehat{\lambda}}}$. For sufficiently large $A$, we have that
\begin{align}
\|\frac{1}{n}\vX\vX^{\tau}\bdelta\|_{\infty} \leq  \mu
\end{align}
holds with high probability.
\proof
Since $\bdelta=\widehat{\bbeta}-\widetilde{\bbeta}$, $\widehat{\boldeta}=\frac{1}{n}\vX\widetilde{\vy}$ and $\bSigma\widetilde{\bbeta}=\widetilde{\boldeta}$,  we have
\begin{equation}\label{decomp:temp}
\begin{aligned}
\frac{1}{n}\|\vX\vX^{\tau}\bdelta\|_{\infty} \leq &
\underbrace{\frac{1}{n}\|\vX\vX^{\tau}\widehat{\bbeta}-\vX\widetilde{\vy}\|_{\infty}}_{I}+
\underbrace{\frac{1}{n}\|\vX\tilde{\vy}-n\widetilde{\boldeta}\|_{\infty}}_{II}
+
\underbrace{\|(\frac{1}{n}\vX\vX^{\tau}-\bSigma)\widetilde{\bbeta}\|_{\infty}}_{III}.
\end{aligned}
\end{equation}
\noindent{\it For I.}
By the definition of $\widehat{\bbeta}$, we have 
$0 \in \frac{1}{n}\left(\vX\vX^{\tau}\widehat{\bbeta}-\vX\widetilde{\vy}\right)+\mu sgn\left(\widehat{\bbeta}\right),
$
i.e., $\frac{1}{n}\|\vX\vX^{\tau}\widehat{\bbeta}-\vX\widetilde{\vy}\|_{\infty}\leq \mu$.

\vspace*{3mm}
\noindent{\it For II.}
Let $\vx=\vz+\vw$ be the orthogonal decomposition with respect to $col(var(\bbE[\vx|y]))$ and its orthogonal complement. 
Recall that we have introduced the decomposition:
\begin{align}\label{eqn:decomposition}
\vX_{H}=\vZ_{H}+\vW_{H}.
\end{align}
It is easy to see that $col(\vZ_{H})=col(\boldeta_{0})$ and  $\sqrt{c}~cov(\vw)^{-1/2}\vW_{H}$ is identically distributed to a matrix, $\mathcal{E}_{1}$, with all the entries are $i.i.d.$ standard normal random variables. Let
\begin{align*}
 \balpha_{1}=\frac{1}{\sqrt{H}}\vZ_{H}^{\tau}\widehat{\boldeta},~~ \balpha_{2}=\frac{1}{\sqrt{H}}\vW_{H}^{\tau}\widehat{\boldeta} \mbox{~~and~~} \balpha=\frac{1}{\sqrt{H}}\vX_{H}^{\tau}\widehat{\boldeta}=\balpha_{1}+\balpha_{2}. 
\end{align*}
Since 
$\frac{1}{H}\vX_{H}\vX^{\tau}_{H}\widehat{\boldeta}=
\widehat{\lambda}\widehat{\boldeta}$, we know that $
\|\balpha\|_{2}^{2}=\widehat{\lambda}
$ and 
$
\widehat{\boldeta}=\frac{1}{\sqrt{H}\widehat{\lambda}}\vZ_{H}\balpha+\frac{1}{\sqrt{H}\widehat{\lambda}}\vW_{H}\balpha.
$
From this, we know 
$
\widetilde{\boldeta}=P_{\eta_{0}}\widehat{\boldeta}=\frac{1}{\sqrt{H}\widehat{\lambda}}\vZ_{H}\balpha$ and $\widehat{\boldeta}-P_{\eta_{0}}\widehat{\boldeta}=\frac{1}{\sqrt{H}\widehat{\lambda}}\vW_{H}\balpha.
$
Since $\|\balpha\|_{1} \leq \sqrt{H}\|\balpha\|_{2}=\sqrt{H\widehat{\lambda}}$, we know that 
$
\|\widehat{\boldeta}-P_{\eta_{0}}\widehat{\boldeta}\|_{\infty}\leq \frac{1}{\sqrt{\widehat{\lambda}}}\|\vW_{H}\|_{\infty,\infty}.
$
It is easy to see that for positive constant $A_{2}(>1)$, we have
\begin{align*}
\bbP\left(\|\frac{cov(\vw)^{1/2}}{\sqrt{c}}\mathcal{E}_{1}\|_{\infty,\infty}>\lambda_{\max}^{1/2}(\bSigma)\sqrt{\frac{A_{2}H\log(pH)}{n}}\right)\leq 2\exp^{-(A_{2}-1)\log(pH)},
\end{align*}
{\it i.e.}, by letting $A>2\lambda_{\max}^{1/2}(\bSigma)\sqrt{A_{2}H}$, we have that
$
\|\widehat{\boldeta}-P_{\eta_{0}}\widehat{\boldeta}\|_{\infty}\leq A\sqrt{\frac{\log(p)}{n\widehat{\lambda}}}
$
holds with high probability.

\vspace*{3mm}

\noindent{\it For III.}
Let $\vx=\bSigma^{1/2}\bepsilon$ where $\bepsilon\sim N(0,\vI_{p})$, then $\vX=\bSigma^{1/2}\mathcal{E}$ where $\mathcal{E}$ is a $p\times n$ matrix  with $i.i.d$ standard normal entries.  
 Since $\widetilde{\bbeta}=\bSigma^{-1}\widetilde{\boldeta}$, $\widetilde{\boldeta}=\|\widetilde{\boldeta}\|\boldeta_{0}$ and $\|\widetilde{\boldeta}\|\leq 1$, we know that for any $0<t<1/2$, 
\begin{align}\label{eqn:key:lower}
\frac{1}{n}\|\left(\vX\vX^{\tau}-n\bSigma\right)\widetilde{\bbeta}\|_{\infty} =\|\widetilde{\boldeta}\| \|\bSigma^{1/2}\left(\frac{1}{n}\mathcal{E}\mathcal{E}^{\tau}-\vI_{p}\right)\bSigma^{-1/2}
\boldeta_{0}\|_{\infty} >t \lambda_{max}(\bSigma)^{1/2}\lambda_{\min}(\bSigma)^{-1/2}  
\end{align}
with probability at most $4\exp^{-\frac{3}{16}nt^{2}+\log(p)}$. 

In fact, it follows from $\|\bSigma^{1/2}_{i,*}\|\leq \lambda_{max}(\bSigma)^{1/2}$ and $\|\bSigma^{-1/2}\boldeta_{0}\|_{2}\leq \lambda_{\min}(\bSigma)^{-1/2}$ and  $
\bbP\left(\widetilde{\mathbf{E}}(\balpha_{1},\balpha_{2})\right)\leq 4\exp^{-\frac{3}{16}nt^{2}}
$ where
for  any two deterministic vectors $\balpha_{1}$ and $\balpha_{2}$,
\begin{align}
\widetilde{\mathbf{E}}(\balpha_{1},\balpha_{2})=\big\{~\omega~|~ |\balpha_{1}\left(\frac{1}{n}\mathcal{E}\mathcal{E}^{\tau}-\vI_{p}\right)\balpha_{2}| > t\|\balpha_{1}\|_{2}\|\balpha_{2}\|_{2}~\big\}.
\end{align}
Let $t=\sqrt{\frac{16A_{3}\log(p)}{3n}}$. Conditioning on $\mathsf{E}_{1}$ and  the events such that equation \eqref{eqn:key:lower} does not hold, we have
\begin{align}\label{eqn:inlince:temp1pp}
\frac{1}{n}\|\left(\vX\vX^{\tau}-n\bSigma\right)\widetilde{\bbeta}\|_{\infty} \leq  C\sqrt{\frac{\lambda}{\widehat{\lambda}}\frac{16A_{3}\log(p)}{3n}} \lambda_{max}(\bSigma)^{1/2}\lambda_{\min}(\bSigma)^{-1/2} 
\end{align}
with high probability.

\vspace*{4mm}

To summarize, we know  that,  for sufficiently large constant $A$,
$
\|\frac{1}{n}\vX\vX^{\tau}\delta\|_{\infty} \leq  A\sqrt{\frac{\log(p)}{n\widehat{\lambda}}}
$
holds with high probability. 
\epf
\end{lem}

\begin{lem}\label{lem:key2}
Let  $\mu=A\sqrt{\frac{\log(p)}{n\widehat{\lambda}}}$. For sufficiently large constant $A$,  we have 
$ 
\|\delta_{T^{c}}\|_{1} \leq 3 \|\delta_{T}\|_{1}
$
holds with high probability.
\proof
Since
$
\frac{1}{n}\|\vX\widetilde{\vy}-\vX\vX^{\tau}\widetilde{\bbeta}\|_{\infty}\leq \|\widehat{\boldeta}-\widetilde{\boldeta}\|_{\infty}+\|(\frac{1}{n}\vX\vX^{\tau}-\bSigma)\widetilde{\bbeta}\|_{\infty}=II+III
$
where $II$ and $III$ are introduced in \eqref{decomp:temp}.
Note that both $II$ and $III$ do not depend on the choice of $\mu$.  Thus there exists $A_{1}$ ( does not depend on the choice of $\mu$) such that, following the argument presented in Lemma \ref{lem:key1}, 
\begin{align}\|\frac{1}{n}\vX\widetilde{y}-\frac{1}{n}\vX\vX^{\tau}\widetilde{\bbeta}\|_{\infty}\leq A_{1}\sqrt{\frac{\log(p)}{n\widehat{\lambda}}}
\end{align}
holds with high probability.

Let us choose a sufficiently large $A>2A_{1}$ such that$ \|\frac{1}{n}\vX\vX^{\tau}\bdelta\|_{\infty} \leq  \mu$ ( by definition $\mu=A\sqrt{\frac{\log(p)}{n\widehat{\lambda}}}$ )  holds high probability. Since
$\mathcal{L}_{\widehat{\bbeta}}\leq \mathcal{L}_{\widetilde{\bbeta}}$, we have
\begin{align}
-\|\delta\|_{1}\|\frac{1}{n}\vX\widetilde{y}-\frac{1}{n}\vX\vX^{\tau}\widetilde{\bbeta}\|_{\infty}\leq  \mu(\|\widetilde{\bbeta}\|_{1}-\|\widehat{\bbeta}\|_{1} ).
\end{align}
Note that $\|\widetilde{\bbeta}\|_{1}-\|\widehat{\bbeta}\|_{1}\leq \|\delta_{T}\|_{1}-\|\delta_{T^{c}}\|_{1}$
and 
$\|\delta\|_{1}=\|\delta_{T}\|_{1}+\|\delta_{T^{c}}\|_{1}$.
Thus, we have
$
\|\delta_{T^{c}}\|_{1}\leq \frac{A+A_{1}}{A-A_{1}}\|\delta_{T}\|_{1}\leq 3 \|\delta_{T}\|_{1}
$ holds with high probability.\epf 
\end{lem}

\vspace*{4mm}

\section{Proof of Lemma \ref{lem:deviation:bound:diagonal}.}

\begin{lem}[Deviation]\label{lem:dev1}
Let $z_{j}=x_{j}^{2}-\sigma_{j}^{2}$, $j=1,....,p$ where $x_{j} \sim  N(0,\sigma^{2}_{j})$. Assume that $a\leq \sigma_{j}^{2} \leq b, j=1,...,p,$ for some positive constants $a,b$. Then for any $\alpha$ satisfying that $\alpha <\frac{1}{2}b^{-2}$, one has
\begin{align}
\bbP\left(\left|\frac{1}{p}\sum_{j}z_{j}\right|>\alpha\right) \leq 2\exp(-Cp\alpha^{2}) 
\end{align}
for some constant $C$.
\proof We have $\bbE[\exp^{tz_{j}}]=\exp\left(-\frac{1}{2}\log(1-2t\sigma^{2}_{j})-t\sigma_{j}^{2} \right)$. Thus, for any $t>0$, one has
\begin{align*}
\bbP\left(\frac{1}{p}\sum_{j}z_{j}>\alpha\right) \leq \exp\left(-\frac{1}{2}\sum_{j}\log(1-2t\sigma^{2}_{j})-t\sum_{j}\sigma_{j}^{2}-pt\alpha \right).
\end{align*}
Let us choose $t=\kappa\alpha$ where $\kappa<1$ is a 
constant to be determined later. Then, we have
\begin{align*}
\frac{1}{2}\sum_{j}\log(1-2t\sigma^{2}_{j})+t\sum_{j}\sigma_{j}^{2}+pt\alpha \geq pt\alpha-2t^{2}\sum_{j}\sigma_{j}^{4}=p\alpha^{2}(\kappa-2\kappa^{2}\sum_{j}\sigma_{j}^{4}/p).
\end{align*}
Thus, if we choose $\kappa$ such that $\kappa'=\kappa-\kappa^{2}\sum_{j}\sigma_{j}^{4}/p>0$, we have
\begin{align}
\bbP\left(\frac{1}{p}\sum_{j}z_{j}>\alpha\right) \leq \exp(-p\kappa'\alpha^{2}).
\end{align} 
Similar argument provides the bound for $\bbP\left(    \frac{1}{p}\sum_{j}z_{j}<-\alpha\right)$.

\end{lem}

\begin{lem}[Deviation II]\label{lem:dev2}
Let $z_{j}$ and $z_{j}'$  $\sim N(0,\sigma_{j}^{2})$ be independent copies for $j=1,...,p$. Assume that $a\leq \sigma_{j}^{2}\leq b$  for some positive constants $a$ and $b$. Then for any $\alpha$ satisfies that $\alpha<\frac{1}{2}b^{-2}$, one has
\begin{align*}
\bbP\left(\left|\frac{1}{p}\sum_{j}z_{j}z_{j}' \right|>\alpha \right)\leq 4\exp(-Cp\alpha^{2})
\end{align*}
for some constant $C$.
\proof
Let $w_{j}=\frac{1}{\sqrt{2}}(z_{j}+z_{j}')$ and  $w_{j}'=\frac{1}{\sqrt{2}}(z_{j}-z_{j}')$, then  $z_{j}z_{j}'=\frac{(w_{j}^{2}-\sigma_{j}^{2})-(w_{j}^{'2}-\sigma_{j}^{2})}{2}$. Lemma \ref{lem:dev1} implies the desired bound.
\end{lem}

\section{ Results of simulations}\label{appendix:simulation}

\begin{table}
    \caption{ \label{tab:sim1} Estimation error: $\bSigma=\bSigma_1$, and  $\rho=0$.}
    \begin{tabular}{|c|c|c|c|c|c|c|c|}
    \hline
    &  p & Lasso-SIR  & DT-SIR & Lasso & M-Lasso & Lasso-SIR(Known $d$) & $\hat{d}$ \\
    \hline
    \hline
    \multirow{4}{*}{I}       & 100  &  0.09 ( 0.02 )&  0.21 ( 0.04 )& 0.08 ( 0.01 )& 0.09 ( 0.01 )& 0.09 ( 0.01 )& 1 \\        & 1000  &  0.12 ( 0.02 )&  0.21 ( 0.04 )& 0.1 ( 0.02 )& 0.22 ( 0.02 )& 0.12 ( 0.02 )& 1 \\        & 2000  &  0.14 ( 0.02 )&  0.22 ( 0.05 )& 0.1 ( 0.01 )& 0.29 ( 0.03 )& 0.14 ( 0.02 )& 1 \\        & 4000  &  0.18 ( 0.09 )&  0.22 ( 0.09 )& 0.11 ( 0.02 )& 0.39 ( 0.09 )& 0.18 ( 0.03 )& 1 \\  
    \hline
    \multirow{4}{*}{II}     & 100  &  0.05 ( 0.01 )&  0.29 ( 0.05 )& 0.23 ( 0.04 )& 0.05 ( 0.01 )& 0.05 ( 0.01 )& 1 \\        & 1000  &  0.09 ( 0.01 )&  0.35 ( 0.06 )& 0.3 ( 0.04 )& 0.12 ( 0.02 )& 0.09 ( 0.01 )& 1 \\        & 2000  &  0.12 ( 0.02 )&  0.38 ( 0.08 )& 0.31 ( 0.04 )& 0.18 ( 0.03 )& 0.11 ( 0.02 )& 1 \\        & 4000  &  0.15 ( 0.03 )&  0.41 ( 0.08 )& 0.33 ( 0.04 )& 0.27 ( 0.06 )& 0.15 ( 0.03 )& 1 \\    
	\hline
    \multirow{4}{*}{III}    & 100  &  0.17 ( 0.03 )&  0.23 ( 0.06 )& 1.14 ( 0.27 )& 0.2 ( 0.03 )& 0.18 ( 0.04 )& 1 \\        & 1000  &  0.27 ( 0.21 )&  0.3 ( 0.22 )& 1.25 ( 0.21 )& 0.63 ( 0.17 )& 0.23 ( 0.04 )& 1.1 \\        & 2000  &  0.35 ( 0.29 )&  0.34 ( 0.3 )& 1.31 ( 0.17 )& 0.77 ( 0.23 )& 0.26 ( 0.05 )& 1.1 \\        & 4000  &  0.45 ( 0.42 )&  0.42 ( 0.39 )& 1.29 ( 0.18 )& 0.93 ( 0.32 )& 0.34 ( 0.14 )& 1.3 \\    
    
    \hline
    \multirow{4}{*}{IV}     & 100  &  0.35 ( 0.03 )&  0.79 ( 0.09 )& 0.41 ( 0.05 )& 0.98 ( 0.3 )& 0.35 ( 0.03 )& 1 \\        & 1000  &  0.59 ( 0.2 )&  0.96 ( 0.17 )& 0.61 ( 0.09 )& 0.84 ( 0.18 )& 0.56 ( 0.05 )& 1.1 \\        & 2000  &  0.72 ( 0.24 )&  1.02 ( 0.2 )& 0.67 ( 0.07 )& 1.01 ( 0.2 )& 0.64 ( 0.06 )& 1.2 \\        & 4000  &  0.95 ( 0.35 )&  1.14 ( 0.28 )& 0.71 ( 0.08 )& 1.23 ( 0.28 )& 0.82 ( 0.16 )& 1.4 \\  
    
    \hline
    \multirow{4}{*}{V}    & 100  &  0.1 ( 0.02 )&  0.18 ( 0.03 )& 0.55 ( 0.22 )& 0.11 ( 0.02 )& 0.09 ( 0.02 )& 1 \\        & 1000  &  0.12 ( 0.03 )&  0.19 ( 0.04 )& 0.69 ( 0.21 )& 0.3 ( 0.02 )& 0.13 ( 0.03 )& 1 \\        & 2000  &  0.15 ( 0.09 )&  0.19 ( 0.09 )& 0.72 ( 0.25 )& 0.37 ( 0.08 )& 0.15 ( 0.03 )& 1 \\        & 4000  &  0.18 ( 0.09 )&  0.2 ( 0.09 )& 0.74 ( 0.21 )& 0.47 ( 0.08 )& 0.18 ( 0.04 )& 1 \\   
    \hline
  \end{tabular}
\end{table}

\begin{table}
  \caption{ \label{tab:sim2} Estimation error: $\bSigma=\bSigma_1$ and $\rho=0.3$.   }
  \begin{tabular}{|c|c|c|c|c|c|c|c|}
    \hline
    &  p & Lasso-SIR  & DT-SIR & Lasso & M-Lasso &Lasso-SIR(Known $d$)& $\hat{d}$ \\
    \hline
    \hline
    \multirow{4}{*}{I}       & 100  &  0.1 ( 0.02 )&  0.44 ( 0.08 )& 0.09 ( 0.01 )& 0.12 ( 0.02 )& 0.1 ( 0.01 )& 1 \\        & 1000  &  0.15 ( 0.02 )&  0.5 ( 0.08 )& 0.12 ( 0.02 )& 0.24 ( 0.02 )& 0.16 ( 0.02 )& 1 \\        & 2000  &  0.18 ( 0.02 )&  0.5 ( 0.07 )& 0.13 ( 0.02 )& 0.31 ( 0.04 )& 0.18 ( 0.02 )& 1 \\        & 4000  &  0.21 ( 0.03 )&  0.49 ( 0.09 )& 0.14 ( 0.02 )& 0.42 ( 0.07 )& 0.22 ( 0.04 )& 1 \\  
    
    \hline
    \multirow{4}{*}{II}     & 100  &  0.06 ( 0.01 )&  0.46 ( 0.08 )& 0.22 ( 0.03 )& 0.08 ( 0.01 )& 0.06 ( 0.01 )& 1 \\        & 1000  &  0.11 ( 0.02 )&  0.55 ( 0.08 )& 0.28 ( 0.04 )& 0.14 ( 0.02 )& 0.11 ( 0.02 )& 1 \\        & 2000  &  0.14 ( 0.02 )&  0.55 ( 0.08 )& 0.3 ( 0.04 )& 0.2 ( 0.03 )& 0.14 ( 0.02 )& 1 \\        & 4000  &  0.19 ( 0.04 )&  0.58 ( 0.09 )& 0.32 ( 0.04 )& 0.32 ( 0.07 )& 0.19 ( 0.04 )& 1 \\    
    
    \hline
    \multirow{4}{*}{III}     & 100  &  0.19 ( 0.03 )&  0.5 ( 0.1 )& 1.18 ( 0.21 )& 0.24 ( 0.03 )& 0.19 ( 0.03 )& 1 \\        & 1000  &  0.29 ( 0.17 )&  0.6 ( 0.16 )& 1.3 ( 0.16 )& 0.58 ( 0.14 )& 0.25 ( 0.04 )& 1 \\        & 2000  &  0.35 ( 0.27 )&  0.63 ( 0.22 )& 1.32 ( 0.14 )& 0.73 ( 0.21 )& 0.29 ( 0.06 )& 1.1 \\        & 4000  &  0.57 ( 0.46 )&  0.75 ( 0.39 )& 1.33 ( 0.14 )& 0.98 ( 0.36 )& 0.4 ( 0.2 )& 1.4 \\    
    
    \hline
    \multirow{4}{*}{IV}    & 100  &  0.38 ( 0.04 )&  0.85 ( 0.1 )& 0.56 ( 0.08 )& 0.54 ( 0.14 )& 0.37 ( 0.04 )& 1 \\        & 1000  &  0.56 ( 0.13 )&  0.97 ( 0.13 )& 0.7 ( 0.08 )& 0.77 ( 0.11 )& 0.54 ( 0.04 )& 1 \\        & 2000  &  0.65 ( 0.21 )&  1.03 ( 0.19 )& 0.76 ( 0.11 )& 0.92 ( 0.19 )& 0.58 ( 0.05 )& 1.1 \\        & 4000  &  0.79 ( 0.3 )&  1.12 ( 0.25 )& 0.79 ( 0.09 )& 1.09 ( 0.25 )& 0.65 ( 0.05 )& 1.3 \\ 
    
    \hline
    \multirow{4}{*}{V}       & 100  &  0.1 ( 0.02 )&  0.47 ( 0.11 )& 0.48 ( 0.22 )& 0.14 ( 0.02 )& 0.1 ( 0.02 )& 1 \\        & 1000  &  0.14 ( 0.03 )&  0.55 ( 0.08 )& 0.6 ( 0.23 )& 0.35 ( 0.04 )& 0.15 ( 0.03 )& 1 \\        & 2000  &  0.17 ( 0.04 )&  0.56 ( 0.08 )& 0.66 ( 0.26 )& 0.44 ( 0.06 )& 0.18 ( 0.04 )& 1 \\        & 4000  &  0.3 ( 0.28 )&  0.6 ( 0.21 )& 0.72 ( 0.27 )& 0.66 ( 0.22 )& 0.25 ( 0.13 )& 1.1 \\    
    
    \hline
  \end{tabular}
\end{table}

\begin{table}
  \caption{ \label{tab:sim4} Estimation error: $\bSigma=\bSigma_1$, and $\rho=0.8$.   }
  \begin{tabular}{|c|c|c|c|c|c|c|c|}
    \hline
    &  p & Lasso-SIR  & DT-SIR & Lasso & M-Lasso &Lasso-SIR(Known $d$) & $\hat{d}$ \\
    \hline
    \hline
    \multirow{4}{*}{I}     & 100  &  0.18 ( 0.02 )&  1.34 ( 0.09 )& 0.16 ( 0.03 )& 1.01 ( 0.04 )& 0.18 ( 0.02 )& 1 \\        & 1000  &  0.24 ( 0.02 )&  1.38 ( 0.05 )& 0.22 ( 0.02 )& 0.79 ( 0.08 )& 0.24 ( 0.02 )& 1 \\        & 2000  &  0.27 ( 0.03 )&  1.39 ( 0.03 )& 0.23 ( 0.02 )& 0.53 ( 0.07 )& 0.27 ( 0.03 )& 1 \\        & 4000  &  0.32 ( 0.04 )&  1.39 ( 0.04 )& 0.25 ( 0.03 )& 0.45 ( 0.05 )& 0.32 ( 0.04 )& 1 \\   
    
    \hline
    \multirow{4}{*}{II}      & 100  &  0.1 ( 0.01 )&  1.34 ( 0.09 )& 0.33 ( 0.06 )& 1.17 ( 0.04 )& 0.11 ( 0.01 )& 1 \\        & 1000  &  0.16 ( 0.01 )&  1.39 ( 0.03 )& 0.55 ( 0.1 )& 1.08 ( 0.03 )& 0.16 ( 0.02 )& 1 \\        & 2000  &  0.19 ( 0.02 )&  1.39 ( 0.05 )& 0.71 ( 0.14 )& 0.92 ( 0.08 )& 0.19 ( 0.02 )& 1 \\        & 4000  &  0.23 ( 0.03 )&  1.4 ( 0.03 )& 0.92 ( 0.14 )& 0.54 ( 0.08 )& 0.23 ( 0.03 )& 1 \\  

    \hline
    \multirow{4}{*}{III}     & 100  &  0.28 ( 0.04 )&  1.34 ( 0.09 )& 1.26 ( 0.22 )& 1 ( 0.06 )& 0.28 ( 0.05 )& 1 \\        & 1000  &  0.45 ( 0.08 )&  1.38 ( 0.05 )& 1.29 ( 0.17 )& 0.92 ( 0.06 )& 0.45 ( 0.09 )& 1 \\        & 2000  &  0.54 ( 0.11 )&  1.39 ( 0.04 )& 1.3 ( 0.16 )& 0.84 ( 0.09 )& 0.54 ( 0.11 )& 1 \\        & 4000  &  0.76 ( 0.28 )&  1.43 ( 0.19 )& 1.29 ( 0.15 )& 0.89 ( 0.28 )& 0.68 ( 0.12 )& 1.1 \\

    \hline
    \multirow{4}{*}{IV}     & 100  &  0.74 ( 0.07 )&  1.4 ( 0.02 )& 1.21 ( 0.09 )& 0.91 ( 0.09 )& 0.72 ( 0.06 )& 1 \\        & 1000  &  0.75 ( 0.07 )&  1.41 ( 0.01 )& 1.23 ( 0.08 )& 0.88 ( 0.08 )& 0.76 ( 0.08 )& 1 \\        & 2000  &  0.79 ( 0.17 )&  1.44 ( 0.1 )& 1.26 ( 0.09 )& 0.94 ( 0.17 )& 0.75 ( 0.07 )& 1.1 \\        & 4000  &  0.93 ( 0.31 )&  1.52 ( 0.22 )& 1.27 ( 0.08 )& 1.09 ( 0.3 )& 0.76 ( 0.06 )& 1.4 \\

    \hline
    \multirow{4}{*}{V}   & 100  &  0.19 ( 0.04 )&  1.31 ( 0.14 )& 0.36 ( 0.14 )& 1.1 ( 0.39 )& 0.19 ( 0.03 )& 1 \\        & 1000  &  0.31 ( 0.1 )&  1.38 ( 0.07 )& 0.56 ( 0.2 )& 0.55 ( 0.12 )& 0.32 ( 0.13 )& 1 \\        & 2000  &  0.5 ( 0.34 )&  1.42 ( 0.15 )& 0.74 ( 0.27 )& 0.71 ( 0.29 )& 0.47 ( 0.26 )& 1.1 \\        & 4000  &  1.15 ( 0.65 )&  1.66 ( 0.36 )& 0.82 ( 0.22 )& 1.25 ( 0.55 )& 0.8 ( 0.38 )& 2.1 \\       \hline
  \end{tabular}
\end{table}

\begin{table}
  \caption{ \label{tab:sim5} Estimation error: $\bSigma=\bSigma_2$ and  $\rho=0.2$.   }
  \begin{tabular}{|c|c|c|c|c|c|c|c|}
    \hline
    &  p & Lasso-SIR  & DT-SIR & Lasso & M-Lasso &Lasso-SIR(Known $d$) & $\hat{d}$ \\
    \hline
    \hline
    \multirow{4}{*}{I}       & 100  &  0.13 ( 0.03 )&  1.22 ( 0.17 )& 0.09 ( 0.01 )& 0.16 ( 0.03 )& 0.13 ( 0.03 )& 1 \\        & 1000  &  0.33 ( 0.25 )&  1.37 ( 0.15 )& 0.11 ( 0.01 )& 0.65 ( 0.18 )& 0.27 ( 0.06 )& 1.1 \\        & 2000  &  0.3 ( 0.16 )&  1.37 ( 0.15 )& 0.12 ( 0.02 )& 0.74 ( 0.16 )& 0.3 ( 0.1 )& 1 \\        & 4000  &  0.36 ( 0.22 )&  1.38 ( 0.17 )& 0.13 ( 0.02 )& 0.81 ( 0.15 )& 0.3 ( 0.09 )& 1.1 \\  
    \hline
    \multirow{4}{*}{II}     & 100  &  0.1 ( 0.02 )&  1.26 ( 0.13 )& 0.24 ( 0.03 )& 0.11 ( 0.02 )& 0.1 ( 0.02 )& 1 \\        & 1000  &  0.25 ( 0.06 )&  1.37 ( 0.12 )& 0.31 ( 0.04 )& 0.47 ( 0.11 )& 0.25 ( 0.05 )& 1 \\        & 2000  &  0.3 ( 0.1 )&  1.38 ( 0.1 )& 0.33 ( 0.04 )& 0.59 ( 0.14 )& 0.29 ( 0.06 )& 1 \\        & 4000  &  0.31 ( 0.08 )&  1.4 ( 0.06 )& 0.35 ( 0.05 )& 0.65 ( 0.16 )& 0.32 ( 0.09 )& 1 \\  
    \hline
    \multirow{4}{*}{III}    & 100  &  0.24 ( 0.12 )&  1.24 ( 0.18 )& 1.19 ( 0.21 )& 0.33 ( 0.11 )& 0.23 ( 0.04 )& 1 \\        & 1000  &  0.55 ( 0.33 )&  1.33 ( 0.23 )& 1.3 ( 0.16 )& 0.98 ( 0.15 )& 0.4 ( 0.1 )& 1.3 \\        & 2000  &  0.59 ( 0.33 )&  1.35 ( 0.24 )& 1.28 ( 0.18 )& 1.08 ( 0.14 )& 0.45 ( 0.15 )& 1.3 \\        & 4000  &  0.58 ( 0.33 )&  1.36 ( 0.22 )& 1.28 ( 0.18 )& 1.14 ( 0.17 )& 0.47 ( 0.17 )& 1.3 \\   
    \hline
    \multirow{4}{*}{IV}    & 100  &  0.54 ( 0.05 )&  1.37 ( 0.06 )& 1.19 ( 0.1 )& 0.6 ( 0.06 )& 0.54 ( 0.05 )& 1 \\        & 1000  &  0.63 ( 0.05 )&  1.41 ( 0.01 )& 1.25 ( 0.1 )& 0.87 ( 0.05 )& 0.63 ( 0.05 )& 1 \\        & 2000  &  0.64 ( 0.05 )&  1.41 ( 0 )& 1.27 ( 0.1 )& 0.99 ( 0.05 )& 0.65 ( 0.05 )& 1 \\        & 4000  &  0.65 ( 0.06 )&  1.41 ( 0 )& 1.27 ( 0.09 )& 1.07 ( 0.04 )& 0.66 ( 0.05 )& 1 \\    
    \hline
    \multirow{4}{*}{V}     & 100  &  0.23 ( 0.29 )&  1.23 ( 0.25 )& 0.49 ( 0.19 )& 0.29 ( 0.27 )& 0.13 ( 0.03 )& 1.1 \\        & 1000  &  1.03 ( 0.24 )&  1.1 ( 0.36 )& 0.61 ( 0.23 )& 1.15 ( 0.17 )& 0.79 ( 0.37 )& 1.6 \\        & 2000  &  1.09 ( 0.2 )&  0.96 ( 0.43 )& 0.67 ( 0.19 )& 1.2 ( 0.16 )& 0.85 ( 0.4 )& 1.6 \\        & 4000  &  1.07 ( 0.27 )&  0.99 ( 0.47 )& 0.71 ( 0.22 )& 1.21 ( 0.16 )& 0.87 ( 0.41 )& 1.7 \\  
    \hline
  \end{tabular}
\end{table}

\begin{table}
  \caption{ \label{tab:sim6} Estimation error: $\bSigma=\bSigma_1$ and $\rho=0$.   }
  \begin{tabular}{|c|c|c|c|c|c|c|}
    \hline
    &  p & Lasso-SIR  & DT-SIR & M-Lasso &Lasso-SIR(Known $d$) &$\hat{d}$ \\
    \hline
    \hline
    \multirow{4}{*}{VI}     & 100  &  0.15 ( 0.03 )&  0.18 ( 0.06 )& 0.23 ( 0.05 )& 0.14 ( 0.04 )& 2 \\           & 1000  &  0.18 ( 0.06 )&  0.17 ( 0.07 )& 0.61 ( 0.03 )& 0.17 ( 0.05 )& 2 \\           & 2000  &  0.22 ( 0.13 )&  0.2 ( 0.14 )& 0.72 ( 0.1 )& 0.2 ( 0.07 )& 2 \\           & 4000  &  0.28 ( 0.13 )&  0.2 ( 0.1 )& 0.86 ( 0.09 )& 0.27 ( 0.14 )& 2 \\   
    \hline
    \multirow{4}{*}{VII}   & 100  &  0.27 ( 0.04 )&  0.35 ( 0.06 )& 0.32 ( 0.06 )& 0.27 ( 0.04 )& 2 \\           & 1000  &  0.37 ( 0.09 )&  0.4 ( 0.1 )& 0.93 ( 0.06 )& 0.38 ( 0.07 )& 2 \\           & 2000  &  0.44 ( 0.14 )&  0.41 ( 0.14 )& 1.09 ( 0.09 )& 0.44 ( 0.09 )& 2 \\           & 4000  &  0.75 ( 0.41 )&  0.64 ( 0.42 )& 1.4 ( 0.27 )& 0.59 ( 0.21 )& 2.4 \\   
    \hline
    \multirow{4}{*}{VIII}    & 100  &  0.87 ( 0.29 )&  0.88 ( 0.27 )& 0.9 ( 0.23 )& 0.23 ( 0.03 )& 1.2 \\           & 1000  &  0.45 ( 0.25 )&  0.44 ( 0.25 )& 0.91 ( 0.21 )& 0.31 ( 0.04 )& 1.8 \\           & 2000  &  0.34 ( 0.05 )&  0.35 ( 0.05 )& 0.8 ( 0.05 )& 0.34 ( 0.05 )& 2 \\           & 4000  &  0.57 ( 0.31 )&  0.53 ( 0.3 )& 1.04 ( 0.18 )& 0.41 ( 0.07 )& 1.9 \\      
    \hline
    \multirow{4}{*}{IX}       & 100  &  0.87 ( 0.28 )&  0.89 ( 0.27 )& 0.91 ( 0.22 )& 0.26 ( 0.08 )& 1.2 \\           & 1000  &  0.6 ( 0.3 )&  0.54 ( 0.33 )& 1.1 ( 0.06 )& 0.39 ( 0.14 )& 1.7 \\           & 2000  &  0.78 ( 0.3 )&  0.71 ( 0.36 )& 1.18 ( 0.14 )& 0.59 ( 0.29 )& 1.6 \\           & 4000  &  0.96 ( 0.26 )&  0.83 ( 0.33 )& 1.25 ( 0.19 )& 0.84 ( 0.36 )& 1.5 \\      
    \hline
    \hline
  \end{tabular}
\end{table}

\begin{table}
  \caption{ \label{tab:sim7} Estimation error: $\bSigma=\bSigma_1$ and $\rho=0.3$.   }
  \begin{tabular}{|c|c|c|c|c|c|c|}
    \hline
    &  p & Lasso-SIR  & DT-SIR & M-Lasso &Lasso-SIR(Known $d$) & $\hat{d}$ \\
    \hline
    \hline
    \multirow{4}{*}{VI}    & 100  &  0.2 ( 0.04 )&  0.34 ( 0.1 )& 0.25 ( 0.05 )& 0.2 ( 0.05 )& 2 \\           & 1000  &  0.24 ( 0.05 )&  0.3 ( 0.11 )& 0.61 ( 0.03 )& 0.23 ( 0.06 )& 2 \\           & 2000  &  0.26 ( 0.11 )&  0.36 ( 0.11 )& 0.71 ( 0.07 )& 0.26 ( 0.08 )& 2 \\           & 4000  &  0.31 ( 0.2 )&  0.41 ( 0.21 )& 0.84 ( 0.14 )& 0.29 ( 0.08 )& 2 \\      
        \hline
    \multirow{4}{*}{VII}    & 100  &  0.28 ( 0.04 )&  0.63 ( 0.09 )& 0.42 ( 0.04 )& 0.28 ( 0.04 )& 2 \\           & 1000  &  0.41 ( 0.12 )&  0.71 ( 0.1 )& 0.95 ( 0.09 )& 0.4 ( 0.12 )& 2 \\           & 2000  &  0.58 ( 0.27 )&  0.78 ( 0.2 )& 1.17 ( 0.18 )& 0.54 ( 0.19 )& 2.1 \\           & 4000  &  0.97 ( 0.41 )&  0.92 ( 0.3 )& 1.46 ( 0.25 )& 0.77 ( 0.31 )& 2.3 \\  
    \hline
    \multirow{4}{*}{VIII}    & 100  &  0.25 ( 0.07 )&  0.55 ( 0.09 )& 0.35 ( 0.09 )& 0.22 ( 0.03 )& 2 \\           & 1000  &  0.32 ( 0.08 )&  0.59 ( 0.13 )& 0.77 ( 0.17 )& 0.29 ( 0.04 )& 2 \\           & 2000  &  0.34 ( 0.14 )&  0.69 ( 0.12 )& 0.81 ( 0.12 )& 0.34 ( 0.06 )& 2 \\           & 4000  &  0.57 ( 0.35 )&  0.77 ( 0.26 )& 1.11 ( 0.24 )& 0.46 ( 0.2 )& 2.2 \\       
    \hline
    \multirow{4}{*}{IX}      & 100  &  0.31 ( 0.07 )&  0.5 ( 0.08 )& 0.43 ( 0.07 )& 0.31 ( 0.07 )& 2 \\           & 1000  &  0.35 ( 0.11 )&  0.47 ( 0.09 )& 0.99 ( 0.05 )& 0.36 ( 0.09 )& 2 \\           & 2000  &  0.42 ( 0.22 )&  0.55 ( 0.2 )& 1.17 ( 0.14 )& 0.4 ( 0.12 )& 2.1 \\           & 4000  &  0.51 ( 0.24 )&  0.56 ( 0.21 )& 1.28 ( 0.11 )& 0.44 ( 0.13 )& 2 \\         
    \hline
    \hline
  \end{tabular}
\end{table}

\begin{table}
  \caption{\label{tab:sim9} Estimation error: $\bSigma=\bSigma_1$ and $\rho=0.8$.   }
  \begin{tabular}{|c|c|c|c|c|c|c|}
    \hline
    &  p & Lasso-SIR  & DT-SIR & M-Lasso &Lasso-SIR(Known $d$) &$\hat{d}$ \\
    \hline
    \hline
    \multirow{4}{*}{VI}     & 100  &  0.52 ( 0.12 )&  1.86 ( 0.13 )& 1.01 ( 0.07 )& 0.51 ( 0.12 )& 2 \\           & 1000  &  0.79 ( 0.11 )&  1.92 ( 0.09 )& 0.93 ( 0.08 )& 0.79 ( 0.12 )& 2 \\           & 2000  &  0.96 ( 0.2 )&  1.94 ( 0.1 )& 1.05 ( 0.17 )& 0.94 ( 0.14 )& 2 \\           & 4000  &  1.14 ( 0.3 )&  2.01 ( 0.18 )& 1.26 ( 0.29 )& 1.06 ( 0.17 )& 2.2 \\          
    \hline
    \multirow{4}{*}{VII}    & 100  &  0.8 ( 0.34 )&  1.77 ( 0.12 )& 1.07 ( 0.24 )& 0.72 ( 0.36 )& 1.6 \\           & 1000  &  1.09 ( 0.2 )&  1.78 ( 0.13 )& 1.23 ( 0.19 )& 1.33 ( 0.21 )& 1.3 \\           & 2000  &  1.09 ( 0.15 )&  1.76 ( 0.12 )& 1.23 ( 0.18 )& 1.38 ( 0.15 )& 1.3 \\           & 4000  &  1.12 ( 0.23 )&  1.76 ( 0.15 )& 1.27 ( 0.23 )& 1.42 ( 0.07 )& 1.2 \\    
    \hline
    \multirow{4}{*}{VIII}     & 100  &  0.42 ( 0.18 )&  1.81 ( 0.14 )& 0.79 ( 0.33 )& 0.34 ( 0.04 )& 2 \\           & 1000  &  1 ( 0.38 )&  1.97 ( 0.14 )& 1.22 ( 0.32 )& 0.86 ( 0.4 )& 2.2 \\           & 2000  &  1.12 ( 0.35 )&  1.93 ( 0.17 )& 1.27 ( 0.29 )& 1.17 ( 0.33 )& 2.1 \\           & 4000  &  1.16 ( 0.28 )&  1.89 ( 0.17 )& 1.27 ( 0.24 )& 1.29 ( 0.24 )& 1.8 \\       
    \hline
    \multirow{4}{*}{IX}   & 100  &  0.78 ( 0.1 )&  1.9 ( 0.09 )& 0.95 ( 0.1 )& 0.79 ( 0.1 )& 2 \\           & 1000  &  0.92 ( 0.13 )&  1.95 ( 0.06 )& 1.08 ( 0.13 )& 0.9 ( 0.09 )& 2 \\           & 2000  &  0.97 ( 0.11 )&  1.97 ( 0.07 )& 1.17 ( 0.08 )& 0.95 ( 0.1 )& 2 \\           & 4000  &  1.12 ( 0.26 )&  2.03 ( 0.12 )& 1.37 ( 0.22 )& 1.01 ( 0.08 )& 2.3 \\ 
    \hline
    \hline
  \end{tabular}
\end{table}

\begin{table}
  \caption{\label{tab:sim10} Estimation error: $\bSigma=\bSigma_2$ and $\rho=0.2$.   }
  \begin{tabular}{|c|c|c|c|c|c|c|}
    \hline
    &  p & Lasso-SIR  & DT-SIR & M-Lasso &Lasso-SIR(Known $d$) &$\hat{d}$ \\
    \hline
    \hline
    \multirow{4}{*}{VI}     & 100  &  0.27 ( 0.21 )&  1.73 ( 0.17 )& 0.43 ( 0.17 )& 0.22 ( 0.04 )& 1.9 \\           & 1000  &  1.01 ( 0.01 )&  1.73 ( 0 )& 1.11 ( 0.03 )& 0.26 ( 0.06 )& 1 \\           & 2000  &  1.01 ( 0.01 )&  1.73 ( 0 )& 1.14 ( 0.04 )& 0.29 ( 0.08 )& 1 \\           & 4000  &  1.02 ( 0.01 )&  1.73 ( 0 )& 1.18 ( 0.04 )& 0.38 ( 0.19 )& 1 \\   
    \hline
    \multirow{4}{*}{VII}      & 100  &  0.39 ( 0.24 )&  1.7 ( 0.18 )& 0.55 ( 0.19 )& 0.31 ( 0.05 )& 1.9 \\           & 1000  &  1.03 ( 0.06 )&  1.73 ( 0.01 )& 1.25 ( 0.05 )& 0.47 ( 0.19 )& 1 \\           & 2000  &  1.04 ( 0.01 )&  1.73 ( 0.01 )& 1.3 ( 0.05 )& 0.55 ( 0.24 )& 1 \\           & 4000  &  1.04 ( 0.02 )&  1.73 ( 0 )& 1.34 ( 0.06 )& 0.69 ( 0.3 )& 1 \\  
    \hline
    \multirow{4}{*}{VIII}     & 100  &  0.24 ( 0.03 )&  1.69 ( 0.17 )& 0.34 ( 0.04 )& 0.24 ( 0.03 )& 2 \\           & 1000  &  0.97 ( 0.21 )&  1.73 ( 0.09 )& 1.15 ( 0.11 )& 0.33 ( 0.05 )& 1.1 \\           & 2000  &  1.03 ( 0.08 )&  1.74 ( 0.05 )& 1.24 ( 0.06 )& 0.35 ( 0.06 )& 1 \\           & 4000  &  1.03 ( 0.07 )&  1.74 ( 0.03 )& 1.26 ( 0.06 )& 0.41 ( 0.12 )& 1 \\       
    \hline
    \multirow{4}{*}{IX}      & 100  &  1 ( 0.12 )&  1.69 ( 0.06 )& 1.04 ( 0.1 )& 0.41 ( 0.07 )& 1 \\           & 1000  &  1.03 ( 0.01 )&  1.73 ( 0 )& 1.22 ( 0.04 )& 0.67 ( 0.2 )& 1 \\           & 2000  &  1.03 ( 0.01 )&  1.73 ( 0 )& 1.27 ( 0.03 )& 0.73 ( 0.21 )& 1 \\           & 4000  &  1.04 ( 0.03 )&  1.73 ( 0 )& 1.3 ( 0.04 )& 0.89 ( 0.25 )& 1 \\   
    \hline
    \hline
  \end{tabular}
\end{table}

\begin{table}
  \caption{\label{tab:sim11} Estimation error: $\bSigma=\bSigma_1$ and $\rho=0$.   }
  \begin{tabular}{|c|c|c|c|c|c|}
    \hline
    &  p & Lasso-SIR  & DT-SIR & M-Lasso & Lasso \\
    \hline
    \hline
    \multirow{4}{*}{X}   & 100  &  0.18 ( 0.03 ) &  0.54 ( 0.04 )& 0.21 ( 0.05 )& 0.18 ( 0.03 ) \\          & 1000  &  0.23 ( 0.03 ) &  1.13 ( 0.01 )& 0.6 ( 0.03 )& 0.25 ( 0.04 ) \\          & 2000  &  0.23 ( 0.04 ) &  1.24 ( 0.01 )& 0.67 ( 0.02 )& 0.27 ( 0.04 ) \\          & 4000  &  0.24 ( 0.03 ) &  1.29 ( 0.01 )& 0.71 ( 0.03 )& 0.28 ( 0.04 ) \\    
        \hline
    \multirow{4}{*}{XI}    & 100  &  0.33 ( 0.09 ) &  0.81 ( 0.05 )& 0.4 ( 0.13 )& 0.34 ( 0.08 ) \\          & 1000  &  0.41 ( 0.1 ) &  1.29 ( 0.01 )& 1.16 ( 0.03 )& 0.44 ( 0.1 ) \\          & 2000  &  0.43 ( 0.1 ) &  1.34 ( 0.01 )& 1.21 ( 0.03 )& 0.45 ( 0.11 ) \\          & 4000  &  0.45 ( 0.11 ) &  1.37 ( 0.01 )& 1.23 ( 0.03 )& 0.48 ( 0.1 ) \\   
        \hline
    \multirow{4}{*}{XII}    & 100  &  0.23 ( 0.03 ) &  0.53 ( 0.04 )& 0.26 ( 0.03 )& 0.2 ( 0.02 ) \\          & 1000  &  0.3 ( 0.03 ) &  1.12 ( 0.01 )& 0.62 ( 0.02 )& 0.3 ( 0.03 ) \\          & 2000  &  0.32 ( 0.03 ) &  1.23 ( 0.01 )& 0.7 ( 0.03 )& 0.33 ( 0.03 ) \\          & 4000  &  0.33 ( 0.03 ) &  1.29 ( 0.01 )& 0.75 ( 0.03 )& 0.36 ( 0.03 ) \\    
        \hline
    \multirow{4}{*}{XIII}      & 100  &  0.36 ( 0.07 ) &  0.92 ( 0.04 )& 0.44 ( 0.1 )& 1.05 ( 0.02 ) \\          & 1000  &  0.44 ( 0.08 ) &  1.69 ( 0.01 )& 1.21 ( 0.04 )& 1.08 ( 0.02 ) \\          & 2000  &  0.44 ( 0.08 ) &  1.8 ( 0.01 )& 1.3 ( 0.03 )& 1.1 ( 0.03 ) \\          & 4000  &  0.45 ( 0.08 ) &  1.85 ( 0.01 )& 1.34 ( 0.02 )& 1.11 ( 0.03 ) \\     
        \hline
    \hline
  \end{tabular}
\end{table}

\begin{table}
  \caption{ \label{tab:sim12} Estimation error: $\bSigma=\bSigma_1$ and  $\rho=0.3$.   }
  \begin{tabular}{|c|c|c|c|c|c|}
    \hline
    &  p & Lasso-SIR  & DT-SIR & M-Lasso & Lasso\\
    \hline
    \hline
    \multirow{4}{*}{X}   & 100  &  0.2 ( 0.03 ) &  0.59 ( 0.04 )& 0.26 ( 0.03 )& 0.19 ( 0.03 ) \\          & 1000  &  0.24 ( 0.03 ) &  1.15 ( 0.02 )& 0.56 ( 0.03 )& 0.26 ( 0.03 ) \\          & 2000  &  0.24 ( 0.03 ) &  1.25 ( 0.01 )& 0.63 ( 0.03 )& 0.27 ( 0.03 ) \\          & 4000  &  0.25 ( 0.03 ) &  1.31 ( 0.01 )& 0.69 ( 0.03 )& 0.29 ( 0.04 ) \\   
    
    \hline
    \multirow{4}{*}{XI}      & 100  &  0.33 ( 0.08 ) &  0.86 ( 0.05 )& 0.58 ( 0.15 )& 0.34 ( 0.08 ) \\          & 1000  &  0.41 ( 0.1 ) &  1.31 ( 0.01 )& 1.12 ( 0.04 )& 0.43 ( 0.09 ) \\          & 2000  &  0.41 ( 0.1 ) &  1.35 ( 0.01 )& 1.18 ( 0.04 )& 0.43 ( 0.1 ) \\          & 4000  &  0.45 ( 0.11 ) &  1.38 ( 0.01 )& 1.22 ( 0.04 )& 0.47 ( 0.12 ) \\    
    \hline
    \multirow{4}{*}{XII}     & 100  &  0.22 ( 0.03 ) &  0.53 ( 0.04 )& 0.3 ( 0.13 )& 0.2 ( 0.02 ) \\          & 1000  &  0.29 ( 0.03 ) &  1.1 ( 0.02 )& 0.58 ( 0.02 )& 0.3 ( 0.03 ) \\          & 2000  &  0.31 ( 0.04 ) &  1.22 ( 0.02 )& 0.66 ( 0.02 )& 0.33 ( 0.03 ) \\          & 4000  &  0.33 ( 0.03 ) &  1.29 ( 0.02 )& 0.72 ( 0.03 )& 0.36 ( 0.03 ) \\  
    \hline
    \multirow{4}{*}{XIII}    & 100  &  0.38 ( 0.07 ) &  1 ( 0.05 )& 0.6 ( 0.08 )& 1.06 ( 0.02 ) \\          & 1000  &  0.39 ( 0.07 ) &  1.73 ( 0.01 )& 1.17 ( 0.04 )& 1.08 ( 0.02 ) \\          & 2000  &  0.39 ( 0.06 ) &  1.84 ( 0.01 )& 1.29 ( 0.03 )& 1.09 ( 0.03 ) \\          & 4000  &  0.42 ( 0.08 ) &  1.88 ( 0.01 )& 1.34 ( 0.03 )& 1.11 ( 0.03 ) \\   
    \hline
    \hline
  \end{tabular}
\end{table}

\begin{table}
  \caption{ \label{tab:sim14} Estimation error: $\bSigma=\bSigma_1$ and $\rho=0.8$.   }
  \begin{tabular}{|c|c|c|c|c|c|}
    \hline
    &  p & Lasso-SIR  & DT-SIR & M-Lasso & Lasso \\
    \hline
    \hline
    \multirow{4}{*}{X}    & 100  &  0.27 ( 0.06 ) &  1.37 ( 0.04 )& 1 ( 0.06 )& 0.26 ( 0.04 ) \\          & 1000  &  0.46 ( 0.08 ) &  1.41 ( 0 )& 0.91 ( 0.06 )& 0.39 ( 0.06 ) \\          & 2000  &  0.53 ( 0.09 ) &  1.41 ( 0.01 )& 0.82 ( 0.09 )& 0.45 ( 0.07 ) \\          & 4000  &  0.64 ( 0.13 ) &  1.41 ( 0 )& 0.76 ( 0.08 )& 0.53 ( 0.1 ) \\     
    \hline
    \multirow{4}{*}{XI}    & 100  &  0.39 ( 0.1 ) &  1.38 ( 0.04 )& 0.81 ( 0.32 )& 0.39 ( 0.09 ) \\          & 1000  &  0.7 ( 0.17 ) &  1.41 ( 0 )& 1.03 ( 0.12 )& 0.7 ( 0.18 ) \\          & 2000  &  0.88 ( 0.17 ) &  1.41 ( 0.01 )& 1.1 ( 0.08 )& 0.86 ( 0.17 ) \\          & 4000  &  1.01 ( 0.17 ) &  1.41 ( 0 )& 1.17 ( 0.08 )& 1.01 ( 0.17 ) \\    
    \hline
    \multirow{4}{*}{XII}     & 100  &  0.36 ( 0.06 ) &  1.37 ( 0.06 )& 1.18 ( 0.04 )& 0.31 ( 0.05 ) \\          & 1000  &  0.63 ( 0.11 ) &  1.41 ( 0.02 )& 1.17 ( 0.04 )& 0.52 ( 0.08 ) \\          & 2000  &  0.8 ( 0.14 ) &  1.41 ( 0 )& 1.14 ( 0.04 )& 0.66 ( 0.11 ) \\          & 4000  &  1.03 ( 0.09 ) &  1.41 ( 0 )& 1.1 ( 0.05 )& 0.92 ( 0.11 ) \\     
    \hline
    \multirow{4}{*}{XIII}    & 100  &  0.51 ( 0.12 ) &  1.93 ( 0.05 )& 0.88 ( 0.25 )& 1.11 ( 0.04 ) \\          & 1000  &  0.55 ( 0.09 ) &  1.99 ( 0.02 )& 1.11 ( 0.09 )& 1.12 ( 0.04 ) \\          & 2000  &  0.56 ( 0.11 ) &  2 ( 0.01 )& 1.2 ( 0.08 )& 1.14 ( 0.04 ) \\          & 4000  &  0.61 ( 0.12 ) &  2 ( 0 )& 1.3 ( 0.05 )& 1.15 ( 0.04 ) \\  
    \hline
    \hline
  \end{tabular}
\end{table}

\begin{table}
  \caption{ \label{tab:sim15} Estimation error: $\bSigma=\bSigma_2$ and $\rho=0.2$.   }
  \begin{tabular}{|c|c|c|c|c|c|}
    \hline
    &  p & Lasso-SIR  & DT-SIR & M-Lasso & Lasso  \\
    \hline
    \hline
    \multirow{4}{*}{X}       & 100  &  0.19 ( 0.03 ) &  1.23 ( 0.16 )& 0.28 ( 0.04 )& 0.19 ( 0.03 ) \\          & 1000  &  0.24 ( 0.04 ) &  1.41 ( 0.01 )& 0.61 ( 0.03 )& 0.26 ( 0.04 ) \\          & 2000  &  0.24 ( 0.04 ) &  1.41 ( 0 )& 0.68 ( 0.03 )& 0.27 ( 0.04 ) \\          & 4000  &  0.26 ( 0.04 ) &  1.41 ( 0 )& 0.73 ( 0.03 )& 0.3 ( 0.04 ) \\    
    \hline
    \multirow{4}{*}{XI}    & 100  &  0.34 ( 0.09 ) &  1.27 ( 0.16 )& 0.59 ( 0.17 )& 0.35 ( 0.09 ) \\          & 1000  &  0.42 ( 0.09 ) &  1.41 ( 0 )& 1.16 ( 0.04 )& 0.44 ( 0.11 ) \\          & 2000  &  0.44 ( 0.12 ) &  1.41 ( 0 )& 1.21 ( 0.03 )& 0.46 ( 0.11 ) \\          & 4000  &  0.47 ( 0.12 ) &  1.41 ( 0 )& 1.24 ( 0.03 )& 0.49 ( 0.12 ) \\  
    \hline
    \multirow{4}{*}{XII}   & 100  &  0.25 ( 0.04 ) &  1.24 ( 0.16 )& 0.3 ( 0.04 )& 0.22 ( 0.02 ) \\          & 1000  &  0.32 ( 0.04 ) &  1.4 ( 0.01 )& 0.62 ( 0.03 )& 0.32 ( 0.03 ) \\          & 2000  &  0.33 ( 0.04 ) &  1.41 ( 0 )& 0.71 ( 0.03 )& 0.35 ( 0.04 ) \\          & 4000  &  0.36 ( 0.05 ) &  1.41 ( 0 )& 0.78 ( 0.03 )& 0.38 ( 0.04 ) \\    
    \hline
    \multirow{4}{*}{XIII}    & 100  &  0.4 ( 0.07 ) &  1.75 ( 0.16 )& 0.68 ( 0.1 )& 1.06 ( 0.02 ) \\          & 1000  &  0.43 ( 0.08 ) &  1.99 ( 0.02 )& 1.21 ( 0.04 )& 1.09 ( 0.03 ) \\          & 2000  &  0.47 ( 0.09 ) &  2 ( 0 )& 1.31 ( 0.03 )& 1.1 ( 0.03 ) \\          & 4000  &  0.48 ( 0.1 ) &  2 ( 0 )& 1.37 ( 0.03 )& 1.11 ( 0.03 ) \\   
    \hline
    \hline
  \end{tabular}
\end{table}

\begin{table}
\centering
  \caption{ Compare Lasso-SIR and Sparse SIR: $\bSigma=\bSigma_1$, $\rho=0$, $n=1,000$. \label{tab:SDR:1}  }
  \begin{tabular}{|c|c|c|c||c|c|c|}
    \hline
    Setting &  p & Lasso-SIR  & Sparse SIR & Setting & Lasso-SIR & Sparse SIR\\
    \hline
    \hline
    \multirow{3}{*}{I}     & 100  &  0.09 ( 0.02 )&  0.16 ( 0.013 )   & \multirow{3}{*}{II}  &  0.05 ( 0.01 )&  0.06 ( 0.006 ) \\        
        & 1000  &  0.12 ( 0.02 )&  1.41 ( 0.001 )   &   &  0.09 ( 0.01 )&  1.41 ( 0.001 )  \\        
        & 2000  &  0.14 ( 0.02 )&  1.41 ( 0 )   &   &  0.12 ( 0.02 )&  1.41 ( 0 ) \\
    \hline
    \hline
    \multirow{3}{*}{III}   & 100  &  0.17 ( 0.03 )&  0.35 ( 0.027 ) &     \multirow{3}{*}{IV}   &  0.35 ( 0.03 )&  0.39 ( 0.031 ) \\        
    & 1000  &  0.27 ( 0.21 )&  1.43 ( 0.1 )     &   &  0.59 ( 0.2 )&  1.44 ( 0.113 ) \\        
    & 2000  &  0.35 ( 0.29 )&  1.45 ( 0.146 ) &    &  0.72 ( 0.24 )&  1.46 ( 0.14 ) \\    
    \hline
    \hline
    \multirow{3}{*}{V}     & 100  &  0.1 ( 0.02 )&  0.22 ( 0.016 ) &    \multirow{3}{*}{VI}    &  0.15 ( 0.03 )&  0.39 ( 0.021 ) \\           
    & 1000  &  0.12 ( 0.03 )&  0.22 ( 0.016 )     &   &  0.18 ( 0.06 )&  2 ( 0.002 ) \\                   
    & 2000  &  0.15 ( 0.09 )&  1.42 ( 0.032 )     &   &  0.22 ( 0.13 )&  2 ( 0.045 ) \\           
    \hline
    \hline
    \multirow{3}{*}{VII}    & 100  &  0.15 ( 0.03 )&  0.59 ( 0.033 )     & \multirow{3}{*}{VIII}   &  0.87 ( 0.28 )&  0.45 ( 0.024 ) \\           

& 1000  &  0.18 ( 0.06 )&  2 ( 0.002 )     &   &  0.6 ( 0.3 )&  1.92 ( 0.128 ) \\                     
    & 2000  &  0.22 ( 0.13 )&  2 ( 0.045 )    &   &  0.78 ( 0.3 )&  1.88 ( 0.161 ) \\      
    \hline
    \hline
    \multirow{3}{*}{VIV}    & 100  &  0.87 ( 0.28 )&  0.96 ( 0.149 ) \\           & 1000  &  0.6 ( 0.3 )&  1.92 ( 0.128 ) \\           & 2000  &  0.78 ( 0.3 )&  1.88 ( 0.161 ) \\      
    \cline{1-4}
  \end{tabular}
\end{table}

\begin{table}
\centering
  \caption{ Compare Lasso-SIR and Sparse SIR: $\bSigma=\bSigma_1$, $\rho=0.3$, $n=1,000$. \label{tab:SDR:2}  }
  \begin{tabular}{|c|c|c|c||c|c|c|}
    \hline
    Setting &  p & Lasso-SIR  & Sparse SIR & Setting & Lasso-SIR & Sparse SIR\\
    \hline
    \hline
    \multirow{3}{*}{I}     & 100  &  0.1 ( 0.02 )&  0.18 ( 0.015 ) & \multirow{3}{*}{II}   &  0.06 ( 0.01 )&  0.07 ( 0.006 )       \\
            & 1000  &  0.15 ( 0.02 )&  1.41 ( 0.001 )     &   &  0.11 ( 0.02 )&  1.41 ( 0.001 ) \\ 
               & 2000  &  0.18 ( 0.02 )&  1.41 ( 0.001 )      &   &  0.14 ( 0.02 )&  0 ( 0 ) \\  
    \hline
    \hline
    \multirow{3}{*}{III}    & 100  &  0.19 ( 0.03 )&  0.36 ( 0.031 )  &   \multirow{3}{*}{IV}      &  0.38 ( 0.04 )&  0.45 ( 0.042 ) \\
      & 1000  &  0.29 ( 0.17 )&  1.43 ( 0.082 ) &  &  0.56 ( 0.13 )&  1.43 ( 0.062 ) \\
        & 2000  &  0.35 ( 0.27 )&  1.41 ( 0.001 )  &  &  0.65 ( 0.21 )&  1.45 ( 0.12 ) \\
    \hline
    \hline
    \multirow{3}{*}{V}     & 100  &  0.1 ( 0.02 )&  0.24 ( 0.02 ) &
    \multirow{3}{*}{VI}  &  0.19 ( 0.04 )&  0.52 ( 0.037 ) \\
       & 1000  &  0.14 ( 0.03 )&  1.41 ( 0.001 ) &  &  0.24 ( 0.05 )&  2 ( 0.001 ) \\ 
       & 2000  &  0.17 ( 0.04 )&  1.41 ( 0.002 )  &  &  0.26 ( 0.11 )&  2 ( 0.024 ) \\   
    \hline
    \hline
    \multirow{3}{*}{VII}    & 100  &  0.28 ( 0.04 )&  0.61 ( 0.039 ) &
    \multirow{3}{*}{VIII} &  0.22 ( 0.03 )&  0.45 ( 0.027 ) \\ 
       & 1000  &  0.41 ( 0.12 )&  2 ( 0.033 ) &  &  0.3 ( 0.04 )&  2 ( 0.001 ) \\
            & 2000  &  0.58 ( 0.29 )&  2.03 ( 0.091 ) &   &  0.34 ( 0.14 )&  2.01 ( 0.05 ) \\ 
    \hline
    \hline
    \multirow{3}{*}{VIV}    & 100  &  0.31 ( 0.07 )&  0.77 ( 0.042 ) \\ 
          & 1000  &  0.35 ( 0.1 )&  2 ( 0.001 ) \\         
          & 2000  &  0.43 ( 0.24 )&  2.02 ( 0.08 ) \\     
    \cline{1-4}
  \end{tabular}
\end{table}

\begin{table}
\centering
  \caption{Compare Lasso-SIR and Sparse SIR: $\bSigma=\bSigma_1$, $\rho=0.5$, and $n=1,000$. \label{tab:SDR:3}  }
   \begin{tabular}{|c|c|c|c||c|c|c|}
    \hline
  Setting  &  p & Lasso-SIR  & Sparse SIR &Setting & Lasso-SIR & Sparse SIR\\
    \hline
    \hline
    \multirow{3}{*}{I}      & 100  &  0.12 ( 0.02 )&  0.21 ( 0.019 ) &    \multirow{3}{*}{II}   &  0.07 ( 0.01 )&  0.08 ( 0.008 ) \\                
    & 1000  &  0.18 ( 0.02 )&  1.41 ( 0.001 )     &   &  0.12 ( 0.02 )&  1.41 ( 0.001 ) \\           
    & 2000  &  0.2 ( 0.02 )&  1.41 ( 0 )     &   &  0.15 ( 0.02 )&  1.41 ( 0.001 ) \\
    \hline
    \hline
    \multirow{3}{*}{III}      & 100  &  0.21 ( 0.03 )&  0.4 ( 0.037 ) &    \multirow{3}{*}{IV}   &  0.46 ( 0.05 )&  0.55 ( 0.05 ) \\                
    & 1000  &  0.28 ( 0.04 )&  1.41 ( 0.002 )     &   &  0.62 ( 0.22 )&  1.44 ( 0.146 ) \\                
    & 2000  &  0.35 ( 0.17 )&  1.43 ( 0.07 )     &   &  0.71 ( 0.34 )&  1.49 ( 0.217 ) \\    
    \hline
    \hline
    \multirow{3}{*}{V}   & 100  &  0.12 ( 0.02 )&  0.28 ( 0.024 )  &    \multirow{3}{*}{VI}        &  0.26 ( 0.06 )&  0.68 ( 0.047 ) \\ 
      & 1000  &  0.2 ( 0.03 )&  1.41 ( 0.001 ) &  &  0.33 ( 0.07 )&  2 ( 0.001 ) \\
       & 2000  &  0.38 ( 0.34 )&  1.47 ( 0.145 ) &  &  0.36 ( 0.11 )&  2 ( 0.024 ) \\ 
    \hline
    \hline
    \multirow{3}{*}{VII}     & 100  &  0.32 ( 0.04 )&  0.68 ( 0.046 )  & 
    \multirow{3}{*}{VIII}     &  0.25 ( 0.03 )&  0.5 ( 0.034 ) \\ 
      & 1000  &  0.6 ( 0.28 )&  2.02 ( 0.113 ) &  &  0.34 ( 0.05 )&  2 ( 0.002 ) \\
         & 2000  &  0.98 ( 0.44 )&  2.08 ( 0.206 ) &  &  0.54 ( 0.35 )&  2.04 ( 0.12 ) \\
    \hline
       \hline
    \multirow{3}{*}{VIV}    & 100  &  0.43 ( 0.06 )&  0.94 ( 0.053 ) \\ 
    & 1000  &  0.47 ( 0.1 )&  2 ( 0.001 ) \\         
    & 2000  &  0.58 ( 0.25 )&  2.03 ( 0.101 ) \\ 
    \cline{1-4}
  \end{tabular}
\end{table}

\begin{table}
\centering
  \caption{Compare Lasso-SIR and Sparse SIR: $\bSigma=\bSigma_1$, $\rho=0.8$, and $n=1,000$. \label{tab:SDR:4}  }
   \begin{tabular}{|c|c|c|c||c|c|c|}
    \hline
  Setting  &  p & Lasso-SIR  & Sparse SIR &Setting & Lasso-SIR & Sparse SIR\\
    \hline
    \hline
    \multirow{3}{*}{I}       & 100  &  0.18 ( 0.02 )&  0.34 ( 0.032 ) &
    \multirow{3}{*}{II}   &  0.1 ( 0.01 )&  0.13 ( 0.013 ) \\  
     & 1000  &  0.24 ( 0.02 )&  1.41 ( 0.001 ) &  &  0.16 ( 0.01 )&  1.41 ( 0.002 ) \\   
       & 2000  &  0.27 ( 0.03 )&  1.41 ( 0 ) &  &  0.19 ( 0.02 )&  1.41 ( 0 ) \\  
    \hline
    \hline
    \multirow{3}{*}{III}    & 100  &  0.28 ( 0.04 )&  0.59 ( 0.048 ) &
    \multirow{3}{*}{IV}    &  0.74 ( 0.07 )&  0.95 ( 0.069 ) \\  
      & 1000  &  0.45 ( 0.08 )&  1.41 ( 0.001 ) &  &  0.75 ( 0.07 )&  1.41 ( 0 ) \\  
        & 2000  &  0.54 ( 0.11 )&  1.41 ( 0 ) &  &  0.79 ( 0.16 )&  1.44 ( 0.102 ) \\   
    \hline
    \hline
    \multirow{3}{*}{V}    & 100  &  0.19 ( 0.04 )&  0.44 ( 0.04 ) &
    \multirow{3}{*}{VI}      &  0.52 ( 0.12 )&  1.1 ( 0.049 ) \\ 
     & 1000  &  0.31 ( 0.1 )&  1.42 ( 0.032 ) &  &  0.79 ( 0.11 )&  2 ( 0.027 ) \\
        & 2000  &  0.5 ( 0.34 )&  1.46 ( 0.135 )  &  &  0.96 ( 0.2 )&  2.01 ( 0.1 ) \\  
    \hline
    \hline
    \multirow{3}{*}{VII}    & 100  &  0.8 ( 0.34 )&  1.13 ( 0.138 ) &
    \multirow{3}{*}{VIII}   &  0.35 ( 0.12 )&  0.72 ( 0.093 ) \\
     & 1000  &  1.09 ( 0.2 )&  1.82 ( 0.136 ) & &  1 ( 0.41 )&  2.06 ( 0.176 ) \\ 
      & 2000  &  1.13 ( 0.19 )&  1.83 ( 0.165 ) & &  1.14 ( 0.36 )&  2.02 ( 0.175 ) \\  
    \hline
       \hline
    \multirow{3}{*}{VIV}   & 100  &  0.78 ( 0.1 )&  1.46 ( 0.052 ) \\           & 1000  &  0.91 ( 0.12 )&  2.01 ( 0.04 ) \\           & 2000  &  0.97 ( 0.14 )&  2.02 ( 0.061 ) \\     
    \cline{1-4}
  \end{tabular}
\end{table}

\begin{table}
\centering
  \caption{Compare Lasso-SIR and Sparse SIR: $\bSigma=\bSigma_2$, $\rho=0.2$, and $n=1,000$. \label{tab:SDR:5}  }
   \begin{tabular}{|c|c|c|c||c|c|c|}
    \hline
  Setting  &  p & Lasso-SIR  & Sparse SIR &Setting & Lasso-SIR & Sparse SIR\\
    \hline
    \hline
    \multirow{3}{*}{I}       & 100  &  0.13 ( 0.03 )&  0.18 ( 0.014 ) & 
    \multirow{3}{*}{II}  &  0.1 ( 0.02 )&  0.07 ( 0.006 ) \\   
     & 1000  &  0.33 ( 0.25 )&  1.45 ( 0.1 ) &  &  0.25 ( 0.06 )&  1.41 ( 0.001 ) \\
       & 2000  &  0.3 ( 0.16 )&  1.43 ( 0.063 ) & &  0.3 ( 0.1 )&  1.41 ( 0 ) \\ 
    \hline
    \hline
    \multirow{3}{*}{III}     & 100  &  0.24 ( 0.12 )&  0.39 ( 0.102 ) &
    \multirow{3}{*}{IV}    &  0.54 ( 0.05 )&  0.64 ( 0.05 ) \\  
     & 1000  &  0.55 ( 0.33 )&  1.5 ( 0.148 )&  &  0.63 ( 0.05 )&  1.41 ( 0 ) \\ 
    & 2000  &  0.59 ( 0.33 )&  1.41 ( 0.01 ) & & 0.64 ( 0.05 )&  1.41 ( 0 ) \\ 
    \hline
    \hline
    \multirow{3}{*}{V}   & 100  &  0.23 ( 0.29 )&  0.33 ( 0.259 ) &
    \multirow{3}{*}{VI}     &  0.27 ( 0.21 )&  0.56 ( 0.129 ) \\
     & 1000  &  1.03 ( 0.24 )&  1.61 ( 0.154 ) &&  1.01 ( 0.01 )&  1.73 ( 0.001 ) \\ 
        & 2000  &  1.09 ( 0.2 )&  1.61 ( 0.156 )& &  1.01 ( 0.01 )&  1.73 ( 0 ) \\
    \hline
    \hline
    \multirow{3}{*}{VII}     & 100  &  0.39 ( 0.24 )&  0.69 ( 0.149 ) &
    \multirow{3}{*}{VIII}  &  0.24 ( 0.03 )&  0.47 ( 0.027 ) \\
      & 1000  &  1.03 ( 0.06 )&  1.73 ( 0.027 ) &  &  0.97 ( 0.21 )&  1.76 ( 0.081 )\\
        & 2000  &  1.04 ( 0.02 )&  1.74 ( 0.071 ) & & 1.03 ( 0.08 )&  1.74 ( 0.062 ) \\      
    \hline
       \hline
    \multirow{3}{*}{VIV}  & 100  &  1 ( 0.12 )&  1.07 ( 0.064 ) \\           & 1000  &  1.03 ( 0.01 )&  1.73 ( 0.001 ) \\           & 2000  &  1.03 ( 0.01 )&  1.73 ( 0 ) \\          
    \cline{1-4}
  \end{tabular}
\end{table}

\end{appendices}

\end{document}